\documentclass[11pt,twoside]{journal}

\newcommand{\myheader}[4]{ \vspace*{-25mm}
	\noindent 
	\parbox[t]{\textwidth}{
	\parbox[t]{6cm}{
	\vspace*{0cm}
	\fontsize{8}{8}\selectfont
	Arch. Comput. Meth. Engng.\\
        Vol. 
   	\bf 
	#1, #2,
	\rm
	#3-#4
	(\the\year)
	}
	\hfill
	\parbox[t]{5.5cm}{\centering
	\fontsize{11}{12}\selectfont \bf      	     
      	     \hrule width \hsize depth 0.3pt height 0pt 
	     \vskip2.5pt%
	Archives of Computational \\ Methods in Engineering \\
\fontsize{8}{8}\selectfont \sl 
	State of the art reviews
      	     \vskip2.5pt%
      	     \hrule width \hsize depth 0.3pt height 0pt 
	}
	}
\vspace{1.0cm}
\\
\addtocounter{page}{#3}
\addtocounter{page}{-1}
}

\newcommand{\mytitle}[1]{\noindent\parbox[t]{\textwidth}{\fontsize{17}{17}\selectfont \bf #1}\vspace{0.5cm}}

 \newcommand{\myauthor}[2]{%
 \noindent\parbox[t]{\textwidth}{\fontsize{11}{11}\selectfont #1}\newline \parbox[t]{13cm}{\fontsize{8}{9.5}\selectfont #2}\vspace{0.4cm}}

\newcommand{\summary}[1]{\vspace{0.4cm}\noindent\parbox[t]{\textwidth}{\fontsize{12}{12}\selectfont Summary} \\[1.5ex] \noindent\parbox[t]{\textwidth}{\fontsize{9}{10}\selectfont #1} \vspace{0.5cm}}

\usepackage{latexsym}
\usepackage{amsbsy}
\usepackage{amsmath}
\usepackage{amssymb}
\usepackage{graphicx}
\usepackage{psfrag}
\usepackage{url}
\usepackage{algorithm}
\usepackage{fancyvrb}
\usepackage{color}

\newcommand{\dx}{\, \mathrm{d}x}
\newcommand{\dX}{\, \mathrm{d}X}
\newcommand{\ds}{\, \mathrm{d}s}
\newcommand{\dt}{\, \mathrm{d}t}
\newcommand{\R}{\mathbb{R}}
\newcommand{\tab}{\hspace*{2em}}

\newtheorem{theorem}{Theorem}

\DefineVerbatimEnvironment{code}{Verbatim}{frame=single,rulecolor=\color{blue}}

\begin{document}

\year2006
\received={March 2006}
\myheader{00}{0}{1}{76}
\markboth{Anders Logg}{Automating the Finite Element Method}
\thispagestyle{myfirstpage}

\mytitle{Automating the Finite Element Method}

\myauthor{Anders Logg}%
{Toyota Technological Institute at Chicago \\
 University Press Building \\
 1427 East 60th Street \\
 Chicago, Illinois 60637 USA \\ \\
 Email: \texttt{logg@tti-c.org} \\
 \\\\}

\summary{%
  The finite element method can be viewed as a machine that automates
  the discretization of differential equations, taking as input a
  variational problem, a finite element and a mesh, and producing as
  output a system of discrete equations. However, the generality of
  the framework provided by the finite element method is seldom
  reflected in implementations (realizations), which are often
  specialized and can handle only a small set of variational problems
  and finite elements (but are typically parametrized over the choice
  of mesh).

  \vspace{0.25cm}

  This paper reviews ongoing research in the direction of a complete
  automation of the finite element method. In particular, this work
  discusses algorithms for the efficient and automatic computation of
  a system of discrete equations from a given variational problem,
  finite element and mesh. It is demonstrated that by automatically
  generating and compiling efficient low-level code, it is possible to
  parametrize a finite element code over variational problem and
  finite element in addition to the mesh.}

%------------------------------------------------------------------------------
\section{INTRODUCTION}

The finite element method (Galerkin's method) has emerged as a
universal method for the solution of differential equations.  Much of
the success of the finite element method can be contributed to its
generality and simplicity, allowing a wide range of differential
equations from all areas of science to be analyzed and solved within a
common framework. Another contributing factor to the success of the
finite element method is the flexibility of formulation, allowing the
properties of the discretization to be controlled by the choice of
finite element (approximating spaces).

At the same time, the generality and flexibility of the finite element
method has for a long time prevented its automation, since any
computer code attempting to automate it must necessarily be
parametrized over the choice of variational problem and finite
element, which is difficult. Consequently, much of the work must still
be done by hand, which is both tedious and error-prone, and results in
long development times for simulation codes.

Automating systems for the solution of differential equations are often
met with skepticism, since it is believed that the generality and
flexibility of such tools cannot be combined with the efficiency of
competing specialized codes that only need to handle one equation for
a single choice of finite element. However, as will be demonstrated in
this paper, by automatically generating and compiling low-level code
for any given equation and finite element, it is possible to develop
systems that realize the generality and flexibility of the finite
element method, while competing with or outperforming specialized and
hand-optimized codes.

\subsection{Automating the Finite Element Method}

To automate the finite element method, we need to build a machine that
takes as input a discrete variational problem posed on a pair of
discrete function spaces defined by a set of finite elements on a
mesh, and generates as output a system of discrete equations for the
degrees of freedom of the solution of the variational problem. In
particular, given a discrete variational problem of the form: Find $U
\in V_h$ such that
\begin{equation} \label{eq:varproblem,first}
  a(U; v) = L(v) \quad \forall v \in \hat{V}_h,
\end{equation}
where $a : V_h \times \hat{V}_h \rightarrow \R$ is a semilinear form
which is linear in its second argument, $L : \hat{V}_h \rightarrow \R$ a linear
form and $(\hat{V}_h, V_h)$ a given pair of discrete function spaces
(the \emph{test} and \emph{trial} spaces), the machine should
automatically generate the discrete system
\begin{equation} \label{eq:system,nonlinear}
  F(U) = 0,
\end{equation}
where $F : V_h \rightarrow \R^N$, $N = |\hat{V}_h| = |V_h|$ and
\begin{equation}
  F_i(U) = a(U; \hat{\phi}_i) - L(\hat{\phi}_i), \quad i=1,2,\ldots,N,
\end{equation}
for $\{\hat{\phi}_i\}_{i=1}^N$ a given basis for $\hat{V}_h$.

Typically, the discrete variational
problem~(\ref{eq:varproblem,first}) is obtained as the discrete
version of a corresponding continuous variational problem: Find $u \in
V$ such that
\begin{equation}
  a(u; v) = L(v) \quad \forall v \in \hat{V},
\end{equation}
where $\hat{V}_h \subset \hat{V}$ and $V_h \subset V$.

The machine should also automatically generate the discrete
representation of the linearization of the given semilinear form~$a$, that
is the matrix $A \in \R^{N\times N}$ defined by
\begin{equation}
  A_{ij}(U) = a'(U; \hat{\phi}_i, \phi_j),
  \quad i,j = 1,2,\ldots,N,
\end{equation}
where $a' : V_h \times \hat{V}_h \times V_h \rightarrow \R$ is the
Fr\'echet derivative of~$a$ with respect to its first
argument and $\{\hat{\phi}_i\}_{i=1}^N$ and $\{\phi_i\}_{i=1}^N$
are bases for $\hat{V}_h$ and $V_h$ respectively.

In the simplest case of a linear variational problem,
\begin{equation}
  a(v, U) = L(v) \quad \forall v \in \hat{V}_h,
\end{equation}
the machine should automatically generate the linear system
\begin{equation} \label{eq:system,linear}
  A U = b,
\end{equation}
where $A_{ij} = a(\hat{\phi}_i, \phi_j)$ and $b_i = L(\hat{\phi}_i)$, and
where $(U_i) \in \R^N$ is the vector of degrees of freedom for the
discrete solution $U$, that is, the expansion coefficients in
the given basis for $V_h$,
\begin{equation}
  U = \sum_{i=1}^N U_i \phi_i.
\end{equation}

We return to this in detail below and identify the key steps towards a
complete automation of the finite element method, including algorithms
and prototype implementations for each of the key steps.

\subsection{The FEniCS Project and the Automation of CMM}

The FEniCS project~\cite{logg:www:03,logg:preprint:10} was initiated
in~2003 with the explicit goal of developing free software for the
Automation of Computational Mathematical Modeling (CMM), including a
complete automation of the finite element method. As such, FEniCS
serves as a prototype implementation of the methods and principles put
forward in this paper.

In~\cite{logg:thesis:03}, an agenda for the automation of CMM is
outlined, including the automation of
(i)~discretization,
(ii)~discrete solution,
(iii)~error control,
(iv)~modeling and
(v)~optimization.
The automation of discretization amounts to the automatic
generation of the system of discrete
equations~(\ref{eq:system,nonlinear}) or~(\ref{eq:system,linear}) from
a given given differential equation or variational problem. Choosing
as the foundation for the automation of discretization the finite
element method, the first step towards the Automation of CMM is thus
the automation of the finite element method. We continue the
discussion on the automation of CMM below in Section~\ref{sec:outlook}.

Since the initiation of the FEniCS project in 2003, much progress has
been made, especially concerning the automation of discretization.  In
particular, two central components that automate central aspects of
the finite element method have been developed. The first of these
components is FIAT, the FInite element Automatic
Tabulator~\cite{www:FIAT,Kir04,Kir06}, which automates the generation
of finite element basis functions for a large class of finite
elements. The second component is FFC, the FEniCS Form
Compiler~\cite{logg:www:04,logg:article:10,logg:submit:toms-ffc-monomial-2006,logg:manual:02},
which automates the evaluation of variational problems by
automatically generating low-level code for the assembly of the system
of discrete equations from given input in the form of a variational
problem and a (set of) finite element(s).

In addition to FIAT and FFC, the FEniCS project develops components
that wrap the functionality of collections of other FEniCS components
(middleware) to provide simple, consistent and intuitive user
interfaces for application programmers. One such example is
DOLFIN~\cite{logg:www:01,logg:preprint:06,logg:manual:01}, which
provides both a C++ and a Python interface (through
SWIG~\cite{www:SWIG,Bea96}) to the basic functionality of FEniCS.

We give more details below in Section~\ref{sec:prototype} on FIAT,
FFC, DOLFIN and other FEniCS components, but list here some of the key
properties of the software components developed as part of the FEniCS
project, as well as the FEniCS system as a whole:
\begin{itemize}
\item
  automatic and efficient evaluation of variational problems through
  FFC~\cite{logg:www:04,logg:article:10,logg:submit:toms-ffc-monomial-2006,logg:manual:02}, including support
  for arbitrary mixed formulations;
\item
  automatic and efficient assembly of systems of discrete equations
  through DOLFIN~\cite{logg:www:01,logg:preprint:06,logg:manual:01};
\item
  support for general families of finite elements, including
  continuous and discontinuous Lagrange finite elements of arbitrary
  degree on simplices through FIAT~\cite{www:FIAT,Kir04,Kir06};
\item
  high-performance parallel linear algebra through
  PETSc~\cite{www:PETSc,BalBus04,BalEij97};
\item
  arbitrary order multi-adaptive
  $\mathrm{mcG}(q)$/$\mathrm{mdG}(q)$ and mono-adaptive
  $\mathrm{cG}(q)$/$\mathrm{dG}(q)$ ODE
  solvers~\cite{logg:article:01,logg:article:02,logg:article:03,logg:article:05,logg:article:08}.
\end{itemize}

\subsection{Automation and Mathematics Education}

By automating mathematical concepts, that is, implementing
corresponding concepts in software, it is possible to raise the level
of the mathematics education. An aspect of this is the possibility of
allowing students to experiment with mathematical concepts and thereby
obtaining an increased understanding (or familiarity) for the
concepts. An illustrative example is the concept of a vector in
$\R^n$, which many students get to know very well through
experimentation and exercises in Octave~\cite{www:Octave} or
MATLAB~\cite{www:MATLAB}. If asked which is the true vector, the~$x$
on the blackboard or the~\texttt{x} on the computer screen, many
students (and the author) would point towards the computer.

By automating the finite element method, much like linear algebra has
been automated before, new advances can be brought to the mathematics
education. One example of this is
Puffin~\cite{logg:www:05,logg:manual:03}, which is a minimal and
educational implementation of the basic functionality of FEniCS for
Octave/MATLAB. Puffin has successfully been used in a number of
courses at Chalmers in G\"oteborg and the Royal Institute of
Technology in Stockholm, ranging from introductory undergraduate
courses to advanced undergraduate/beginning graduate courses. Using
Puffin, first-year undergraduate students are able to design and
implement solvers for coupled systems of convection--diffusion--reaction
equations, and thus obtaining important understanding of mathematical
modeling, differential equations, the finite element method and
programming, without reducing the students to button-pushers.

Using the computer as an integrated part of the mathematics education
constitutes a change of paradigm~\cite{logg:book:02}, which will have
profound influence on future mathematics education.

\subsection{Outline}

This paper is organized as follows. In the next section, we first
present a survey of existing finite element software that automate
particular aspects of the finite element method. In
Section~\ref{sec:fem}, we then give an introduction to the finite
element method with special emphasis on the process of generating the
system of discrete equations from a given variational problem, finite
element(s) and mesh. A summary of the notation can be found at the end
of this paper.

Having thus set the stage for our main topic, we next identify in
Sections~\ref{sec:automation,tabulation}--\ref{sec:automation,assembly}
the key steps towards an automation of the finite element method and
present algorithms and systems that accomplish (in part) the
automation of each of these key steps. We also discuss a framework for
generating an optimized computation from these algorithms in
Section~\ref{sec:optimizations}. In Section~\ref{sec:se}, we then
highlight a number of important concepts and techniques from software
engineering that play an important role for the automation of the
finite element method.

Prototype implementations of the algorithms are then discussed in
Section~\ref{sec:prototype}, including benchmark results that
demonstrate the efficiency of the algorithms and their
implementations. We then, in Section~\ref{sec:examples}, present a
number of examples to illustrate the benefits of a system automating
the finite element method. As an outlook towards further research, we
present in Section~\ref{sec:outlook} an agenda for the development of
an extended automating system for the Automation of CMM, for which the
automation of the finite element method plays a central role. Finally,
we summarize our findings in Section~\ref{sec:conclusion}.

%------------------------------------------------------------------------------
\section{SURVEY OF CURRENT FINITE ELEMENT SOFTWARE}

\label{sec:survey}

There exist today a number of projects that seek to create systems
that (in part) automate the finite element method. In this section, we
survey some of these projects. A complete survey is difficult to make
because of the large number of such projects. The survey is instead
limited to a small set of projects that have attracted the attention
of the author. In effect, this means that most proprietary systems
have been excluded from this survey.

It is instructional to group the systems both by their level of
automation and their design. In particular, a number of systems
provide automated generation of the system of discrete equations from
a given variational problem, which we in this paper refer to as the
automation of the finite element method or automatic assembly, while
other systems only provide the user with a basic toolkit for this
purpose.  Grouping the systems by their design, we shall differentiate
between systems that provide their functionality in the form of a
library in an existing language and systems that implement new
domain-specific languages for finite element computation. A summary
for the surveyed systems is given in Table~\ref{tab:projects}.

It is also instructional to compare the basic specification of a
simple test problem such as Poisson's equation, $-\Delta u = f$ in
some domain $\Omega \subset \R^d$, for the surveyed systems, or more
precisely, the specification of the corresponding discrete variational
problem $a(v, U) = L(v)$ for all $v$ in some suitable test space, with
the bilinear form~$a$ given by
\begin{equation} \label{eq:poisson,a}
  a(v, U) = \int_{\Omega} \nabla v \cdot \nabla U \dx,
\end{equation}
and the linear form $L$ given by
\begin{equation} \label{eq:poisson,L}
  L(v) = \int_{\Omega} v \, f \dx.
\end{equation}

Each of the surveyed systems allow the specification of the
variational problem for Poisson's equation with varying degree of
automation. Some of the systems provide a high level of automation and
allow the variational problem to be specified in a notation that is
very close to the mathematical notation used in~(\ref{eq:poisson,a})
and~(\ref{eq:poisson,L}), while others require more user-intervention.
In connection to the presentation of each of the surveyed systems
below, we include as an illustration the specification of the
variational problem for Poisson's equation in the notation employed by
the system in question. In all cases, we include only the part of the
code essential to the specification of the variational problem. Since
the different systems are implemented in different languages,
sometimes even providing new domain-specific languages, and since
there are differences in philosophies, basic concepts and
capabilities, it is difficult to make a uniform presentation. As a
result, not all the examples specify exactly the same problem.

\begin{table}[htbp]
  \begin{center}
    \begin{tabular}{|l||c|c|c|}
      \hline
      Project & automatic assembly & library / language & license \\
      \hline
      \hline
      Analysa & yes & language & proprietary \\
      \hline
      deal.II & no & library & QPL\footnotemark{} \\
      \hline
      Diffpack & no & library & proprietary \\
      \hline
      FEniCS & yes & both & GPL, LGPL \\
      \hline
      FreeFEM & yes & language & LGPL \\
      \hline
      GetDP & yes & language & GPL \\
      \hline
      Sundance & yes & library & LGPL \\
      \hline
    \end{tabular}
    \caption{Summary of projects seeking to automate the finite element method.}
    \label{tab:projects}
  \end{center}
\end{table}

\footnotetext{In addition to the terms imposed by the QPL, the deal.II
  license imposes a form of advertising clause, requiring the citation
  of certain publications. See \cite{www:deal.II} for details.}

\subsection{Analysa}
\label{sec:analysa}

Analysa~\cite{www:Analysa,BagSco04} is a domain-specific language and
problem-solving environment (PSE) for partial differential
equations. Analysa is based on the functional language Scheme and
provides a language for the definition of variational problems.
Analysa thus falls into the category of domain-specific languages.

Analysa puts forward the idea that it is sometimes desirable to
compute the \emph{action} of a bilinear form, rather than assembling the
matrix representing the bilinear form in the current basis.  In the
notation of~\cite{BagSco04}, the action of a bilinear form~$a :
\hat{V}_h \times V_h \rightarrow \R$ on a given discrete function $U
\in V_h$ is
\begin{equation}
  w = a(\hat{V}_h, U) \in \R^N,
\end{equation}
where
\begin{equation}
  w_i = a(\hat{\phi}_i, U), \quad i = 1, 2, \ldots, N.
\end{equation}
Of course, we have $w = A U$, where $A$ is the matrix representing the
bilinear form, with $A_{ij} = a(\hat{\phi}_i, \phi_j)$, and
$(U_i) \in \R^N$ is the vector of expansion coefficients for $U$ in the
basis of $V_h$. It follows that
\begin{equation}
  w = a(\hat{V}_h, U) = a(\hat{V}_h, V_h) U.
\end{equation}

If the action only needs to be evaluated a few times for different
discrete functions~$U$ before updating a linearization (reassembling the
matrix~$A$), it might be more efficient to compute each action
directly than first assembling the matrix~$A$ and applying it to
each~$U$.

To specify the variational problem for Poisson's equation with
Analysa, one specifies a pair of bilinear forms~\texttt{a}
and~\texttt{m}, where~\texttt{a} represents the bilinear form~$a$
in~(\ref{eq:poisson,a}) and~\texttt{m} represents the bilinear
form
\begin{equation}
  m(v, U) = \int_{\Omega} v \, U \dx,
\end{equation}
corresponding to a mass matrix. In the language of Analysa, the linear form~$L$
in~(\ref{eq:poisson,L}) is represented as the application of the bilinear
form~$m$ on the test space~$\hat{V}_h$ and the right-hand side~$f$,
\begin{equation}
  L(\hat{\phi}_i) = m(\hat{V}_h, f)_i, \quad i=1,2,\ldots,N,
\end{equation}
as shown in Table~\ref{tab:poisson,analysa}. Note that Analysa thus
defers the coupling of the forms and the test and trial spaces until
the computation of the system of discrete equations.

\begin{table}[htbp]
  \begin{code}
 (integral-forms
   ((a v U) (dot (gradient v) (gradient U)))
   ((m v U) (* v U))
 )

 (elements
   (element (lagrange-simplex 1))
 )

 (spaces
   (test-space (fe element (all mesh) r:))
   (trial-space (fe element (all mesh) r:))
 )

 (functions
   (f (interpolant test-space (...)))
 )

 (define A-matrix (a testspace trial-space))
 (define b-vector (m testspace f))
  \end{code}
  \vspace{-0.75cm}
  \caption{Specifying the variational problem for Poisson's equation
  with Analysa using piecewise linear elements on simplices (triangles
  or tetrahedra).}
  \label{tab:poisson,analysa}
\end{table}

\subsection{deal.II}
\label{sec:deal}

deal.II~\cite{www:deal.II,BanKan99,Ban00} is a C++ library for finite
element computation. While providing tools for finite elements, meshes
and linear algebra, deal.II does not provide support for automatic
assembly. Instead, a user needs to supply the complete code for the
assembly of the system~(\ref{eq:system,linear}), including the
explicit computation of the element stiffness matrix (see Section
\ref{sec:fem} below) by quadrature, and the insertion of each element
stiffness matrix into the global matrix, as illustrated in
Table~\ref{tab:poisson,deal}. This is a common design for many finite
element libraries, where the ambition is not to automate the finite
element method, but only to provide a set of basic tools.

\begin{table}[htbp]
  \begin{code}
 ...
 for (dof_handler.begin_active(); cell! = dof_handler.end(); ++cell)
 {
   ...
   for (unsigned int i = 0; i < dofs_per_cell; ++i)
     for (unsigned int j = 0; j < dofs_per_cell; ++j)
       for (unsigned int q_point = 0; q_point < n_q_points; ++q_point)
         cell_matrix(i, j) += (fe_values.shape_grad (i, q_point) *
                               fe_values.shape_grad (j, q_point) *
                               fe_values.JxW(q_point));

   for (unsigned int i = 0; i < dofs_per_cell; ++i)
     for (unsigned int q_point = 0; q_point < n_q_points; ++q_point)
       cell_rhs(i) += (fe_values.shape_value (i, q_point) *
                       <value of right-hand side f> *
                       fe_values.JxW(q_point));

   cell->get_dof_indices(local_dof_indices);

   for (unsigned int i = 0; i < dofs_per_cell; ++i)
     for (unsigned int j = 0; j < dofs_per_cell; ++j)
       system_matrix.add(local_dof_indices[i],
                         local_dof_indices[j],
                         cell_matrix(i, j));

   for (unsigned int i = 0; i < dofs_per_cell; ++i)
     system_rhs(local_dof_indices[i]) += cell_rhs(i);
 }
 ...
  \end{code}
  \vspace{-0.75cm}
  \caption{Assembling the linear system~(\ref{eq:system,linear}) for
    Poisson's equation with deal.II.}
  \label{tab:poisson,deal}
\end{table}

\subsection{Diffpack}
\label{sec:diffpack}

Diffpack~\cite{www:diffpack,Lan99} is a C++ library for finite element
and finite difference solution of partial differential
equations. Initiated in 1991, in a time when most finite element codes
were written in FORTRAN, Diffpack was one of the pioneering
libraries for scientific computing with C++. Although originally
released as free software, Diffpack is now a proprietary product.

Much like deal.II, Diffpack requires the user to supply the code for
the computation of the element stiffness matrix, but automatically
handles the loop over quadrature points and the insertion of the
element stiffness matrix into the global matrix, as illustrated in
Table~\ref{tab:poisson,diffpack}.

\begin{table}[htbp]
  \begin{code}
 for (int i = 1; i <= nbf; i++)
   for (int j = 1; j <= nbf; j++)
     elmat.A(i, j) += (fe.dN(i, 1) * fe.dN(j, 1) +
                       fe.dN(i, 2) * fe.dN(j, 2) +
                       fe.dN(i, 3) * fe.dN(j, 3)) * detJxW;

 for (int i = 1; i <= nbf; i++)
   elmat.b(i) += fe.N(i)*<value of right-hand side f>*detJxW;
  \end{code}
  \vspace{-0.75cm}
  \caption{Computing the element stiffness matrix and element
    load vector for Poisson's equation with Diffpack.}
  \label{tab:poisson,diffpack}
\end{table}

\subsection{FEniCS}

The FEniCS project~\cite{logg:www:03,logg:preprint:10} is structured
as a system of interoperable components that automate central aspects
of the finite element method. One of these components is the form
compiler
FFC~\cite{logg:www:04,logg:article:10,logg:submit:toms-ffc-monomial-2006,logg:manual:02},
which takes as input a variational problem together with a set of
finite elements and generates low-level code for the automatic
computation of the system of discrete equations. In this regard, the
FEniCS system implements a domain-specific language for finite element
computation, since the form is entered in a special language
interpreted by the compiler. On the other hand, the form compiler FFC
is also available as a Python module and can be used as a just-in-time
(JIT) compiler, allowing variational problems to be specified and
computed with from within the Python scripting environment. The FEniCS
system thus falls into both categories of being a library and a
domain-specific language, depending on which interface is used.

To specify the variational problem for Poisson's equation with
FEniCS, one must specify a pair of basis functions \texttt{v} and
\texttt{U}, the right-hand side function \texttt{f}, and of course the
bilinear form~\texttt{a} and the linear form~\texttt{L}, as shown in
Table~\ref{tab:poisson,fenics}.

\begin{table}[htbp]
  \begin{code}
 element = FiniteElement(``Lagrange'', ``tetrahedron'', 1)

 v = BasisFunction(element)
 U = BasisFunction(element)
 f = Function(element)

 a = dot(grad(v), grad(U))*dx
 L = v*f*dx
  \end{code}
  \vspace{-0.75cm}
  \caption{Specifying the variational problem for Poisson's equation
  with FEniCS using piecewise linear elements on tetrahedra.}
  \label{tab:poisson,fenics}
\end{table}

Note in Table~\ref{tab:poisson,fenics} that the function spaces
(finite elements) for the test and trial functions~\texttt{v}
and~\texttt{U} together with all additional functions/coefficients (in
this case the right-hand side~\texttt{f}) are fixed at compile-time,
which allows the generation of very efficient low-level code since the
code can be generated for the specific given variational problem and
the specific given finite element(s).

Just like Analysa, FEniCS (or FFC) supports the specification of
actions, but while Analysa allows the specification of a general
expression that can later be treated as a bilinear form, by applying
it to a pair of function spaces, or as a linear form, by applying it
to a function space and a given fixed function, the arity of the form
must be known at the time of specification in the form language of
FFC. As an example, the specification of a linear form~\texttt{a}
representing the action of the bilinear form~(\ref{eq:poisson,a}) on a
function~\texttt{U} is given in Table~\ref{tab:poisson,fenics,action}.

\begin{table}[htbp]
  \begin{code}
 element = FiniteElement(``Lagrange'', ``tetrahedron'', 1)

 v = BasisFunction(element)
 U = Function(element)

 a = dot(grad(v), grad(U))*dx
  \end{code}
  \vspace{-0.75cm}
  \caption{Specifying the linear form for the action of the bilinear
    form~(\ref{eq:poisson,a}) with FEniCS using piecewise linear
    elements on tetrahedra.}
  \label{tab:poisson,fenics,action}
\end{table}

A more detailed account of the various components of the FEniCS
project is given below in Section~\ref{sec:prototype}.

\subsection{FreeFEM}

FreeFEM \cite{www:FreeFEM,HecPir05} implements a domain-specific
language for finite element solution of partial differential
equations. The language is based on C++, extended with a special
language that allows the specification of variational problems. In
this respect, FreeFEM is a compiler, but it also provides an
integrated development environment (IDE) in which programs can be
entered, compiled (with a special compiler) and
executed. Visualization of solutions is also provided.

FreeFEM comes in two flavors, the current version FreeFEM++ which
only supports 2D problems and the 3D version FreeFEM3D. Support for 3D
problems will be added to FreeFEM++ in the future.~\cite{www:FreeFEM}.

To specify the variational problem for Poisson's equation with
FreeFEM++, one must first define the test and trial spaces (which we
here take to be the same space~\texttt{V}), and then the test and
trial functions~\texttt{v} and~\texttt{U}, as well as the
function~\texttt{f} for the right-hand side. One may then define the
bilinear form~\texttt{a} and linear form~\texttt{L} as illustrated in
Table~\ref{tab:poisson,freefem}.

\begin{table}[htbp]
  \begin{code}
 fespace V(mesh, P1);
 V v, U;
 func f = ...;

 varform a(v, U) = int2d(mesh)(dx(v)*dx(U) + dy(v)*dy(U));
 varform L(v) = int2d(mesh)(v*f);
  \end{code}
  \vspace{-0.75cm}
  \caption{Specifying the variational problem for Poisson's equation
    with FreeFEM++ using piecewise linear elements on triangles (as
    determined by the mesh).}
  \label{tab:poisson,freefem}
\end{table}

\subsection{GetDP}

GetDP~\cite{www:GetDP,DulGeu05} is a finite element solver which
provides a special declarative language for the specification of
variational problems. Unlike FreeFEM, GetDP is not a compiler, nor is
it a library, but it will be classified here under the category of
domain-specific languages. At start-up, GetDP parses a problem
specification from a given input file and then proceeds according to
the specification.

To specify the variational problem for Poisson's equation with
GetDP, one must first give a definition of a function space, which may
include constraints and definition of sub spaces. A variational
problem may then be specified in terms of functions from the
previously defined function spaces, as illustrated in
Table~\ref{tab:poisson,getdp}.

\begin{table}[htbp]
  \begin{code}
 FunctionSpace {
   { Name V; Type Form0;
     BasisFunction {
       { ... }
     }
   }
 }

 Formulation {
   { Name Poisson; Type FemEquation;
     Quantity {
       { Name v; Type Local; NameOfSpace V; }
     }
     Equation {
       Galerkin { [Dof{Grad v}, {Grad v}];
	          ....
       }
     }
   }
 }
  \end{code}
  \vspace{-0.75cm}
  \caption{Specifying the bilinear form for Poisson's equation
    with GetDP.}
  \label{tab:poisson,getdp}
\end{table}

\subsection{Sundance}

Sundance~\cite{www:Sundance,Lon03,Lon04} is a C++ library for finite element
solution of partial differential equations (PDEs), with special
emphasis on large-scale PDE-constrained optimization.

Sundance supports automatic generation of the system of discrete
equations from a given variational problem and has a powerful symbolic
engine, which allows variational problems to be specified and
differentiated symbolically natively in C++. Sundance thus falls into
the category of systems providing their functionality in the form of
library.

To specify the variational problem for Poisson's equation with
Sundance, one must specify a test function \texttt{v}, an unknown
function \texttt{U}, the right-hand side \texttt{f}, the differential
operator \texttt{grad} and the variational problem written in the form
$a(v, U) - L(v) = 0$, as shown in Table~\ref{tab:poisson,sundance}.

\begin{table}[htbp]
  \begin{code}
 Expr v = new TestFunction(new Lagrange(1));
 Expr U = new UnknownFunction(new Lagrange(1));
 Expr f = new DiscreteFunction(...);

 Expr dx = new Derivative(0);
 Expr dy = new Derivative(1);
 Expr dz = new Derivative(2);
 Expr grad = List(dx, dy, dz);

 Expr poisson = Integral((grad*v)*(grad*U) - v*f);
  \end{code}
  \vspace{-0.75cm}
  \caption{Specifying the variational problem for Poisson's equation
    with Sundance using piecewise linear elements on tetrahedra (as
    determined by the mesh).}
  \label{tab:poisson,sundance}
\end{table}

%------------------------------------------------------------------------------
\section{THE FINITE ELEMENT METHOD}
\label{sec:fem}

\begin{center}
  \begin{minipage}{10cm}
  \emph{It once happened that a man thought he had written original
    verses, and was then found to have read them word for word, long
    before, in some ancient poet.}

  \vspace{0.5cm}

  \raggedleft
  Gottfried Wilhelm Leibniz \\
  \emph{Nouveaux essais sur l'entendement humain} (1704/1764)
  \end{minipage}
\end{center}

In this section, we give an overview of the finite element method,
with special focus on the general algorithmic aspects that form the
basis for its automation. In many ways, the material is
standard~\cite{ZieTay67,StrFix73,Cia76,Cia78,BecCar81,Hug87,BreSco94,EriEst96,SamWib98},
but it is presented here to give a background for the continued
discussion on the automation of the finite element method and to
summarize the notation used throughout the remainder of this
paper. The purpose is also to make precise what we set out to
automate, including assumptions and limitations.

\subsection{Galerkin's Method}

Galerkin's method (the weighted residual method) was originally
formulated with global polynomial spaces~\cite{Gal15} and goes back
to the variational principles of Leibniz, Euler, Lagrange, Dirichlet,
Hamilton, Castigliano \cite{Cas1879}, Rayleigh~\cite{Ray1870} and Ritz
\cite{Rit08}. Galerkin's method with piecewise polynomial spaces
$(\hat{V}_h,V_h)$ is known as the \emph{finite element method}. The finite
element method was introduced by engineers for structural analysis in
the 1950s and was independently proposed by Courant in 1943
\cite{Cou43}. The exploitation of the finite element method among
engineers and mathematicians exploded in the 1960s. In addition to the
references listed above, we point to the following general
references:~\cite{EriEst95,EriJoh91,EriJoh95a,EriJohIII,EriJoh95b,EriJoh95c,EriJoh98,BecRan01}.

We shall refer to the family of Galerkin methods (weighted residual
methods) with piecewise (polynomial) function spaces as the finite
element method, including Petrov-Galerkin methods (with different test
and trial spaces) and Galerkin/least-squares methods.

\subsection{Finite Element Function Spaces}
\label{sec:femfunctions}

A central aspect of the finite element method is the construction of
discrete function spaces by piecing together local function spaces on
the cells $\{K\}_{K\in\mathcal{T}}$ of a mesh~$\mathcal{T}$ of a domain
$\Omega = \cup_{K \in \mathcal{T}} \subset \R^d$, with each local
function space defined by a \emph{finite element}.

\subsubsection{The finite element}

We shall use the standard Ciarlet~\cite{Cia78,BreSco94}
definition of a finite element, which reads as follows.
A finite element is a triple $(K, \mathcal{P}_K, \mathcal{N}_K)$, where
\begin{itemize}
\item
  $K \subset \R^d$ is a bounded closed subset of $\R^d$ with nonempty
  interior and piecewise smooth boundary;
\item
  $\mathcal{P}_K$ is a function space on $K$ of dimension $n_K < \infty$;
\item
  $\mathcal{N}_K = \{\nu^K_1, \nu^K_2, \ldots, \nu^K_{n_K}\}$ is a basis for
  $\mathcal{P}_K'$ (the bounded linear functionals on $\mathcal{P}_K$).
\end{itemize}
We shall further assume that we are given a nodal basis $\{\phi^K_i\}_{i=1}^{n_K}$
for~$\mathcal{P}_K$ that for each \emph{node} $\nu^K_i \in \mathcal{N}_K$
satisfies $\nu^K_i(\phi^K_j) = \delta_{ij}$
for $j = 1,2,\ldots,n_K$. Note that this implies that for any $v
\in \mathcal{P}_K$, we have
\begin{equation}
  v = \sum_{i=1}^{n_K} \nu^K_i(v) \phi^K_i.
\end{equation}
In the simplest case, the nodes are given by evaluation of function
values or directional derivatives at a set of points
$\{x^K_i\}_{i=1}^{n_K}$, that is,
\begin{equation}
  \nu^K_i(v) = v(x^K_i), \quad i = 1,2,\ldots,n_K.
\end{equation}

\subsubsection{The local-to-global mapping}

Now, to define a global function space $V_h = \mathrm{span}
\{\phi_i\}_{i=1}^N$ on $\Omega$ and a set of global nodes $\mathcal{N}
= \{\nu_i\}_{i=1}^N$ from a given set $\{(K,
\mathcal{P}_K,\mathcal{N}_K)\}_{K\in\mathcal{T}}$ of finite elements,
we also need to specify how the local function spaces are pieced
together. We do this by specifying for each cell $K \in \mathcal{T}$ a
\emph{local-to-global mapping},
\begin{equation}
  \iota_K : [1,n_K] \rightarrow N,
\end{equation}
that specifies how the local nodes~$\mathcal{N}_K$ are mapped to global
nodes~$\mathcal{N}$, or more precisely,
\begin{equation} \label{eq:nodemapping}
  \nu_{\iota_K(i)}(v) = \nu^K_i(v|_K),
  \quad i = 1,2,\ldots,n_K,
\end{equation}
for any $v\in V_h$,
that is, each local node~$\nu^K_i \in \mathcal{N}_K$ corresponds to a global
node~$\nu_{\iota_K(i)} \in \mathcal{N}$ determined by the local-to-global mapping
$\iota_K$.

\subsubsection{The global function space}

We now define the global function space $V_h$ as the set of functions
on $\Omega$ satisfying
\begin{equation}
  v |_K \in \mathcal{P}_K \quad \forall K \in \mathcal{T},
\end{equation}
and furthermore satisfying the constraint that if for any
pair of cells $(K, K') \in \mathcal{T} \times \mathcal{T}$ and
local node numbers~$(i, i') \in [1,n_K] \times
[1,n_{K'}]$, we have
\begin{equation}
  \iota_K(i) = \iota_{K'}(i'),
\end{equation}
then
\begin{equation} \label{eq:constraint}
  \nu^K_i(v|_K) = \nu^{K'}_{i'}(v|_{K'}),
\end{equation}
where $v|_K$ denotes the continuous extension to $K$ of the
restriction of $v$ to the interior of $K$, that is, if two local nodes
$\nu^K_i$ and $\nu^{K'}_{i'}$ are mapped to the same global node, then
they must agree for each function $v \in V_h$.

Note that by this construction, the functions of $V_h$ are undefined on
cell boundaries, unless the constraints (\ref{eq:constraint}) force
the (restrictions of) functions of $V_h$ to be continuous on cell
boundaries, in which case we may uniquely define the functions of $V_h$
on the entire domain~$\Omega$. However, this is usually not a problem,
since we can perform all operations on the restrictions of functions
to the local cells.

\subsubsection{Lagrange finite elements}

\begin{figure}[htbp]
  \begin{center}
    \includegraphics[width=7cm]{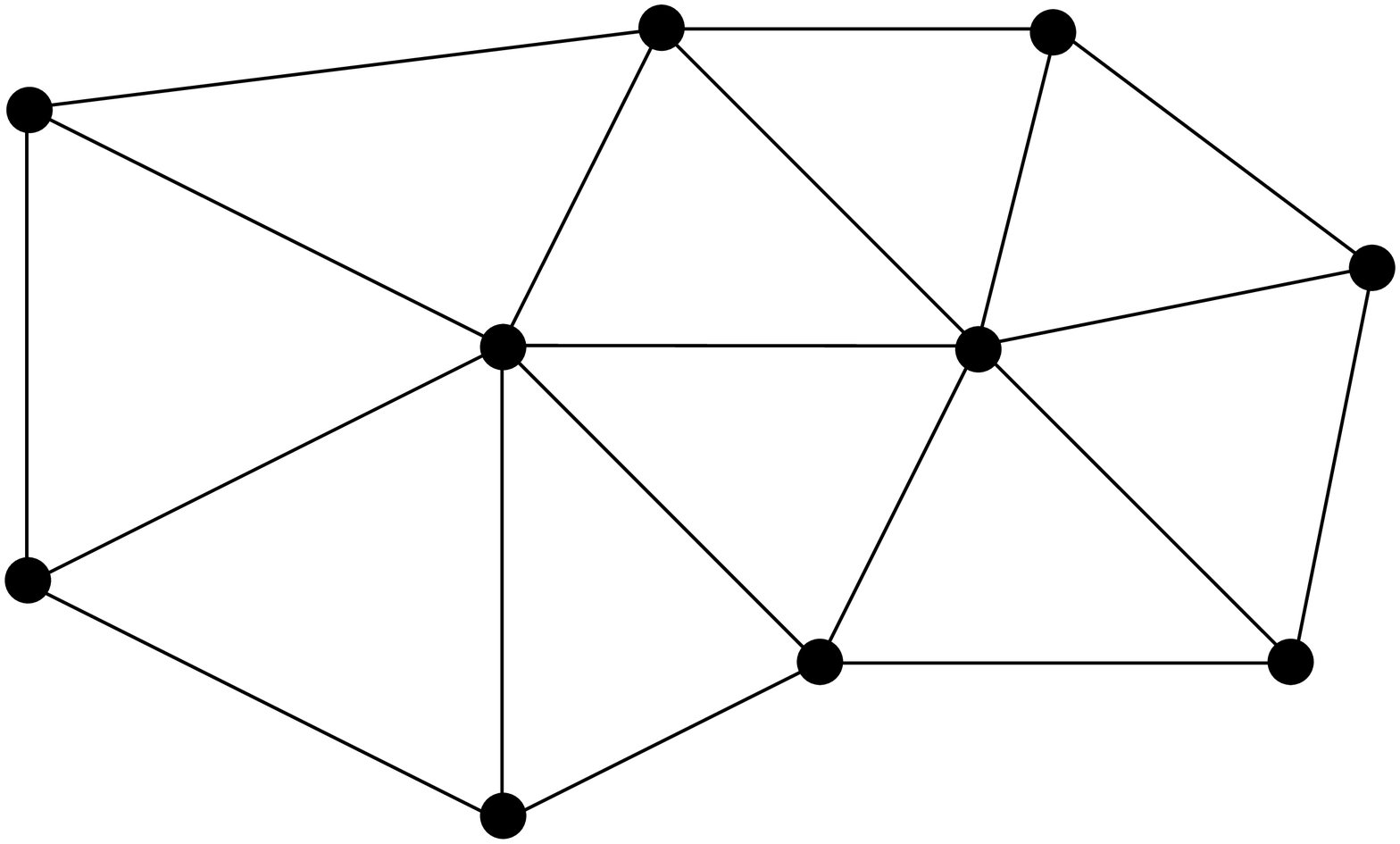}\quad\quad
    \includegraphics[width=7cm]{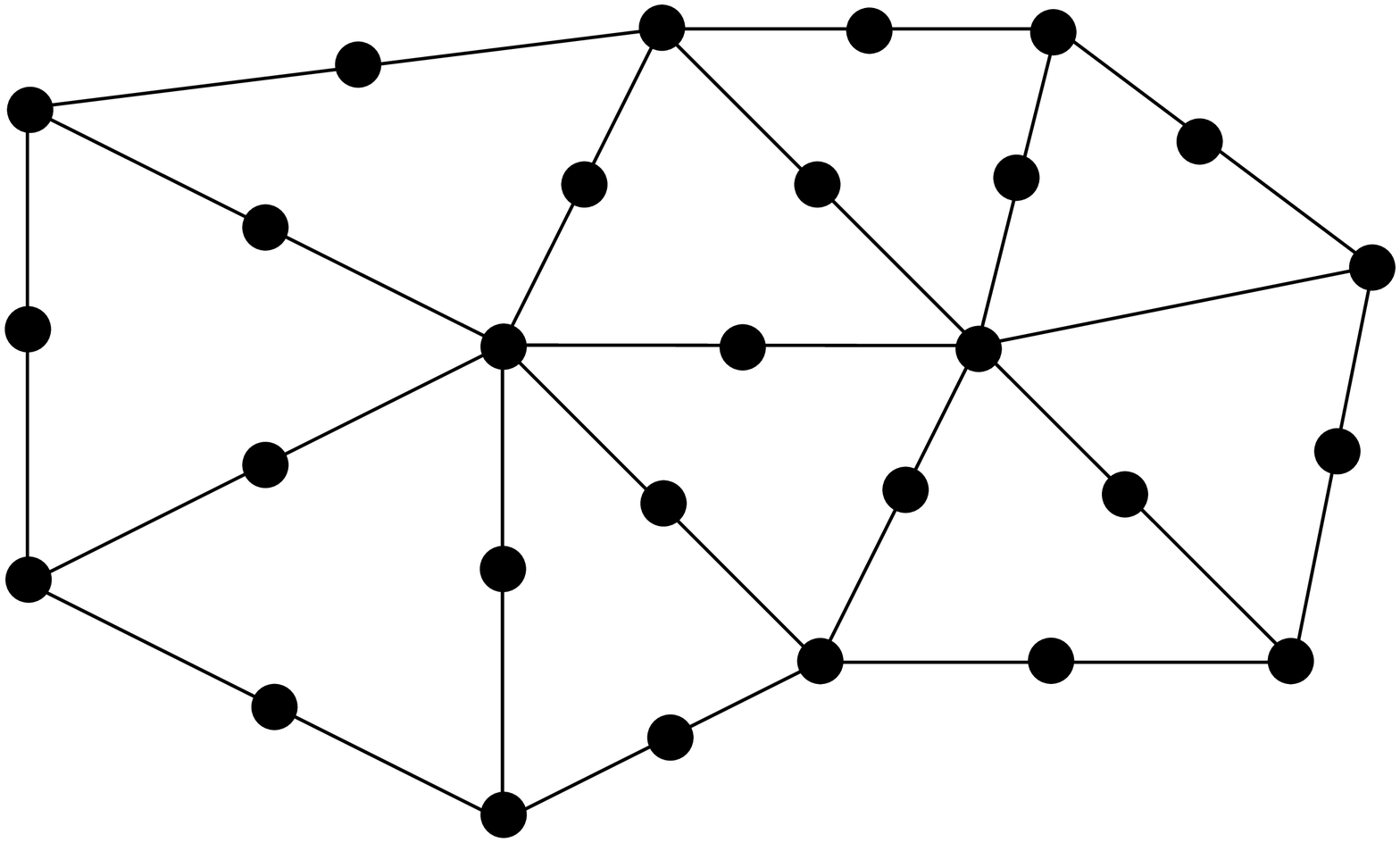}
    \caption{Distribution of the nodes on a triangulation of a domain
      $\Omega \subset \R^2$ for Lagrange finite elements of degree $q = 1$ (left) and
      $q = 2$ (right).}
    \label{fig:P1P2}
  \end{center}
\end{figure}

The basic example of finite element function spaces is given by the
family of Lagrange finite elements on simplices in $\R^d$. A Lagrange
finite element is given by a triple $(K,
\mathcal{P}_K,\mathcal{N}_K)$, where the $K$ is a simplex in $\R^d$ (a
line in $\R$, a triangle in $\R^2$, a tetrahedron in $\R^3$),
$\mathcal{P}_K$ is the space~$P_q(K)$ of scalar polynomials of
degree~$\leq q$ on~$K$ and each $\nu^K_i \in \mathcal{N}_K$ is given
by point evaluation at some point $x^K_i \in K$, as illustrated in
Figure~\ref{fig:P1P2} for $q = 1$ and $q = 2$ on a triangulation of
some domain~$\Omega \subset \R^2$. Note that by the placement of the
points $\{x^K_i\}_{i=1}^{n_K}$ at the vertices and edge midpoints of
each cell~$K$, the global function space is the set of continuous
piecewise polynomials of degree $q = 1$ and $q = 2$ respectively.

\subsubsection{The reference finite element}

\begin{figure}[htbp]
  \begin{center}
    \psfrag{p0}{$X^1 = (0,0)$}
    \psfrag{p1}{$X^2 = (1,0)$}
    \psfrag{p2}{$X^3 = (0,1)$}
    \psfrag{xi}{$X$}
    \psfrag{x}{$x = F_K(X)$}
    %\psfrag{F=}{$F_K(X) = x^1 \Phi_1(X) + x^2 \Phi_2(X) + x^3 \Phi_3(X)$}
    \psfrag{F=}{}
    \psfrag{F}{$F_K$}
    \psfrag{x0}{$x^1$}
    \psfrag{x1}{$x^2$}
    \psfrag{x2}{$x^3$}
    \psfrag{K0}{$K_0$}
    \psfrag{K}{$K$}
    \includegraphics[width=12cm]{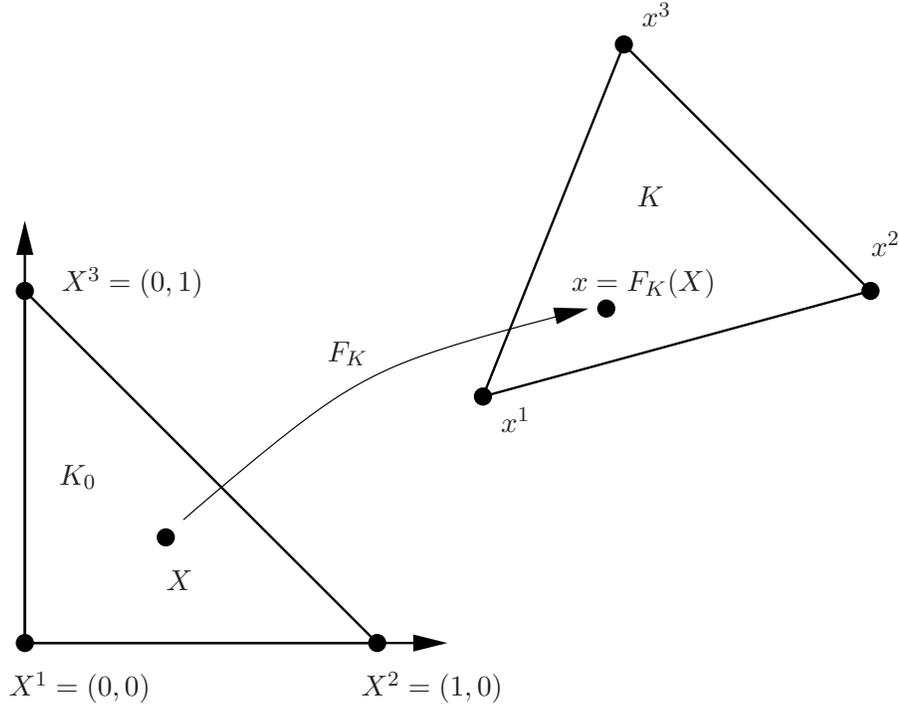}
    \caption{The (affine) mapping $F_K$ from a reference cell $K_0$
      to some cell $K \in \mathcal{T}$.}
    \label{fig:affinemap}
  \end{center}
\end{figure}

As we have seen, a global discrete function space $V_h$ may be
described by a mesh $\mathcal{T}$, a set of finite elements
$\{(K,\mathcal{P}_K,\mathcal{N}_K)\}_{K\in\mathcal{T}}$ and a set of
local-to-global mappings $\{\iota_K\}_{K\in\mathcal{T}}$. We may
simplify this description further by introducing a \emph{reference
finite element} $(K_0,\mathcal{P}_0,\mathcal{N}_0)$, where
$\mathcal{N}_0 = \{\nu_1^0,\nu_2^0,\ldots,\nu_{n_0}^0\}$, and a set of
invertible mappings $\{F_K\}_{K\in\mathcal{T}}$ that map the reference
cell~$K_0$ to the cells of the mesh,
\begin{equation}
  K = F_K(K_0) \quad \forall K \in \mathcal{T},
\end{equation}
as illustrated in Figure~\ref{fig:affinemap}. Note that $K_0$ is
generally not part of the mesh.  Typically, the mappings
$\{F_K\}_{K\in\mathcal{T}}$ are \emph{affine}, that is, each $F_K$ can
be written in the form $F_K(X) = A_K X + b_K$ for some matrix $A_K \in
\R^{d\times d}$ and some vector $b_K \in \R^d$, or
\emph{isoparametric}, in which case the components of $F_K$ are
functions in $\mathcal{P}_0$.

For each cell $K \in \mathcal{T}$, the mapping $F_K$ generates a
function space on $F_K$ given by
\begin{equation}
  \mathcal{P}_K = \{ v = v_0 \circ F_K^{-1} : v_0 \in \mathcal{P}_0 \},
\end{equation}
that is, each function $v = v(x)$ may be written in the form
$v(x) = v_0(F_K^{-1}(x)) = v_0 \circ F_K^{-1} (x)$ for some
$v_0 \in \mathcal{P}_0$.

Similarly, we may also generate for each $K \in \mathcal{T}$ a set of
nodes $\mathcal{N}_K$ on $\mathcal{P}_K$ given by
\begin{equation}
  \mathcal{N}_K = \{ \nu^K_i : \nu^K_i(v) = \nu^0_i(v \circ F_K),
  \quad i=1,2,\ldots,n_0 \}.
\end{equation}
Using the set of mappings $\{F_K\}_{K\in\mathcal{T}}$, we may thus
generate from the reference finite
element~$(K_0,\mathcal{P}_0,\mathcal{N}_0)$
a set of finite elements $\{(K,\mathcal{P}_K,\mathcal{N}_K)\}_{K\in\mathcal{T}}$
given by
\begin{equation} \label{eq:elementgeneration}
  \begin{split}
  K &= F_K(K_0), \\
  \mathcal{P}_K &= \{ v = v_0 \circ F_K^{-1} : v_0 \in \mathcal{P}_0 \}, \\
  \mathcal{N}_K &= \{ \nu^K_i : \nu^K_i(v) = \nu^0_i(v \circ F_K),
  \quad i=1,2,\ldots,n_0 = n_K \}.
  \end{split}
\end{equation}
With this construction, it is also simple to generate a set of nodal basis
functions $\{\phi^K_i\}_{i=1}^{n_K}$ on $K$ from a
set of nodal basis functions $\{\Phi_i\}_{i=1}^{n_0}$ on the
reference
element satisfying
$\nu^0_i(\Phi_j) = \delta_{ij}$. Noting that
if $\phi^K_i = \Phi_i \circ F_K^{-1}$ for
$i=1,2,\ldots,n_K$, then
\begin{equation}
  \nu^K_i(\phi^K_j) = \nu^0_i(\phi^K_j \circ F_K) = \nu^0_i(\Phi_j) = \delta_{ij},
\end{equation}
so $\{\phi^K_i\}_{i=1}^{n_K}$ is a nodal basis for $\mathcal{P}_K$.

Note that not all finite elements may be generated from a
reference finite element using this simple construction. For example,
this construction fails for the family of Hermite finite
elements~\cite{Cia76,Cia78,BreSco94}. Other examples include
$H(\mathrm{div})$ and $H(\mathrm{curl})$ conforming finite elements
(preserving the divergence and the curl respectively over cell
boundaries) which require a special mapping of the basis functions
from the reference element.

However, we shall limit our current discussion to finite elements that
can be generated from a reference finite element according
to~(\ref{eq:elementgeneration}), which includes all affine and
isoparametric finite elements with nodes given by point evaluation
such as the family of Lagrange finite elements on simplices.

\begin{figure}[htbp]
  \begin{center}
    \includegraphics[width=11cm]{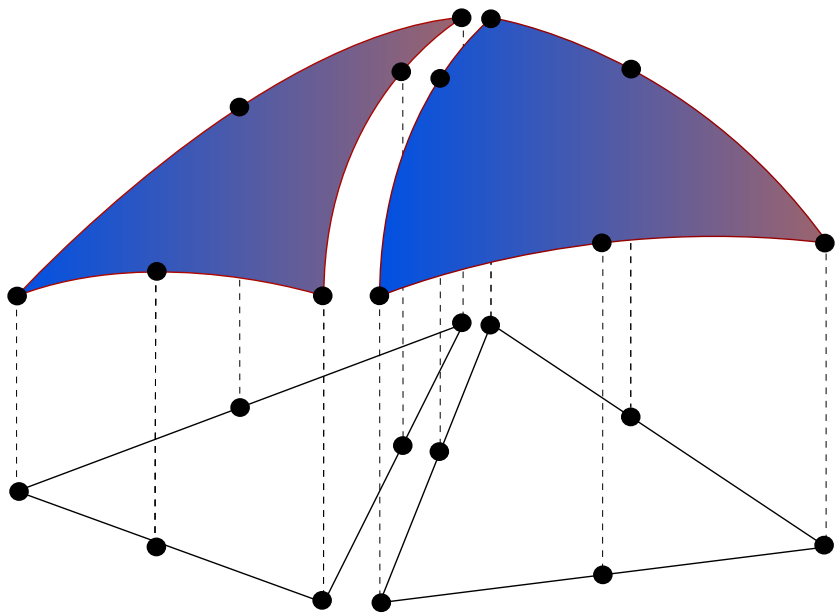}
    \caption{Piecing together local function spaces on the cells of a
    mesh to form a discrete function space on~$\Omega$, generated by
    a reference finite element $(K_0, \mathcal{P}_0, \mathcal{N}_0)$,
    a set of local-to-global mappings~$\{\iota_K\}_{K\in\mathcal{T}}$
    and a set of mappings~$\{F_K\}_{K\in\mathcal{T}}$.}
    \label{fig:femspace}
  \end{center}
\end{figure}

We may thus define a discrete function space by specifying a
mesh~$\mathcal{T}$, a reference finite element $(K, \mathcal{P}_0,
\mathcal{N}_0)$, a set of local-to-global mappings
$\{\iota_K\}_{K\in\mathcal{T}}$ and a set of mappings
$\{F_K\}_{K\in\mathcal{T}}$ from the reference cell $K_0$, as
demonstrated in Figure~\ref{fig:femspace}. Note that in general, the
mappings need not be of the same type for all cells $K$ and not all
finite elements need to be generated from the same reference finite
element. In particular, one could employ a different (higher-degree)
isoparametric mapping for cells on a curved boundary.

\subsection{The Variational Problem}

We shall assume that we are given a set of discrete function spaces
defined by a corresponding set of finite elements on some
triangulation~$\mathcal{T}$ of a domain~$\Omega
\subset \R^d$. In particular, we are given a pair of function spaces,
\begin{equation}
  \begin{split}
    \hat{V}_h &= \mathrm{span} \{\hat{\phi_i}\}_{i=1}^N, \\
    V_h &= \mathrm{span} \{\phi_i\}_{i=1}^N,
  \end{split}
\end{equation}
which we refer to as the test and trial spaces respectively.

We shall also assume that we are given a variational
problem of the form: Find $U \in V_h$ such that
\begin{equation} \label{eq:varproblem}
  a(U; v) = L(v) \quad \forall v \in \hat{V}_h,
\end{equation}
where $a : V_h \times \hat{V}_h \rightarrow \R$ is a semilinear form
which is linear in its second argument\footnotemark{} and $L : \hat{V}_h
\rightarrow \R$ is a linear form (functional). Typically, the
forms~$a$ and~$L$ of~(\ref{eq:varproblem}) are defined in terms of
integrals over the domain~$\Omega$ or subsets of the boundary
$\partial \Omega$ of $\Omega$.

\footnotetext{We shall use the convention that a semilinear form is
  linear in each of the arguments appearing after the semicolon.
  Furthermore, if a semilinear form~$a$ with two arguments is linear
  in both its arguments, we shall use the notation
  \begin{equation}
    a(v, U) = a'(U; v, U) = a'(U; v) \, U,
  \end{equation}
  where $a'$ is the Fr\'echet derivative of $a$ with respect to $U$,
  that is, we write the bilinear form with the test function as its
  first argument.}

\subsubsection{Nonlinear variational problems}

The variational problem~(\ref{eq:varproblem}) gives rise to
a system of discrete equations,
\begin{equation} \label{eq:system,nonlinear,again}
  F(U) = 0,
\end{equation}
for the vector~$(U_i) \in \R^N$ of degrees of freedom of the solution $U =
\sum_{i=1}^N U_i \phi_i \in V_h$, where
\begin{equation}
  F_i(U) = a(U; \hat{\phi}_i) - L(\hat{\phi}_i), \quad i = 1,2,\ldots,N.
\end{equation}

It may also be desirable to compute the Jacobian $A = F'$ of the
nonlinear system~(\ref{eq:system,nonlinear,again}) for use in a Newton's
method. We note that if the semilinear form~$a$ is differentiable in~$U$,
then the entries of the Jacobian~$A$ are given by
\begin{equation} \label{eq:jacobian}
  \begin{split}
    A_{ij}
    &= \frac{\partial F_i(U)}{\partial U_j}
    = \frac{\partial}{\partial U_j} a(U; \hat{\phi}_i)
    = a'(U; \hat{\phi}_i) \frac{\partial U}{\partial U_j}
    = a'(U; \hat{\phi}_i) \phi_j
    = a'(U; \hat{\phi}_i,\phi_j).
  \end{split}
\end{equation}

As an example, consider the nonlinear Poisson's equation
\begin{equation} \label{eq:poisson,nonlinear}
  \begin{split}
    - \nabla \cdot ((1+u) \nabla u) &= f \quad \mbox{ in } \Omega, \\
    u &= 0 \quad \mbox{ on } \partial\Omega.
  \end{split}
\end{equation}
Multiplying (\ref{eq:poisson,nonlinear}) with a test function~$v$ and
integrating by parts, we obtain
\begin{equation}
  \int_{\Omega} \nabla v \cdot ((1+u) \nabla u) \dx =
  \int_{\Omega} v f \dx,
\end{equation}
and thus a discrete nonlinear variational
problem of the form~(\ref{eq:varproblem}), where
\begin{equation}
  \begin{split}
    a(U; v) &= \int_{\Omega} \nabla v \cdot ((1+U) \nabla U) \dx, \\
    L(v) &= \int_{\Omega} v \, f \dx.
  \end{split}
\end{equation}
Linearizing the semilinear form~$a$ around $U$, we obtain
\begin{equation}
  a'(U; v, w) =
  \int_{\Omega} \nabla v \cdot (w \nabla U) \dx +
  \int_{\Omega} \nabla v \cdot ((1+U) \nabla w) \dx,
\end{equation}
for any $w \in V_h$. In particular, the entries of the Jacobian
matrix~$A$ are given by
\begin{equation}
  A_{ij} = a'(U; \hat{\phi}_i, \phi_j) =
  \int_{\Omega} \nabla \hat{\phi}_i \cdot (\phi_j \nabla U) \dx +
  \int_{\Omega} \nabla \hat{\phi}_i \cdot ((1+U) \nabla \phi_j) \dx.
\end{equation}

\subsubsection{Linear variational problems}

If the variational problem~(\ref{eq:varproblem}) is linear, the
nonlinear system~(\ref{eq:system,nonlinear,again}) is reduced to the linear system
\begin{equation}
  A U = b,
\end{equation}
for the degrees of freedom~$(U_i) \in \R^N$, where
\begin{equation}
  \begin{split}
    A_{ij} &= a(\hat{\phi}_i, \phi_j), \\
    b_i &= L(\hat{\phi}_i).
  \end{split}
\end{equation}
Note the relation to~(\ref{eq:jacobian}) in that
$A_{ij} = a(\hat{\phi}_i, \phi_j) = a'(U; \hat{\phi}_i, \phi_j)$.

In Section~\ref{sec:survey}, we saw the canonical example of a linear
variational problem with Poisson's equation,
\begin{equation}
  \begin{split}
    - \Delta u &= f \quad \mbox{ in } \Omega, \\
             u &= 0 \quad \mbox{ on } \partial\Omega,
  \end{split}
\end{equation}
corresponding to a discrete linear variational
problem of the form~(\ref{eq:varproblem}), where
\begin{equation}
  \begin{split}
    a(v, U) &= \int_{\Omega} \nabla v \cdot \nabla U \dx, \\
       L(v) &= \int_{\Omega} v \, f \dx.
  \end{split}
\end{equation}

\subsection{Multilinear Forms}

We find that for both nonlinear and linear problems, the system of
discrete equations is obtained from the given variational problem by
evaluating a set of multilinear forms on the set of basis
functions. Noting that the semilinear form~$a$ of the nonlinear
variational problem~(\ref{eq:varproblem}) is a linear form for any
given fixed~$U \in V_h$ and that the form~$a$ for a linear variational
problem can be expressed as $a(v, U) = a'(U; v, U)$, we thus need to be
able to evaluate the following multilinear forms:
\begin{equation}
  \begin{split}
    a(U; \cdot) : \hat{V}_h \rightarrow \R, \\
    L : \hat{V}_h \rightarrow \R, \\
    a'(U, \cdot, \cdot) : \hat{V}_h \times V_h \rightarrow \R.
  \end{split}
\end{equation}
We shall therefore consider the evaluation of general multilinear forms
of arity~$r > 0$,
\begin{equation} \label{eq:multilinear}
  a : V_h^1 \times V_h^2 \times \cdots \times V_h^r \rightarrow \R,
\end{equation}
defined on the product space $V_h^1 \times V_h^2 \times \cdots \times V_h^r$
of a given set $\{V_h^j\}_{j=1}^r$ of discrete function spaces on a
triangulation $\mathcal{T}$ of a domain $\Omega \subset \R^d$. In the
simplest case, all function spaces are equal but there are many
important examples, such as mixed methods, where it is important to
consider arguments coming from different function spaces. We shall
restrict our attention to multilinear forms expressed as integrals
over the domain~$\Omega$ (or subsets of its boundary).

Let now
$\{\phi_i^1\}_{i=1}^{N^1},
 \{\phi_i^2\}_{i=1}^{N^2}, \ldots,
 \{\phi_i^r\}_{i=1}^{N^r}$
be bases of $V_h^1, V_h^2, \ldots, V_h^r$ respectively and let $i =
(i_1, i_2, \ldots, i_r)$ be a multiindex of length $|i| = r$. The
multilinear form $a$ then defines a rank~$r$ tensor given by
\begin{equation} \label{eq:tensor}
  A_i = a(\phi_{i_1}^1, \phi_{i_2}^2, \ldots, \phi_{i_r}^r)
  \quad \forall i \in \mathcal{I},
\end{equation}
where $\mathcal{I}$ is the index set
\begin{equation}
  \mathcal{I} =
  \prod_{j=1}^r[1,|V^j_h|] = \{(1,1,\ldots,1), (1,1,\ldots,2), \ldots,
  (N^1,N^2,\ldots,N^r)\}.
\end{equation}
For any given multilinear form of arity~$r$, the tensor~$A$ is a
(typically sparse) tensor of rank~$r$ and dimension
$(|V_h^1|, |V_h^2|, \ldots, |V_h^r|) = (N^1, N^2,\ldots,N^r)$.

Typically, the arity of the multilinear form~$a$ is $r = 2$, that is,
$a$ is a bilinear form, in which case the corresponding tensor~$A$ is
a matrix (the ``stiffness matrix''), or the arity of the multilinear
form~$a$ is $r = 1$, that is, $a$ is a linear form, in which case the
corresponding tensor~$A$ is a vector (``the load vector'').

Sometimes it may also be of interest to consider forms of higher
arity. As an example, consider the discrete trilinear form $a : V_h^1
\times V_h^2 \times V_h^3 \rightarrow \R$ associated with the weighted
Poisson's equation~$-\nabla \cdot (w \nabla u) = f$. The trilinear
form~$a$ is given by
\begin{equation} \label{eq:poisson,weighted}
  a(v, U, w) = \int_{\Omega} w \nabla v \cdot \nabla U \dx, \\
\end{equation}
for $w = \sum_{i=1}^{N^3} w_i \phi_i^3 \in V_h^3$ a given discrete
weight function. The corresponding rank~three tensor is given by
\begin{equation}
  A_i = \int_{\Omega} \phi^3_{i_3}
  \nabla \phi^1_{i_1} \cdot \nabla \phi^2_{i_2} \dx.
\end{equation}
Noting that for any $w = \sum_{i=1}^{N^3} w_i \phi_i^3$, the tensor
contraction $A : w = \left(\sum_{i_3=1}^{N^3} A_{i_1 i_2 i_3}
w_{i_3}\right)_{i_1 i_2}$ is a matrix, we may thus obtain the solution $U$
by solving the linear system
\begin{equation}
  (A : w) U = b,
\end{equation}
where $b_i = L(\phi^1_{i_1}) = \int_{\Omega} \phi^1_{i_1} f \dx$.  Of
course, if the solution is needed only for one single weight
function~$w$, it is more efficient to consider $w$ as a fixed function
and directly compute the matrix~$A$ associated with the bilinear
form~$a(\cdot, \cdot, w)$. In some cases, it may even be desirable to
consider the function~$U$ as being fixed and directly compute a
vector~$A$ (the action) associated with the linear form~$a(\cdot, U,
w)$, as discussed above in Section~\ref{sec:analysa}. It is thus
important to consider multilinear forms of general arity~$r$.

\subsection{Assembling the Discrete System}
\label{sec:assembly}

The standard algorithm~\cite{ZieTay67,Hug87,Lan99} for computing the
tensor~$A$ is known as \emph{assembly}; the tensor is
computed by iterating over the cells of the mesh
$\mathcal{T}$ and adding from each
cell the local contribution to the global tensor $A$.

To explain how the standard assembly algorithm applies to the
computation of the tensor~$A$ defined in~(\ref{eq:tensor}) from a
given multilinear form~$a$, we note that if the multilinear form $a$ is
expressed as an integral over the domain $\Omega$, we can write the
multilinear form as a sum of element multilinear forms,
\begin{equation}
  a = \sum_{K\in\mathcal{T}} a_K,
\end{equation}
and thus
\begin{equation}
  A_i = \sum_{K\in\mathcal{T}}
  a_K(\phi_{i_1}^1, \phi_{i_2}^2, \ldots, \phi_{i_r}^r).
\end{equation}
We note that in the case of Poisson's equation, $-\Delta u = f$, the
element bilinear form $a_K$ is given by $a_K(v, U) = \int_{K}
\nabla v \cdot \nabla U \dx$.

We now let $\iota_K^j : [1,n_K^j] \rightarrow [1,N^j]$ denote the
local-to-global mapping introduced above in
Section~\ref{sec:femfunctions} for each discrete function space $V_h^j$,
$j=1,2,\ldots,r$, and define for each $K \in \mathcal{T}$ the
collective local-to-global mapping $\iota_K : \mathcal{I}_K
\rightarrow \mathcal{I}$ by
\begin{equation}
  \iota_K(i) =
  (\iota_K^1(i_1),\iota_K^2(i_2),\ldots,\iota_K^3(i_3))
  \quad \forall i \in \mathcal{I}_K,
\end{equation}
where $\mathcal{I}_K$ is the index set
\begin{equation}
  \mathcal{I}_K =
  \prod_{j=1}^r[1,|\mathcal{P}_K^j|] = \{(1,1,\ldots,1), (1,1,\ldots,2), \ldots,
  (n_K^1,n_K^2,\ldots,n_K^r)\}.
\end{equation}

Furthermore, for each $V_h^j$ we let $\{\phi^{K,j}_i\}_{i=1}^{n_K^j}$
denote the restriction to an element $K$ of the subset of the basis
$\{\phi_i^j\}_{i=1}^{N^j}$ of $V_h^j$ supported on $K$, and for
each $i \in \mathcal{I}$ we let
$\mathcal{T}_i \subset \mathcal{T}$ denote the subset of
cells on which all of the basis functions
$\{\phi_{i_j}^j\}_{j=1}^r$ are supported.

We may now compute the tensor~$A$ by summing the contributions from
each local cell~$K$,
\begin{equation}
  \begin{split}
  A_i
  &=
  \sum_{K\in\mathcal{T}}
  a_K(\phi_{i_1}^1, \phi_{i_2}^2, \ldots, \phi_{i_r}^r)
  =
  \sum_{K\in\mathcal{T}_i}
  a_K(\phi_{i_1}^1, \phi_{i_2}^2, \ldots, \phi_{i_r}^r) \\
  &=
  \sum_{K\in\mathcal{T}_i}
  a_K(\phi_{(\iota_K^1)^{-1}(i_1)}^{K,1},
      \phi_{(\iota_K^2)^{-1}(i_2)}^{K,2}, \ldots,
      \phi_{(\iota_K^r)^{-1}(i_r)}^{K,r}).
  \end{split}
\end{equation}
This computation may be carried out efficiently by iterating once over
all cells~$K \in \mathcal{T}$ and adding the contribution from each
$K$ to every entry $A_i$ of $A$ such that $K \in \mathcal{T}_i$, as
illustrated in Algorithm~\ref{alg:assembly,1}. In particular, we never
need to form the set $\mathcal{T}_i$, which is implicit through the
set of local-to-global mappings~$\{\iota_K\}_{K\in\mathcal{T}}$.

\begin{algorithm}
  \begin{tabbing}
    $A = 0$\\
    \textbf{for}  {$K \in \mathcal{T}$}\\
    \tab \textbf{for} $i \in \mathcal{I}_K$ \\
    \tab \tab %
	 {$A_{\iota_K(i)} =
	   A_{\iota_K(i)} +
           a_K(\phi_{i_1}^{K,1}, \phi_{i_2}^{K,2}, \ldots, \phi_{i_r}^{K,r})$} \\
    \tab \textbf{end for} \\
    \textbf{end for}
  \end{tabbing}
  \caption{$A$ = Assemble($a$, $\{V_h^j\}_{j=1}^r$, $\{\iota_K\}_{K\in\mathcal{T}}$, $\mathcal{T}$)}
  \label{alg:assembly,1}
\end{algorithm}

The assembly algorithm may be improved
by defining the \emph{element tensor}~$A^K$ by
\begin{equation} \label{eq:elementtensor}
  A^K_i = a_K(\phi_{i_1}^{K,1}, \phi_{i_2}^{K,2}, \ldots,
  \phi_{i_r}^{K,r}) \quad \forall i \in \mathcal{I}_K.
\end{equation}
For any multilinear form of arity~$r$, the element tensor $A^K$ is a
(typically dense) tensor of rank~$r$ and dimension $(n_K^1, n_K^2, \ldots,
n_K^r)$.

By computing first on each cell~$K$ the element tensor~$A^K$ before
adding the entries to the tensor~$A$ as in
Algorithm~\ref{alg:assembly,2}, one may take advantage of optimized
library routines for performing each of the two steps. Note that
Algorithm~\ref{alg:assembly,2} is independent of the algorithm used to
compute the element tensor.

\begin{algorithm}
  \begin{tabbing}
    {$A = 0$} \\
    \textbf{for} {$K \in \mathcal{T}$} \\
    \tab Compute $A^K$ according to (\ref{eq:elementtensor}) \\
    \tab Add $A^K$ to $A$ according to $\iota_K$ \\
    \textbf{end for}
  \end{tabbing}
  \caption{$A$ = Assemble($a$, $\{V_h^j\}_{j=1}^r$, $\{\iota_K\}_{K\in\mathcal{T}}$, $\mathcal{T}$)}
  \label{alg:assembly,2}
\end{algorithm}

\begin{figure}[htbp]
  \begin{center}
    \psfrag{i0}{\hspace{-0.5cm}$\iota_K^1(1)$}
    \psfrag{i1}{\hspace{-0.5cm}$\iota_K^1(2)$}
    \psfrag{i2}{\hspace{-0.5cm}$\iota_K^1(3)$}
    \psfrag{j0}{\hspace{-0.3cm}$\iota_K^2(1)$}
    \psfrag{j1}{\hspace{-0.5cm}$\iota_K^2(2)$}
    \psfrag{j2}{\hspace{-0.1cm}$\iota_K^2(3)$}
    \psfrag{A21}{$A^K_{32}$}
    \psfrag{1}{$1$}
    \psfrag{2}{$2$}
    \psfrag{3}{$3$}
    \includegraphics[width=12cm]{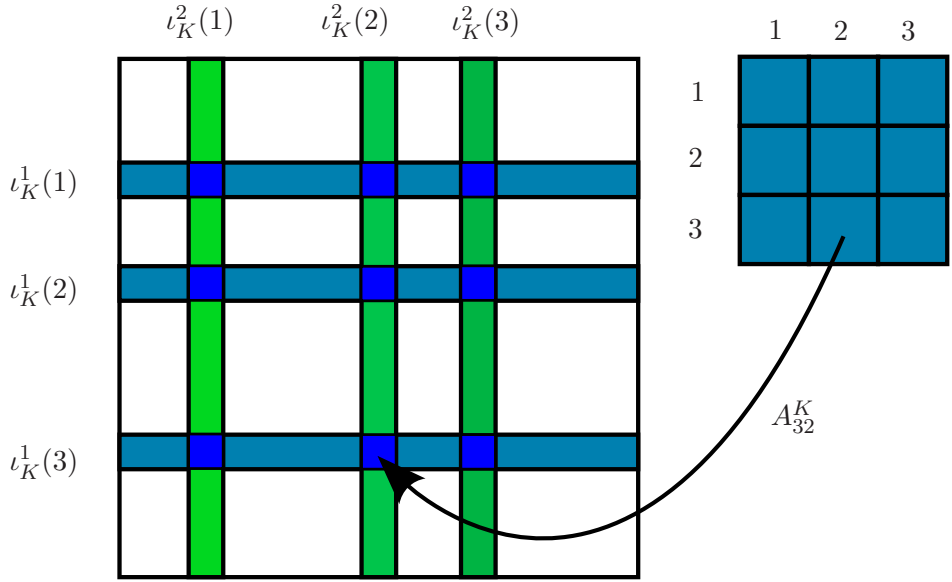}
    \caption{Adding the entries of the element tensor~$A^K$ to the
      global tensor~$A$ using the  local-to-global mapping
      $\iota_K$, illustrated here for a rank two
      tensor (a matrix).}
    \label{fig:insertion}
  \end{center}
\end{figure}

Considering first the second operation of inserting (adding) the
entries of~$A^K$ into the global sparse tensor~$A$, this may in
principle be accomplished by iterating over all $i \in I_K$ and adding
the entry $A^K_i$ at position $\iota_K(i)$ of $A$ as illustrated in
Figure~\ref{fig:insertion}. However, sparse matrix libraries such as
PETSc \cite{www:PETSc,BalBus04,BalEij97} often provide optimized
routines for this type of operation, which may significantly improve
the performance compared to accessing each entry of~$A$ individually
as in Algorithm~\ref{alg:assembly,1}. Even so, the cost of adding
$A^K$ to $A$ may be substantial even with an efficient implementation
of the sparse data structure for $A$, see \cite{logg:article:07}.

A similar approach can be taken to the first step of computing the
element tensor, that is, an optimized library routine is called to
compute the element tensor. Because of the wide variety of multilinear
forms that appear in applications, a separate implementation is needed
for any given multilinear form. Therefore, the implementation of this
code is often left to the user, as illustrated above in
Section~\ref{sec:deal} and Section~\ref{sec:diffpack}, but the code in
question may also be automatically generated and optimized for each
given multilinear form. We shall return to this question below in
Section~\ref{sec:automation,forms} and Section~\ref{sec:prototype}.

\subsection{Summary}

If we thus view the finite element method as a machine that automates
the discretization of differential equations, or more precisely, a
machine that generates the system of discrete
equations~(\ref{eq:system,nonlinear,again}) from a given variational
problem~(\ref{eq:varproblem}), an automation of the finite element
method is straightforward up to the point of computing the element
tensor for any given multilinear form and the local-to-global mapping
for any given discrete function space; if the element tensor $A^K$ and
the local-to-global mapping $\iota_K$ can be computed on any given
cell $K$, the global tensor~$A$ may be computed by
Algorithm~\ref{alg:assembly,2}.

Assuming now that each of the discrete function spaces involved in the
definition of the variational problem~(\ref{eq:varproblem}) is generated
on some mesh $\mathcal{T}$ of the domain $\Omega$ from
some reference finite element
$(K_0, \mathcal{P}_0,\mathcal{N}_0)$ by a set of local-to-global
mappings $\{\iota_K\}_{K\in\mathcal{T}}$ and a set of mappings
$\{F_K\}_{K\in\mathcal{T}}$ from the reference cell $K_0$, as
discussed in Section~\ref{sec:femfunctions}, we identify the following
key steps towards an automation of the finite element method:
\begin{itemize}
\item
  the automatic and efficient tabulation of the nodal basis functions
  on the reference finite element $(K_0,\mathcal{P}_0,\mathcal{N}_0)$;
\item
  the automatic and efficient evaluation of the element tensor~$A^K$ on each
  cell $K \in \mathcal{T}$;
\item
  the automatic and efficient assembly of the global tensor $A$ from
  the set of element tensors $\{A^K\}_{K\in\mathcal{T}}$ and the set
  of local-to-global mappings $\{\iota_K\}_{K\in\mathcal{T}}$.
\end{itemize}
We discuss each of these key steps below.

%------------------------------------------------------------------------------
\section{AUTOMATING THE TABULATION OF BASIS FUNCTIONS}
\label{sec:automation,tabulation}

\begin{figure}[htbp]
  \begin{center}
    \psfrag{p0}{}
    \psfrag{p1}{}
    \psfrag{p2}{}
    \psfrag{p3}{}
    \psfrag{v0}{$v^1$}
    \psfrag{v1}{$v^2$}
    \psfrag{v2}{$v^3$}
    \psfrag{v3}{$v^4$}
    \psfrag{x0}{$X_1$}
    \psfrag{x1}{$X_2$}
    \psfrag{x2}{$X_3$}
    \includegraphics[width=6.8cm]{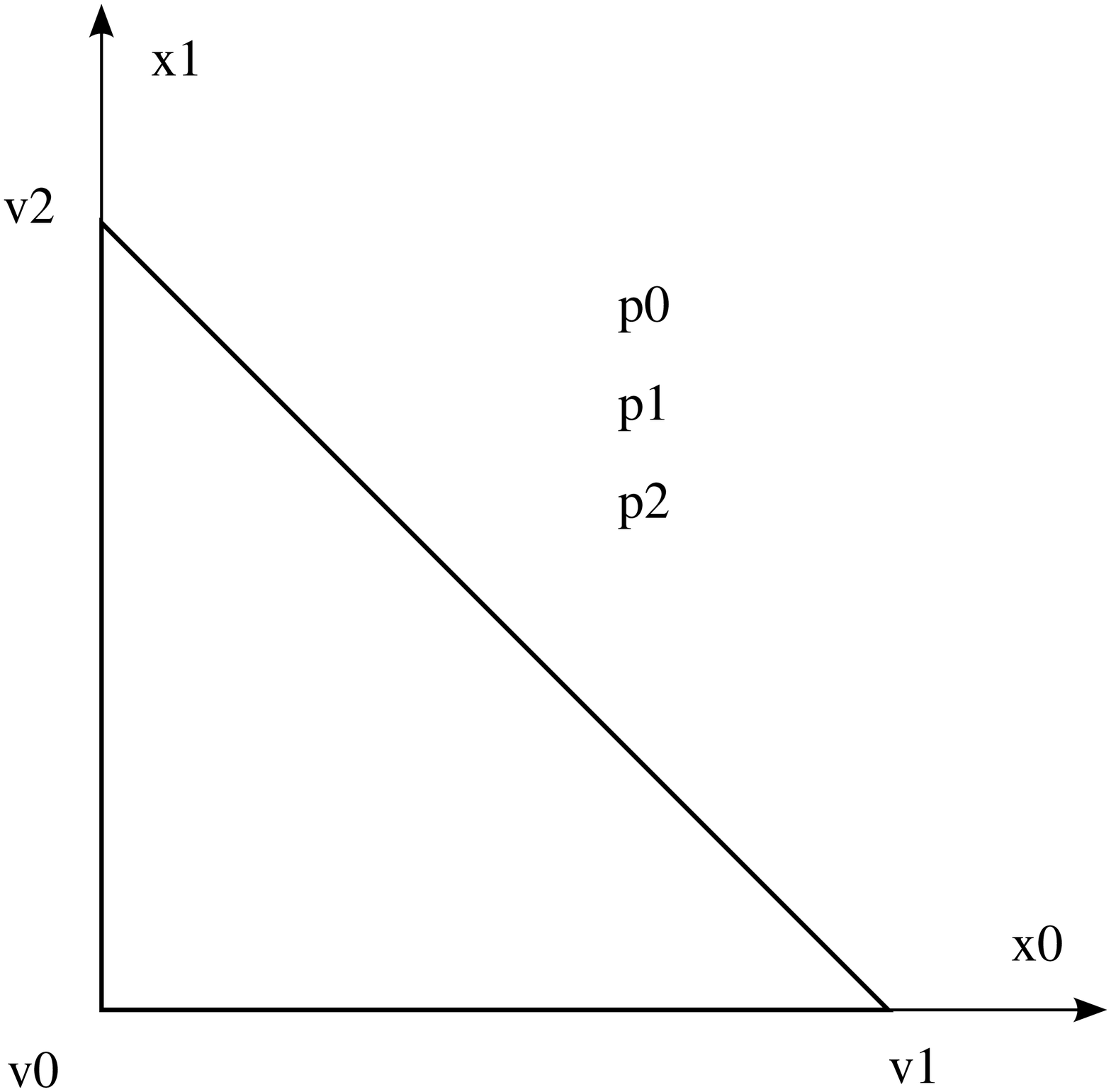} \quad\quad
    \includegraphics[width=6.8cm]{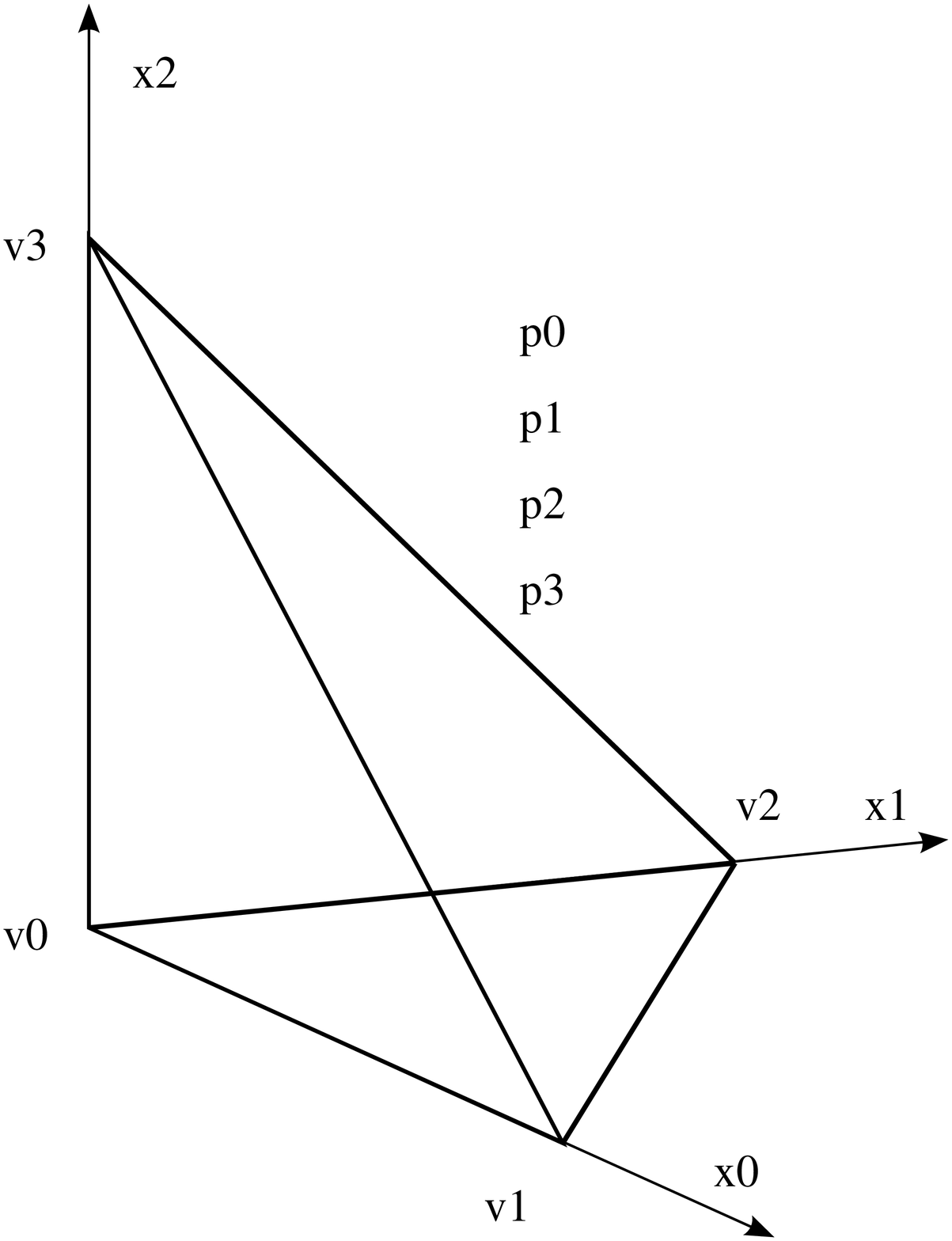}
    \caption{The reference triangle (left) with vertices at
      $v^1 = (0,0)$, $v^2 = (1,0)$ and $v^3 = (0,1)$, and the reference
      tetrahedron (right) with vertices at
      $v^1 = (0,0,0)$, $v^2 = (1,0,0)$, $v^3 = (0,1,0)$ and $v^4 = (0,0,1)$.}
    \label{fig:referenceshapes}
  \end{center}
\end{figure}

Given a reference finite element $(K_0,\mathcal{P}_0,\mathcal{N}_0)$,
we wish to generate the unique nodal basis
$\{\Phi_i\}_{i=1}^{n_0}$ for $\mathcal{P}_0$ satisfying
\begin{equation} \label{eq:nodalcondition}
  \nu^0_i(\Phi_j) = \delta_{ij}, \quad i,j = 1,2,\ldots,n_0.
\end{equation}
In some simple cases, these nodal basis functions can be worked out
analytically by hand or found in the literature,
see for example~\cite{ZieTay67,Hug87}. As a concrete example,
consider the nodal basis functions in the case when $\mathcal{P}_0$ is
the set of quadratic polynomials on the reference triangle $K_0$ with
vertices at $v^1 = (0,0)$, $v^2 = (1,0)$ and $v^3 = (0,1)$ as in
Figure~\ref{fig:referenceshapes} and nodes $\mathcal{N}_0 =
\{\nu^0_1,\nu^0_2,\ldots,\nu^0_6\}$ given by point evaluation at the
vertices and edge midpoints. A basis for $\mathcal{P}_0$ is then
given by
\begin{equation} \label{eq:P2}
  \begin{split}
    \Phi_1(X) &= (1 - X_1 - X_2)(1 - 2X_1 - 2X_2), \\
    \Phi_2(X) &= X_1 (2X_1 - 1), \\
    \Phi_3(X) &= X_2 (2X_2 - 1), \\
    \Phi_4(X) &= 4 X_1 X_2, \\
    \Phi_5(X) &= 4 X_2 (1 - X_1 - X_2), \\
    \Phi_6(X) &= 4 X_1 (1 - X_1 - X_2),
  \end{split}
\end{equation}
and it is easy to verify that this is the nodal basis.
However, in the general case, it may be very difficult
to obtain analytical expressions for the nodal basis
functions. Furthermore, copying the often complicated analytical
expressions into a computer program is prone to errors and may even
result in inefficient code.

In recent work, Kirby~\cite{www:FIAT,Kir04,Kir06} has proposed a
solution to this problem; by expanding the nodal basis functions
for~$\mathcal{P}_0$ as linear combinations of another (non-nodal)
basis for $\mathcal{P}_0$ which is easy to compute, one may translate
operations on the nodal basis functions, such as evaluation and
differentiation, into linear algebra operations on the expansion
coefficients.

This new linear algebraic approach to computing and representing finite
element basis functions removes the need for having explicit
expressions for the nodal basis functions, thus simplifying or
enabling the implementation of complicated finite elements.

\subsection{Tabulating Polynomial Spaces}

To generate the set of nodal basis functions~$\{\Phi_i\}_{i=1}^{n_0}$
for $\mathcal{P}_0$, we must first identify some other known basis
$\{\Psi_i\}_{i=1}^{n_0}$ for $\mathcal{P}_0$, referred to
in~\cite{Kir04} as the \emph{prime basis}. We return to the question
of how to choose the prime basis below.

Writing now each $\Phi_i$ as a linear combination of the prime basis
functions with $\alpha \in \R^{d\times d}$ the matrix of coefficients,
we have
\begin{equation}
  \Phi_i = \sum_{j=1}^{n_0} \alpha_{ij} \Psi_j,
  \quad i = 1,2,\ldots,n_0.
\end{equation}
The conditions~(\ref{eq:nodalcondition}) thus translate into
\begin{equation}
  \delta_{ij} = \nu^0_i(\Phi_j) =
  \sum_{k=1}^{n_0} \alpha_{jk} \nu^0_i(\Psi_k),
  \quad i,j = 1,2,\ldots,n_0,
\end{equation}
or
\begin{equation} \label{eq:fiatsystem}
  \mathcal{V} \alpha^{\top} = I,
\end{equation}
where $\mathcal{V} \in \R^{n_0\times n_0}$ is the (Vandermonde) matrix
with entries $\mathcal{V}_{ij} = \nu^0_i(\Psi_j)$ and $I$ is the
$n_0\times n_0$ identity matrix. Thus, the nodal
basis~$\{\Phi_i\}_{i=1}^{n_0}$ is easily computed by first computing
the matrix~$\mathcal{V}$ by evaluating the nodes at the prime basis
functions and then solving the linear system~(\ref{eq:fiatsystem}) to
obtain the matrix~$\alpha$ of coefficients.

In the simplest case, the space~$\mathcal{P}_0$ is the set
$P_q(K_0)$ of polynomials of degree~$\leq q$ on $K_0$. For
typical reference cells, including the reference triangle and the
reference tetrahedron shown in Figure~\ref{fig:referenceshapes},
orthogonal prime bases are available with simple recurrence
relations for the evaluation of the basis functions and their derivatives,
see for example~\cite{Dub91}. If $\mathcal{P}_0 = P_q(K_0)$,
it is thus straightforward to evaluate the prime basis and thus to
generate and solve the linear system~$(\ref{eq:fiatsystem})$ that
determines the nodal basis.

\subsection{Tabulating Spaces with Constraints}

In other cases, the space~$\mathcal{P}_0$ may be defined as some
subspace of $P_q(K_0)$, typically by constraining certain
derivatives of the functions in $\mathcal{P}_0$ or the functions
themselves to lie in $P_{q'}(K_0)$ for some $q ' < q$ on
some part of $K_0$. Examples include the
the Raviart--Thomas~\cite{RavTho77a},
Brezzi--Douglas--Fortin--Marini~\cite{BreFor91}
and Arnold--Winther \cite{ArnWin02} elements, which put constraints on
the derivatives of the functions in $\mathcal{P}_0$.

Another more obvious example, taken from~\cite{Kir04}, is the case
when the functions in $\mathcal{P}_0$ are constrained to
$P_{q-1}(\gamma_0)$ on some part $\gamma_0$ of the boundary
of $K_0$ but are otherwise in $P_q(K_0)$, which may be used
to construct the function space on a $p$-refined cell $K$ if the
function space on a neighboring cell $K'$ with common boundary
$\gamma_0$ is only $P_{q-1}(K')$. We may then define the
space~$\mathcal{P}_0$ by
\begin{equation}
  \mathcal{P}_0 = \{ v \in P_q(K_0) :
  v |_{\gamma_0} \in P_{q-1}(\gamma_0) \}
  =
  \{ v \in P_q(K_0) : l(v) = 0 \},
\end{equation}
where the linear functional $l$ is given by integration against the
$q$th degree Legendre polynomial along the boundary~$\gamma_0$.

In general, one may define a set $\{l_i\}_{i=1}^{n_c}$ of linear
functionals (constraints) and define $\mathcal{P}_0$ as the
intersection of the null spaces of these linear functionals on
$P_q(K_0)$,
\begin{equation}
  \mathcal{P}_0 = \{ v \in P_q(K_0) :
  l_i(v) = 0, \quad i = 1,2,\ldots,n_c \}.
\end{equation}
To find a prime basis~$\{\Psi_i\}_{i=1}^{n_0}$ for
$\mathcal{P}_0$, we note that any function in $\mathcal{P}_0$
may be expressed as a linear combination of some basis functions
$\{\bar{\Psi}_i\}_{i=1}^{|P_q(K_0)|}$ for $P_q(K_0)$, which we may
take as the orthogonal basis discussed above. We find that if $\Psi =
\sum_{i=1}^{|P_q(K_0)|} \beta_i \bar{\Psi}_i$, then
\begin{equation}
  0 = l_i(\Psi) = \sum_{j=1}^{|P_q(K_0)|}
  \beta_j l_i(\bar{\Psi}_j), \quad i=1,2,\ldots,n_c,
\end{equation}
or
\begin{equation}
  L \beta = 0,
\end{equation}
where $L$ is the $n_c \times |P_q(K_0)|$ matrix with
entries
\begin{equation}
  L_{ij} = l_i(\bar{\Psi}_j),
  \quad i=1,2,\ldots,n_c, \quad j=1,2,\ldots,|P_q(K_0)|.
\end{equation}
A prime basis for $\mathcal{P}_0$ may thus be found by computing
the nullspace of the matrix $L$, for example by computing its
singular value decomposition (see~\cite{GolLoa96}).
Having thus found the prime basis~$\{\Psi_i\}_{i=1}^{n_0}$, we may
proceed to compute the nodal basis as before.

%------------------------------------------------------------------------------
\section{AUTOMATING THE COMPUTATION OF THE ELEMENT TENSOR}
\label{sec:automation,forms}

As we saw in Section~\ref{sec:assembly}, given a multilinear form~$a$
defined on the product space $V_h^1 \times V_h^2 \times \ldots \times
V_h^r$, we need to compute for each cell $K \in \mathcal{T}$ the
rank~$r$ element tensor~$A^K$ given by
\begin{equation}
  A^K_i = a_K(\phi_{i_1}^{K,1}, \phi_{i_2}^{K,2}, \ldots,
  \phi_{i_r}^{K,r}) \quad \forall i \in \mathcal{I}_K,
\end{equation}
where $a_K$ is the local contribution to the multilinear form~$a$ from
the cell~$K$.

We investigate below two very different ways to compute the element
tensor, first a modification of the standard approach based on
quadrature and then a novel approach based on a special tensor
contraction representation of the element tensor, yielding speedups of
several orders of magnitude in some cases.

\subsection{Evaluation by Quadrature}
\label{sec:quadrature}

The element tensor~$A^K$ is typically evaluated by quadrature on the
cell~$K$. Many finite element libraries like
Diffpack~\cite{www:diffpack,Lan99} and
deal.II~\cite{www:deal.II,BanKan99,Ban00} provide the values of
relevant quantities like basis functions and their
derivatives at the quadrature points on~$K$ by mapping precomputed
values of the corresponding basis functions on the reference cell
$K_0$ using the mapping $F_K : K_0 \rightarrow K$.

Thus, to evaluate the element tensor~$A^K$ for Poisson's equation by
quadrature on $K$, one computes
\begin{equation} \label{eq:poisson,quadrature}
  A^K_i =
  \int_K
  \nabla \phi_{i_1}^{K,1} \cdot
  \nabla \phi_{i_2}^{K,2} \dx
  \approx
  \sum_{k=1}^{N_q}
  w_k
  \nabla \phi_{i_1}^{K,1}(x^k) \cdot
  \nabla \phi_{i_2}^{K,2}(x^k)
  \det F_K'(x^k),
\end{equation}
for some suitable set of quadrature points $\{x^i\}_{i=1}^{N_q}
\subset K$ with corresponding quadrature weights
$\{w_i\}_{i=1}^{N_q}$, where we assume that the quadrature weights are
scaled so that $\sum_{i=1}^{N_q} w_i = |K_0|$. Note that the
approximation~(\ref{eq:poisson,quadrature}) can be made exact for a
suitable choice of quadrature points if the basis functions are
polynomials.

Comparing~(\ref{eq:poisson,quadrature}) to the example codes in
Table~\ref{tab:poisson,deal} and Table~\ref{tab:poisson,diffpack}, we
note the similarities between~(\ref{eq:poisson,quadrature}) and the
two codes.  In both cases, the gradients of the basis functions as
well as the products of quadrature weight and the determinant of
$F_K'$ are precomputed at the set of quadrature points and then
combined to produce the integral~(\ref{eq:poisson,quadrature}).

If we assume that the two discrete spaces $V_h^1$ and $V_h^2$ are
equal, so that the local basis functions
$\{\phi^{K,1}_i\}_{i=1}^{n_K^1}$ and
$\{\phi^{K,2}_i\}_{i=1}^{n_K^2}$ are all generated from the same
basis $\{\Phi_i\}_{i=1}^{n_0}$ on the reference cell~$K_0$, the work
involved in precomputing the gradients of the basis functions at the
set of quadrature points amounts to computing for each quadrature point
$x_k$ and each basis function $\phi^K_i$ the matrix--vector product
$\nabla_x \phi^K_i(x_k) = (F_K')^{-\top}(x_k) \nabla_X \Phi_i(X_k)$, that is,
\begin{equation}
  \frac{\partial \phi^K_i}{\partial x_j}(x^k)
  = \sum_{l=1}^d
  \frac{\partial X_l}{\partial x_j} (x^k)
  \frac{\partial \Phi^K_i}{\partial X_l} (X^k),
\end{equation}
where $x^k = F_K(X^k)$ and $\phi^K_i = \Phi_i \circ F_K^{-1}$. Note
that the the gradients $\{\nabla_X \Phi_i(X_k)\}_{i=1,k=1}^{n_0,N_q}$
of the reference element basis functions at the set of quadrature
points on the reference element remain constant throughout the
assembly process and may be pretabulated and stored. Thus, the
gradients of the basis functions on $K$ may be computed in $N_q n_0
d^2$ multiply--add pairs (MAPs) and the total work to compute the
element tensor~$A^K$ is $N_q n_0 d^2 + N_q n_0^2 (d+2) \sim N_q n_0^2
d$, if we ignore that we also need to compute the mapping~$F_K$, and
the determinant and inverse of $F_K'$.  In
Section~\ref{sec:representation} and Section~\ref{sec:optimizations}
below, we will see that this operation count may be significantly
reduced.

\subsection{Evaluation by Tensor Representation}
\label{sec:representation}

It has long been known that it is sometimes possible to speed up the
computation of the element tensor by precomputing certain integrals on
the reference element. Thus, for any specific multilinear form, it may
be possible to find quantities that can be precomputed in order to
optimize the code for the evaluation of the element tensor. These
ideas were first introduced in a general setting in
\cite{logg:article:07,KirKne05} and later formalized and automated in
\cite{logg:article:10,logg:submit:toms-ffc-monomial-2006}. A similar
approach was implemented in early versions of
DOLFIN~\cite{logg:www:01,logg:preprint:06,logg:manual:01}, but only
for piecewise linear elements.

We first consider the case when the mapping~$F_K$ from the reference cell
is affine, and then discuss possible extensions to non-affine
mappings such as when $F_K$ is the isoparametric mapping. As a first
example, we consider again the computation of the element
tensor~$A^K$ for Poisson's equation. As before, we have
\begin{equation} \label{eq:poisson,A}
  A^K_i =
  \int_K
  \nabla \phi_{i_1}^{K,1} \cdot
  \nabla \phi_{i_2}^{K,2} \dx
  =
  \int_K
  \sum_{\beta=1}^d
  \frac{\partial \phi_{i_1}^{K,1}}{\partial x_{\beta}}
  \frac{\partial \phi_{i_2}^{K,2}}{\partial x_{\beta}} \dx,
\end{equation}
but instead of evaluating the gradients on $K$ and then proceeding to
evaluate the integral by quadrature, we make a change of
variables to write
\begin{equation}
  A^K_i =
  \int_{K_0}
  \sum_{\beta=1}^d
  \sum_{\alpha_1=1}^d
  \frac{\partial X_{\alpha_1}}{\partial x_{\beta}}
  \frac{\partial \Phi^1_{i_1}}{\partial X_{\alpha_1}}
  \sum_{\alpha_2=1}^d
  \frac{\partial X_{\alpha_2}}{\partial x_{\beta}}
  \frac{\partial \Phi^2_{i_2}}{\partial X_{\alpha_2}}
  \det F_K'
  \dX,
\end{equation}
and thus, if the mapping $F_K$ is affine so that the transforms
$\partial X / \partial x$ and the determinant~$\det F_K'$ are
constant, we obtain
\begin{equation}
  A^K_i =
  \det F_K'
  \sum_{\alpha_1=1}^d
  \sum_{\alpha_2=1}^d
  \sum_{\beta=1}^d
  \frac{\partial X_{\alpha_1}}{\partial x_{\beta}}
  \frac{\partial X_{\alpha_2}}{\partial x_{\beta}}
  \int_{K_0}
  \frac{\partial \Phi^1_{i_1}}{\partial X_{\alpha_1}}
  \frac{\partial \Phi^2_{i_2}}{\partial X_{\alpha_2}}
  \dX
  =
  \sum_{\alpha_1=1}^d
  \sum_{\alpha_2=1}^d
  A^0_{i\alpha} G_K^{\alpha},
\end{equation}
or
\begin{equation} \label{eq:poisson,contraction}
  A^K = A^0 : G_K,
\end{equation}
where
\begin{equation} \label{eq:poisson,AandG}
  \begin{split}
    A^0_{i\alpha}
    &=
    \int_{K_0}
    \frac{\partial \Phi^1_{i_1}}{\partial X_{\alpha_1}}
    \frac{\partial \Phi^2_{i_2}}{\partial X_{\alpha_2}}
    \dX, \\
    G_K^{\alpha}
    &=
    \det F_K'
    \sum_{\beta=1}^d
    \frac{\partial X_{\alpha_1}}{\partial x_{\beta}}
    \frac{\partial X_{\alpha_2}}{\partial x_{\beta}}.
  \end{split}
\end{equation}
We refer to the tensor~$A^0$ as the \emph{reference tensor}
and to the tensor~$G_K$ as the \emph{geometry tensor}.

Now, since the reference tensor is constant and does not depend on the
cell $K$, it may be precomputed before the assembly of the global
tensor~$A$. For the current example, the work on each cell~$K$ thus
involves first computing the rank~two geometry tensor $G^K$, which may
be done in $d^3$ multiply--add pairs, and then computing the rank~two
element tensor~$A^K$ as the tensor
contraction~(\ref{eq:poisson,contraction}), which may be done in
$n_0^2 d^2$ multiply--add pairs. Thus, the total operation count is
$d^3 + n_0^2 d^2 \sim n_0^2 d^2$, which should be compared to $N_q
n_0^2 d$ for the standard quadrature-based approach.  The speedup in
this particular case is thus roughly a factor $N_q/d$, which may be a
significant speedup, in particular for higher order elements.

As we shall see, the tensor representation~(\ref{eq:poisson,contraction})
generalizes to other multilinear forms as well. To see this,
we need to make some assumptions about the structure of the
multilinear form~(\ref{eq:multilinear}). We shall assume that the
multilinear form~$a$ is expressed as an integral over~$\Omega$ of a
weighted sum of products of basis functions or derivatives of basis
functions. In particular, we shall assume that the element
tensor~$A^K$ can be expressed as a sum, where each term takes the
following canonical form,
\begin{equation} \label{eq:canonical}
  A^K_i =
  \sum_{\gamma\in\mathcal{C}}
  \int_K
  \prod_{j=1}^m
  c_j(\gamma)
  D_x^{\delta_j(\gamma)}
  \phi^{K,j}_{\iota_j(i,\gamma)}[\kappa_j(\gamma)] \dx,
\end{equation}
where $\mathcal{C}$ is some given set of multiindices, each
coefficient $c_j$ maps the multiindex $\gamma$ to a real
number, $\iota_j$ maps $(i, \gamma)$ to a basis function index,
$\kappa_j$ maps $\gamma$ to a component index (for vector or tensor
valued basis functions) and $\delta_j$ maps $\gamma$ to a derivative
multiindex. To distinguish component indices from indices for basis
functions, we use~$[\cdot]$ to denote a component index and subscript
to denote a basis function index. In the simplest case, the number of
factors $m$ is equal to the arity~$r$ of the multilinear form (rank of
the tensor), but in general, the canonical form~(\ref{eq:canonical})
may contain factors that correspond to additional functions which are
not arguments of the multilinear form. This is the case for the
weighted Poisson's equation~(\ref{eq:poisson,weighted}), where $m = 3$
and $r = 2$. In general, we thus have $m > r$.

As an illustration of this notation, we consider again the bilinear
form for Poisson's equation and write it in the notation
of~(\ref{eq:canonical}). We will also consider a more involved example
to illustrate the generality of the notation.  From
(\ref{eq:poisson,A}), we have
\begin{equation}
  A^K_i
  =
  \int_K
  \sum_{\gamma=1}^d
  \frac{\partial \phi_{i_1}^{K,1}}{\partial x_{\gamma}}
  \frac{\partial \phi_{i_2}^{K,2}}{\partial x_{\gamma}} \dx
  =
  \sum_{\gamma=1}^d
  \int_K
  \frac{\partial \phi_{i_1}^{K,1}}{\partial x_{\gamma}}
  \frac{\partial \phi_{i_2}^{K,2}}{\partial x_{\gamma}} \dx,
\end{equation}
and thus, in the notation of~(\ref{eq:canonical}),
\begin{equation}
  \begin{split}
    m               &= 2, \\
    \mathcal{C}     &= [1,d], \\
    c(\gamma)       &= (1,1), \\
    \iota(i,\gamma) &= (i_1, i_2), \\
    \kappa(\gamma)  &= (\emptyset, \emptyset), \\
    \delta(\gamma)  &= (\gamma, \gamma),
  \end{split}
\end{equation}
where $\emptyset$ denotes an empty component index (the basis functions
are scalar).

As another example, we consider the bilinear form for a
stabilization term appearing in a least-squares stabilized
$\mathrm{cG}(1)\mathrm{cG}(1)$ method for the incompressible
Navier--Stokes equations~\cite{EriEst03c,HofJoh04a,HofJoh04b,HofJoh06},
\begin{equation} \label{eq:stabilization}
  a(v, U) = \int_{\Omega}
  (w \cdot \nabla v) \cdot
  (w \cdot \nabla U) \dx
  =
  \int_{\Omega}
  \sum_{\gamma_1,\gamma_2,\gamma_3=1}^d
  w[\gamma_2] \frac{\partial v[\gamma_1]}{\partial x_{\gamma_2}}
  w[\gamma_3] \frac{\partial U[\gamma_1]}{\partial x_{\gamma_3}} \dx,
\end{equation}
where $w \in V_h^3 = V_h^4$ is a given approximation of the velocity,
typically obtained from the previous iteration in an iterative method
for the nonlinear Navier--Stokes equations.  To write the element
tensor for~(\ref{eq:stabilization}) in the canonical
form~(\ref{eq:canonical}), we expand $w$ in the nodal basis for
$\mathcal{P}_K^3 = \mathcal{P}_K^4$ and note that
\begin{equation}
  A^K_i
  =
  \sum_{\gamma_1,\gamma_2,\gamma_3=1}^d
  \sum_{\gamma_4=1}^{n_K^3}
  \sum_{\gamma_5=1}^{n_K^4}
  \int_K
  \frac{\partial \phi^{K,1}_{i_1}[\gamma_1]}{\partial x_{\gamma_2}}
  \frac{\partial \phi^{K,2}_{i_2}[\gamma_1]}{\partial x_{\gamma_3}}
  w_{\gamma_4}^K \phi^{K,3}_{\gamma_4}[\gamma_2]
  w_{\gamma_5}^K \phi^{K,4}_{\gamma_5}[\gamma_3]
  \dx,
\end{equation}
We may then write the element tensor~$A^K$ for the bilinear
form~(\ref{eq:stabilization}) in the canonical
form~(\ref{eq:canonical}), with
\begin{equation}
  \begin{split}
    m &= 4, \\
    \mathcal{C}      &= [1,d]^3 \times [1,n_K^3] \times [1,n_K^4], \\
    c(\gamma)        &= (1, 1, w_{\gamma_4}^K, w_{\gamma_5}^K), \\
    \iota(i, \gamma) &= (i_1, i_2, \gamma_4, \gamma_5), \\
    \kappa(\gamma)   &= (\gamma_1,\gamma_1,\gamma_2,\gamma_3), \\
    \delta(\gamma)   &= (\gamma_2,\gamma_3,\emptyset,\emptyset),
\end{split}
\end{equation}
where $\emptyset$ denotes an empty derivative multiindex (no
differentiation).

In \cite{logg:submit:toms-ffc-monomial-2006}, it is proved that any
element tensor~$A^K$ that can be expressed in the general canonical
form~(\ref{eq:canonical}), can be represented as a tensor
contraction~$A^K = A^0 : G_K$ of a reference tensor~$A^0$ independent
of $K$ and a geometry tensor~$G_K$. A similar result is also presented
in~\cite{logg:article:10} but in less formal notation. As noted
above, element tensors that can be expressed in the general canonical
form correspond to multilinear forms that can be expressed as
integrals over~$\Omega$ of linear combinations of products of basis
functions and their derivatives. The representation theorem reads as
follows.

% Exercise: Do the calculation for the stabilization term

\newpage
\begin{theorem}[Representation theorem] \label{th:representation}
  If $F_K$ is a given affine mapping from a reference cell $K_0$ to
  a cell $K$ and $\{\mathcal{P}_K^j\}_{j=1}^m$ is a given set of
  discrete function spaces on $K$, each generated by a discrete
  function space $\mathcal{P}_0^j$ on the reference cell~$K_0$ through
  the affine mapping, that is, for each $\phi \in \mathcal{P}_K^j$
  there is some $\Phi \in \mathcal{P}_0^j$ such that $\Phi = \phi
  \circ F_K$, then the element tensor~\emph{(\ref{eq:canonical})} may
  be represented as the tensor contraction of a \emph{reference
  tensor} $A^0$ and a \emph{geometry tensor} $G_K$,
  \begin{equation} \label{eq:representation}
    A^K = A^0 : G_K,
  \end{equation}
  that is,
  \begin{equation}
    A^K_i = \sum_{\alpha\in\mathcal{A}} A^0_{i\alpha} G_K^{\alpha}
    \quad \forall i \in \mathcal{I}_K,
  \end{equation}
  where the reference tensor $A^0$ is independent of $K$. In
  particular, the reference tensor $A^0$ is given by
  \begin{equation}
    A^0_{i\alpha}
    =
    \sum_{\beta\in\mathcal{B}}
    \int_{K_0}
    \prod_{j=1}^m
    D_X^{\delta'_j(\alpha,\beta)}
    \Phi^j_{\iota_j(i,\alpha,\beta)}[\kappa_j(\alpha,\beta)]
    \dX,
  \end{equation}
  and the geometry tensor $G_K$ is the outer product of the
  coefficients of any weight functions with a tensor that depends only
  on the Jacobian $F_K'$,
  \begin{equation}
    G_K^{\alpha}
    =
    \prod_{j=1}^m
    c_{j}(\alpha) \,
    \det F_K'
    \sum_{\beta\in\mathcal{B'}}
    \prod_{j'=1}^m
    \prod_{k=1}^{|\delta_{j'}(\alpha,\beta)|}
    \frac{\partial X_{\delta'_{j'k}(\alpha,\beta)}}{\partial x_{\delta_{j'k}(\alpha,\beta)}
},
  \end{equation}
  for some appropriate index sets $\mathcal{A}$, $\mathcal{B}$ and
  $\mathcal{B}'$. We refer to the index set~$\mathcal{I}_K$ as
  the set of \emph{primary indices}, the index set~$\mathcal{A}$ as
  the set of \emph{secondary indices}, and to the index sets~$\mathcal{B}$ and
  $\mathcal{B}'$ as sets of \emph{auxiliary indices}.
\end{theorem}

The ranks of the tensors $A^0$ and $G_K$ are determined by the
properties of the multilinear form~$a$, such as the number of
coefficients and derivatives. Since the rank of the element tensor
$A^K$ is equal to the arity~$r$ of the multilinear form~$a$, the rank
of the reference tensor~$A^0$ must be $|i\alpha| = r + |\alpha|$, where
$|\alpha|$ is the rank of the geometry tensor. For the examples
presented above, we have $|i\alpha| = 4$ and $|\alpha| = 2$ in the case of
Poisson's equation and $|i\alpha| = 8$ and $|\alpha| = 6$ for the
Navier--Stokes stabilization term.

The proof of Theorem~\ref{th:representation} is constructive and gives
an algorithm for computing the
representation~(\ref{eq:representation}). A number of concrete examples
with explicit formulas for the reference and geometry tensors are
given in Tables~\ref{tab:testcase1}--\ref{tab:testcase4}.
We return to these test cases below
in Section~\ref{sec:ffc}, when we discuss the implementation of
Theorem~\ref{th:representation} in the form compiler FFC and present
benchmark results for the test cases.

We remark that in general, a multilinear form will correspond to a sum
of tensor contractions, rather than a single tensor contraction as
in~(\ref{eq:representation}), that is,
\begin{equation}
  A^K = \sum_k A^{0,k} : G_{K,k}.
\end{equation}
One such example is the computation of the element tensor for the
convection--reaction problem $-\Delta u + u = f$, which may be computed
as the sum of a tensor contraction of a rank~four reference
tensor~$A^{0,1}$ with a rank~two geometry tensor~$G_{K,1}$ and a
rank~two reference tensor~$A^{0,2}$ with a rank~zero geometry
tensor~$G_{K,2}$.

\linespread{2.0}
\begin{table}[htbp]
  \begin{center}
    \begin{tabular}{|rcl|c|}
      \hline
      $a(v,U)$ &$=$& $ \int_{\Omega} v \, U \dx$ & rank \\
      \hline
      \hline
      $A^0_{i\alpha}$ &$=$& $\int_{K_0} \Phi_{i_1}^{1} \Phi_{i_2}^{2} \dX$
      & $|i\alpha| = 2$ \\
      \hline
      $G_K^{\alpha}$ &$=$& $\det F_K'$
      & $|\alpha| = 0$ \\
      \hline
    \end{tabular}
    \linespread{1.0}
    \caption{The tensor contraction representation $A^K = A^0 : G_K$ of the
      element tensor~$A^K$ for the bilinear form associated with a
      mass matrix (test case~1).}
    \label{tab:testcase1}
  \end{center}
\end{table}

\linespread{2.0}
\begin{table}[htbp]
  \begin{center}
    \begin{tabular}{|rcl|c|}
      \hline
      $a(v,U)$ &$=$& $\int_{\Omega} \nabla v \cdot \nabla U \dx$
      & rank \\
      \hline
      \hline
      $A^0_{i\alpha}$ &$=$&
      $\int_{K_0}
      \frac{\partial \Phi^1_{i_1}}{\partial X_{\alpha_1}}
      \frac{\partial \Phi^2_{i_2}}{\partial X_{\alpha_2}} \dX$
      & $|i\alpha| = 4$ \\
      \hline
      $G_K^{\alpha}$ &$=$&
      $\det F_K' \sum_{\beta=1}^d
      \frac{\partial X_{\alpha_1}}{\partial x_{\beta}}
      \frac{\partial X_{\alpha_2}}{\partial x_{\beta}}$
      & $|\alpha| = 2$ \\
      \hline
    \end{tabular}
    \linespread{1.0}
    \caption{The tensor contraction representation $A^K = A^0 : G_K$ of the
      element tensor~$A^K$ for the bilinear form associated with
      Poisson's equation (test case~2).}
    \label{tab:testcase2}
  \end{center}
\end{table}

\linespread{2.0}
\begin{table}[htbp]
  \begin{center}
    \begin{tabular}{|rcl|c|}
      \hline
      $a(v,U)$ &$=$& $\int_{\Omega} v \cdot (w \cdot \nabla) U \dx$
      & rank \\
      \hline
      \hline
      $A^0_{i\alpha}$ &$=$&
      $\sum_{\beta=1}^d
      \int_{K_0}
      \Phi^1_{i_1}[\beta]
      \frac{\partial \Phi^2_{i_2}[\beta]}{\partial X_{\alpha_3}}
      \Phi^3_{\alpha_1}[\alpha_2]
      \dX$
      & $|i\alpha| = 5$ \\
      \hline
      $G_K^{\alpha}$ &$=$&
      $w^K_{\alpha_1} \det F_K'
      \frac{\partial X_{\alpha_3}}{\partial x_{\alpha_2}}$
      & $|\alpha| = 3$ \\
      \hline
    \end{tabular}
    \linespread{1.0}
    \caption{The tensor contraction representation $A^K = A^0 : G_K$ of the
      element tensor~$A^K$ for the bilinear form associated with a
      linearization of the nonlinear term $u\cdot \nabla u$ in the
      incompressible Navier--Stokes equations (test case~3).}
    \label{tab:testcase3}
  \end{center}
\end{table}

\linespread{2.0}
\begin{table}[htbp]
  \begin{center}
    \begin{tabular}{|rcl|c|}
      \hline
      $a(v,U)$ &$=$& $\int_{\Omega} \epsilon(v) : \epsilon(U) \dx$
      & rank \\
      \hline
      \hline
      $A^0_{i\alpha}$ &$=$&
      $\sum_{\beta=1}^d
      \int_{K_0}
      \frac{\partial \Phi^1_{i_1}[\beta]}{\partial X_{\alpha_1}}
      \frac{\partial \Phi^2_{i_2}[\beta]}{\partial X_{\alpha_2}}
      \dX$
      & $|i\alpha| = 4$ \\
      \hline
      $G_K^{\alpha}$ &$=$&
      $\frac{1}{2}
      \det F_K'
      \sum_{\beta=1}^d
      \frac{\partial X_{\alpha_1}}{\partial x_{\beta}}
      \frac{\partial X_{\alpha_2}}{\partial x_{\beta}}$
      & $|\alpha| = 2$ \\
      \hline
    \end{tabular}
    \linespread{1.0}
    \caption{The tensor contraction representation $A^K = A^0 : G_K$ of the
      element tensor~$A^K$ for the bilinear form
      $\int_{\Omega} \epsilon(v) : \epsilon(U) \dx =
      \int_{\Omega} \frac{1}{4}
      (\nabla v + (\nabla v)^{\top}) : (\nabla U + (\nabla U)^{\top}) \dx$
      associated with the strain-strain term of linear
      elasticity (test case~4).
      Note that the product expands into four terms which can be grouped
     in pairs of two. The representation is given only for the first
      of these two terms.}
    \label{tab:testcase4}
  \end{center}
\end{table}

\linespread{1.0}

\subsection{Extension to Non-Affine Mappings}
\label{sec:extension}

The tensor contraction representation~(\ref{eq:representation}) of
Theorem~\ref{th:representation} assumes that the mapping~$F_K$ from
the reference cell is affine, allowing the transforms $\partial X /
\partial x$ and the determinant to be pulled out of the integral.  To
see how to extend this result to the case when the mapping~$F_K$ is
non-affine, such as in the case of an isoparametric mapping for a
higher-order element used to map the reference cell to a curvilinear
cell on the boundary of~$\Omega$, we consider again the computation of
the element tensor~$A^K$ for Poisson's equation.  As in
Section~\ref{sec:quadrature}, we use quadrature to evaluate the
integral, but take advantage of the fact that the discrete function
spaces $\mathcal{P}_K^1$ and $\mathcal{P}_K^2$ on $K$ may be generated
from a pair of reference finite elements as discussed in
Section~\ref{sec:femfunctions}. We have
\begin{equation} \label{eq:poisson,quadrature,tensor}
  \begin{split}
  A^K_i
  &=
  \int_K
  \nabla \phi_{i_1}^{K,1} \cdot
  \nabla \phi_{i_2}^{K,2} \dx
  =
  \int_K
  \sum_{\beta=1}^d
  \frac{\partial \phi_{i_1}^{K,1}}{\partial x_{\beta}}
  \frac{\partial \phi_{i_2}^{K,2}}{\partial x_{\beta}} \dx \\
  &=
  \sum_{\alpha_1=1}^d
  \sum_{\alpha_2=1}^d
  \sum_{\beta=1}^d\
  \int_{K_0}
  \frac{\partial X_{\alpha_1}}{\partial x_{\beta}}
  \frac{\partial X_{\alpha_2}}{\partial x_{\beta}}
  \frac{\partial \Phi_{i_1}^{1}}{\partial X_{\alpha_1}}
  \frac{\partial \Phi_{i_2}^{2}}{\partial X_{\alpha_2}} \det F_K' \dX \\
  &\approx
  \sum_{\alpha_1=1}^d
  \sum_{\alpha_2=1}^d
  \sum_{\alpha_3=1}^{N_q}
  w_{\alpha_3}
  \frac{\partial \Phi_{i_1}^{1}}{\partial X_{\alpha_1}}(X_{\alpha_3})
  \frac{\partial \Phi_{i_2}^{2}}{\partial X_{\alpha_2}}(X_{\alpha_3})
  \sum_{\beta=1}^d
  \frac{\partial X_{\alpha_1}}{\partial x_{\beta}}(X_{\alpha_3})
  \frac{\partial X_{\alpha_2}}{\partial x_{\beta}}(X_{\alpha_3})
  \det F_K'(X_{\alpha_3}).
  \end{split}
\end{equation}
As before, we thus obtain a representation of the form
\begin{equation}
  A^K = A^0 : G_K,
\end{equation}
where the reference tensor~$A^0$ is now given by
\begin{equation} \label{eq:quadrature,A}
  A^0_{i\alpha}
  = w_{\alpha_3}
  \frac{\partial \Phi_{i_1}^{1}}{\partial X_{\alpha_1}}(X_{\alpha_3})
  \frac{\partial \Phi_{i_2}^{2}}{\partial X_{\alpha_2}}(X_{\alpha_3}),
\end{equation}
and the geometry tensor~$G_K$ is given by
\begin{equation} \label{eq:quadrature,G}
  G_K^{\alpha}
  =
  \det F_K'(X_{\alpha_3})
  \sum_{\beta=1}^d
  \frac{\partial X_{\alpha_1}}{\partial x_{\beta}}(X_{\alpha_3})
  \frac{\partial X_{\alpha_2}}{\partial x_{\beta}}(X_{\alpha_3}).
\end{equation}
We thus note that a (different) tensor contraction representation
of the element tensor~$A^K$ is possible even if the mapping~$F_K$ is
non-affine. One may also prove a representation theorem similar to
Theorem~\ref{th:representation} for non-affine mappings.

% Exercise: Find the representations for the other three test cases

Comparing the representation
(\ref{eq:quadrature,A})--(\ref{eq:quadrature,G}) with the affine
representation~(\ref{eq:poisson,AandG}), we note that the ranks of
both $A^0$ and $G_K$ have increased by one. As before, we may
precompute the reference tensor $A^0$ but the number of multiply--add
pairs to compute the element tensor~$A^K$ increase by a factor $N_q$
from $n_0^2 d^2$ to $N_q n_0^2 d^2$ (if again we ignore the cost of
computing the geometry tensor).

We also note that the cost has increased by a factor $d$ compared to
the cost of a direct application of quadrature as described in
Section~\ref{sec:quadrature}. However, by expressing the element
tensor~$A^K$ as a tensor contraction, the evaluation of the element
tensor is more readily optimized than if expressed as a triply nested
loop over quadrature points and basis functions as in
Table~\ref{tab:poisson,deal} and Table~\ref{tab:poisson,diffpack}.

As demonstrated below in Section~\ref{sec:optimizations}, it may in
some cases be possible to take advantage of special structures such as
dependencies between different entries in the tensor~$A^0$ to
significantly reduce the operation count. Another more straightforward
approach is to use an optimized library routine such as
a BLAS call to compute the tensor contraction as we shall see below in
Section~\ref{sec:matrixvector}.

\subsection{A Language for Multilinear Forms}

To automate the process of evaluating the element tensor~$A^K$, we
must create a system that takes as input a multilinear form~$a$ and
automatically computes the corresponding element tensor~$A^K$. We do
this by defining a \emph{language} for multilinear forms and
automatically translating any given string in the language to the
canonical form~(\ref{eq:canonical}). From the canonical
form, we may then compute the element tensor $A^K$
by the tensor contraction $A^K = A^0 : G_K$.

When designing such a language for multilinear forms, we have two
things in mind. First, the multilinear forms specified in the language
should be ``close'' to the corresponding mathematical notation (taking
into consideration the obvious limitations of specifying the form as a
string in the ASCII character set). Second, it should be
straightforward to translate a multilinear form specified in the
language to the canonical form~(\ref{eq:canonical}).

A language may be specified formally by defining a \emph{formal
grammar} that generates the language. The grammar specifies a set of
rewrite rules and all strings in the language can be generated by
repeatedly applying the rewrite rules. Thus, one may specify a
language for multilinear forms by defining a suitable grammar
(such as a standard EBNF grammar~\cite{ISO96}), with
basis functions and multiindices as the terminal symbols. One could
then use an automating tool (a compiler-compiler) to create a compiler
for multilinear forms.

However, since a closed canonical form is available for the set of
possible multilinear forms, we will take a more explicit approach.
We fix a small set of operations, allowing only multilinear forms
that have a corresponding canonical form~(\ref{eq:canonical}) to
be expressed through these operations, and observe how the canonical
form transforms under these operations.

\subsubsection{An algebra for multilinear forms}

Consider the set of local finite element spaces
$\{\mathcal{P}_K^j\}_{j=1}^m$ on a cell $K$ corresponding to a set of
global finite element spaces $\{V_h^j\}_{j=1}^m$. The set of local
basis functions $\{\phi_i^{K,j}\}_{i,j=1}^{n_K^j,m}$ span a vector
space $\overline{\mathcal{P}}_K$ and each
function $v$ in this vector space may be expressed as a linear
combination of the basis functions, that is, the set of functions
$\overline{\mathcal{P}}_K$ may be generated from the basis functions
through addition $v + w$ and multiplication with scalars $\alpha
v$. Since $v - w = v + (-1) w$ and $v/\alpha = (1/\alpha) v$, we can
also easily equip the vector space with subtraction and division by
scalars. Informally, we may thus write
\begin{equation}
  \overline{\mathcal{P}}_K = \left\{ v : v = \sum c_{(\cdot)} \phi^K_{(\cdot)} \right\}.
\end{equation}

We next equip our vector space $\overline{\mathcal{P}}_K$ with multiplication between
elements of the vector space. We thus obtain an \emph{algebra} (a
vector space with multiplication) of linear combinations of products
of basis functions. Finally, we extend our algebra $\overline{\mathcal{P}}_K$ by
differentiation $\partial / \partial x_i$ with respect to the coordinate
directions on~$K$, to obtain
\begin{equation}
  \overline{\mathcal{P}}_K = \left\{ v : v = \sum
  c_{(\cdot)} \prod
  \frac{\partial^{|(\cdot)|} \phi^K_{(\cdot)}}
       {\partial x_{(\cdot)}} \right\},
\end{equation}
where $(\cdot)$ represents some multiindex.

To summarize, $\overline{\mathcal{P}}_K$ is the algebra of linear
combinations of products of basis functions or derivatives of basis
functions that is generated from the set of basis functions through
addition ($+$), subtraction ($-$), multiplication $(\cdot)$, including
multiplication with scalars, division by scalars $(/)$, and
differentiation $\partial / \partial x_i$. We note that the algebra is
closed under these operations, that is, applying any of the operators
to an element $v \in \overline{\mathcal{P}}_K$ or a pair of elements
$v, w \in \overline{\mathcal{P}}_K$ yields a member of
$\overline{\mathcal{P}}_K$.

If the basis functions are vector-valued (or tensor-valued), the
algebra is instead generated from the set of scalar components of the
basis functions. Furthermore, we may introduce linear algebra
operators, such as inner products and matrix--vector products, and
differential operators, such as the gradient, the divergence and
rotation, by expressing these compound operators in terms of the basic
operators (addition, subtraction, multiplication and differentiation).

We now note that the algebra $\overline{\mathcal{P}}_K$ corresponds
precisely to the canonical form~(\ref{eq:canonical}) in that the
element tensor~$A^K$ for any multilinear form on $K$ that can be
expressed as an integral over $K$ of an element $v \in
\overline{\mathcal{P}}_K$ has an immediate representation as a sum of
element tensors of the canonical form $(\ref{eq:canonical})$.
We demonstrate this below.

\subsubsection{Examples}

As an example, consider the bilinear form
\begin{equation} \label{eq:mass,again}
  a(v, U) = \int_{\Omega} v \, U \dx,
\end{equation}
with corresponding element tensor canonical form
\begin{equation}
  A^K_i = \int_K \phi^{K,1}_{i_1} \phi^{K,2}_{i_2} \dx.
\end{equation}
If we now let $v = \phi^{K,1}_{i_1} \in \overline{\mathcal{P}}_K$ and
$U = \phi^{K,2}_{i_2} \in \overline{\mathcal{P}}_K$, we note that $v
\, U \in \overline{\mathcal{P}}_K$ and we may thus express the element
tensor as an integral over~$K$ of an element
in~$\overline{\mathcal{P}}_K$,
\begin{equation}
  A^K_i = \int_K v \, U \dx,
\end{equation}
which is close to the notation of~(\ref{eq:mass,again}).
As another example, consider the bilinear form
\begin{equation} \label{eq:reactiondiffusion}
  a(v, U) = \int_{\Omega} \nabla v \cdot \nabla U + v \, U \dx,
\end{equation}
with corresponding element tensor canonical form\footnotemark{}
\begin{equation}
  A^K_i =
  \sum_{\gamma=1}^d
  \int_K
  \frac{\partial \phi_{i_1}^{K,1}}{\partial x_{\gamma}}
  \frac{\partial \phi_{i_2}^{K,2}}{\partial x_{\gamma}}
  \dx
  +
  \int_K \phi^{K,1}_{i_1} \phi^{K,2}_{i_2} \dx.
\end{equation}
As before, we let $v = \phi^{K,1}_{i_1} \in \overline{\mathcal{P}}_K$
and $U = \phi^{K,2}_{i_2} \in \overline{\mathcal{P}}_K$ and note that
$\nabla v \cdot \nabla U + v \, U \in \overline{\mathcal{P}}_K$. It thus
follows that the element tensor~$A^K$ for the bilinear
form~(\ref{eq:reactiondiffusion}) may be expressed as an integral over~$K$
of an element in~$\overline{\mathcal{P}}_K$,
\begin{equation}
  A^K_i = \int_K \nabla v \cdot \nabla U + v \, U \dx,
\end{equation}
which is close to the notation of~(\ref{eq:reactiondiffusion}).
Thus, by a suitable definition of $v$ and $U$ as local basis functions
on~$K$, the canonical form~(\ref{eq:canonical}) for the element tensor
of a given multilinear form may be expressed in a notation that is
close to the notation for the multilinear form itself.

\footnotetext{To be precise, the element tensor is the sum of two
  element tensors, each written in the canonical
  form~(\ref{eq:canonical}) with a suitable definition of
  multiindices~$\iota$, $\kappa$ and $\delta$.}

\subsubsection{Implementation by operator-overloading}
\label{sec:operatoroverloading}

It is now straightforward to implement the
algebra~$\overline{\mathcal{P}}_K$ in any object-oriented language
with support for operator overloading, such as Python or C++. We first
implement a class \texttt{BasisFunction}, representing (derivatives
of) basis functions of some given finite element space. Each
\texttt{BasisFunction} is associated with a particular finite element
space and different \texttt{BasisFunction}s may be associated with
different finite element spaces. Products of scalars and (derivatives
of) basis functions are represented by the class \texttt{Product},
which may be implemented as a list of \texttt{BasisFunction}s.  Sums
of such products are represented by the class \texttt{Sum}, which may
be implemented as a list of \texttt{Product}s. We then define an
operator for differentiation of basis functions and overload the
operators addition, subtraction and multiplication, to generate the
algebra of \texttt{BasisFunction}s, \texttt{Product}s and
\texttt{Sum}s, and note that any combination of such operators and
objects ultimately yields an object of class~\texttt{Sum}. In
particular, any object of class \texttt{BasisFunction} or
\texttt{Product} may be cast to an object of class
\texttt{Sum}.

By associating with each object one or more \emph{indices},
implemented by a class \texttt{Index}, an object of class
\texttt{Product} automatically represents a tensor expressed in the
canonical form~(\ref{eq:canonical}). Finally, we note that we may
introduce compound operators such as \texttt{grad}, \texttt{div},
\texttt{rot}, \texttt{dot} etc. by expressing these operators in terms
of the basic operators.

Thus, if \texttt{v} and \texttt{U} are objects of
class~\texttt{BasisFunction}, the integrand of the bilinear
form~(\ref{eq:reactiondiffusion}) may be given as the string
\begin{equation}
  \mbox{\texttt{dot(grad(v), grad(U)) + v*U}}.
\end{equation}
In Table~\ref{tab:poisson,fenics} we saw a similar example of how the
bilinear form for Poisson's equation is specified in the language of
the FEniCS Form Compiler FFC. Further examples will be given below in
Section~\ref{sec:ffc} and Section~\ref{sec:examples}.

%------------------------------------------------------------------------------
\section{AUTOMATING THE ASSEMBLY OF THE DISCRETE SYSTEM}
\label{sec:automation,assembly}

In Section~\ref{sec:fem}, we reduced the task of automatically
generating the discrete system $F(U) = 0$ for a given nonlinear
variational problem $a(U;v) = L(v)$ to the automatic assembly of  the
tensor~$A$ that represents a given multilinear form~$a$ in a given
finite element basis. By Algorithm~\ref{alg:assembly,2}, this
process may be automated by automating first the computation of the
element tensor~$A^K$, which we discussed in the previous section, and
then automating the addition of the element tensor~$A^K$ into the global
tensor~$A$, which is the topic of the current section.

\subsection{Implementing the Local-to-Global Mapping}

With $\{\iota_K^j\}_{j=1}^{r}$ the local-to-global mappings for
a set of discrete function spaces, $\{V_h^j\}_{j=1}^r$,
we evaluate for each $j$ the local-to-global mapping $\iota_K^j$ on
the set of local node numbers $\{1,2,\ldots,n_K^j\}$, thus obtaining
for each $j$ a tuple
\begin{equation}
  \iota_K^j([1,n_K^j]) = (\iota_K^j(1), \iota_K^j(2), \ldots, \iota_K^j(n_K^j)).
\end{equation}
The entries of the element tensor $A^K$ may then be added to the
global tensor~$A$ by an optimized low-level library
call\footnotemark{} that takes as input the two tensors $A$ and $A^K$
and the set of tuples (arrays) that determine how each dimension of
$A^K$ should be distributed onto the global tensor~$A$. Compare
Figure~\ref{fig:insertion} with the two tuples given by
$(\iota_K^1(1), \iota_K^1(2), \iota_K^1(3))$
and
$(\iota_K^2(1), \iota_K^2(2), \iota_K^2(3))$ respectively.

\footnotetext{If PETSc~\cite{www:PETSc,BalBus04,BalEij97} is used as
  the linear algebra backend, such a library call is available with
  the call \texttt{VecSetValues()} for a rank~one tensor (a vector)
  and \texttt{MatSetValues()} for a rank~two tensor (a matrix).}

Now, to compute the set of tuples~$\{\iota_K^j([1,n_K^j])\}_{j=1}^r$,
we may consider implementing for each $j$ a function that takes as
input the current cell $K$ and returns the corresponding tuple
$\iota_K^j([1,n_K])$. Since the local-to-global mapping
may look very different for different function spaces, in
particular for different degree Lagrange elements, a different
implementation is needed for each different function space. Another
option is to implement a general purpose function that handles a range
of function spaces, but this quickly becomes inefficient. From the example
implementations given in Table~\ref{tab:mapping,P1} and
Table~\ref{tab:mapping,P2} for continuous linear and quadratic
Lagrange finite elements on tetrahedra, it is further clear that if
the local-to-global mappings are implemented individually for each
different function space, the mappings can be implemented
very efficiently, with minimal need for arithmetic or branching.

\begin{table}[htbp]
  \begin{code}
 void nodemap(int nodes[], const Cell& cell, const Mesh& mesh)
 {
   nodes[0] = cell.vertexID(0);
   nodes[1] = cell.vertexID(1);
   nodes[2] = cell.vertexID(2);
   nodes[3] = cell.vertexID(3);
 }
 \end{code}
  \caption{A C++ implementation of the mapping from local to global
    node numbers for continuous linear Lagrange finite elements on
    tetrahedra. One node is associated with each vertex of a local cell
    and the local node number for each of the four nodes is mapped to
    the global number of the associated vertex.}
  \label{tab:mapping,P1}
\end{table}

\begin{table}[htbp]
  \begin{code}
  void nodemap(int nodes[], const Cell& cell, const Mesh& mesh)
  {
    nodes[0] = cell.vertexID(0);
    nodes[1] = cell.vertexID(1);
    nodes[2] = cell.vertexID(2);
    nodes[3] = cell.vertexID(3);
    int offset = mesh.numVertices();
    nodes[4] = offset + cell.edgeID(0);
    nodes[5] = offset + cell.edgeID(1);
    nodes[6] = offset + cell.edgeID(2);
    nodes[7] = offset + cell.edgeID(3);
    nodes[8] = offset + cell.edgeID(4);
    nodes[9] = offset + cell.edgeID(5);
 }
  \end{code}
  \caption{A C++ implementation of the mapping from local to global
    node numbers for continuous quadratic Lagrange finite elements on
    tetrahedra. One node is associated with each vertex and also each
    edge of a local cell. As for linear Lagrange elements, local
    vertex nodes are mapped to the global number of the associated
    vertex, and the remaining six edge nodes are given global numbers
    by adding to the global edge number an offset given
    by the total number of vertices in the mesh.}
  \label{tab:mapping,P2}
\end{table}

\subsection{Generating the Local-to-Global Mapping}

We thus seek a way to automatically generate the code for the
local-to-global mapping from a simple description of the distribution
of nodes on the mesh. As before, we restrict our attention to elements
with nodes given by point evaluation. In that case, each node can be
associated with a geometric entity, such as a vertex, an edge, a face
or a cell. More generally, we may order the geometric entities by
their topological dimension to make the description independent of
dimension-specific notation (compare~\cite{KneKar05}); for a
two-dimensional triangular mesh, we may refer to a (topologically
two-dimensional) triangle as a \emph{cell}, whereas for a
three-dimensional tetrahedral mesh, we would refer to a (topologically
two-dimensional) triangle as a \emph{face}. We may thus for each
topological dimension list the nodes associated with the geometric
entities within that dimension. More specifically, we may list for
each topological dimension and each geometric entity within that
dimension a tuple of nodes associated with that geometric entity.
This approach is used by the FInite element Automatic
Tabulator~FIAT~\cite{www:FIAT,Kir04,Kir06}.

As an example, consider the local-to-global mapping for the linear
tetrahedral element of Table~\ref{tab:mapping,P1}. Each cell has four
nodes, one associated with each vertex. We may then describe the nodes
by specifying for each geometric entity of dimension zero (the
vertices) a tuple containing one local node number, as demonstrated in
Table~\ref{tab:entityids,P1}. Note that we may specify the nodes for a
\emph{discontinuous} Lagrange finite element on a tetrahedron
similarly by associating all for nodes with topological dimension
three, that is, with the cell itself, so that no nodes are shared
between neighboring cells.

As a further illustration, we may describe the nodes for the quadratic
tetrahedral element of Table~\ref{tab:mapping,P2} by associating the
first four nodes with topological dimension zero (vertices) and the
remaining six nodes with topological dimension one (edges), as
demonstrated in Table~\ref{tab:entityids,P2}.

Finally, we present in Table~\ref{tab:entityids,P5} the specification
of the nodes for fifth-degree Lagrange finite elements on
tetrahedra. Since there are now multiple nodes associated with some
entities, the ordering of nodes becomes important. In particular, two
neighboring tetrahedra sharing a common edge (face) must agree on the
global node numbering of edge (face) nodes. This can be accomplished by
checking the orientation of geometric entities with respect to some
given convention.\footnotemark{} For each edge, there are two possible
orientations and for each face of a tetrahedron, there are six
possible orientations. In Table~\ref{tab:mapping,P5}, we present the
local-to-global mapping for continuous fifth-degree Lagrange finite
elements, generated automatically from the description
of~Table~\ref{tab:entityids,P5} by the FEniCS Form
Compiler~FFC~\cite{logg:www:04,logg:article:10,logg:submit:toms-ffc-monomial-2006,logg:manual:02}.

\footnotetext{For an example of such a convention,
  see~\cite{logg:manual:01} or~\cite{logg:manual:02}.}

\begin{table}[htbp]
  \begin{center}
    \begin{tabular}{|l|l|}
      \hline
      $d = 0$ & $(1)$ -- $(2)$ -- $(3)$ -- $(4)$ \\
      \hline
    \end{tabular}
    \caption{Specifying the nodes for continuous linear Lagrange
      finite elements on tetrahedra.}
    \label{tab:entityids,P1}
  \end{center}
\end{table}

\begin{table}[htbp]
  \begin{center}
    \linespread{1.5}
    \begin{tabular}{|l|l|}
      \hline
      $d = 0$ & $(1)$ -- $(2)$ -- $(3)$ -- $(4)$ \\
      \hline
      $d = 1$ & $(5)$ -- $(6)$ -- $(7)$ -- $(8)$ -- $(9)$ -- $(10)$ \\
      \hline
    \end{tabular}
    \linespread{1.0}
    \caption{Specifying the nodes for continuous quadratic
      Lagrange finite elements on tetrahedra.}
    \label{tab:entityids,P2}
  \end{center}
\end{table}

We may thus think of the local-to-global mapping as a function that
takes as input the current cell~$K$ (\texttt{cell}) together with the
mesh~$\mathcal{T}$ (\texttt{mesh}) and generates a tuple
(\texttt{nodes}) that maps the local node numbers on~$K$ to global
node numbers. For finite elements with nodes given by point
evaluation, we may similarly generate a function that interpolates any
given function to the current cell~$K$ by evaluating it at the nodes.

\begin{table}[htbp]
  \begin{center}
    \begin{tabular}{|l|l|}
      \hline
      $d = 0$ & $(1)$ -- $(2)$ -- $(3)$ -- $(4)$ \\
      \hline
      $d = 1$ & $(5,6,7,8)$ -- $(9,10,11,12)$ -- $(13,14,15,16)$ -- \\
              & $(17,18,19,20)$ -- $(21,22,23,24)$ -- $(25,26,27,28)$
      \\
      \hline
      $d = 2$ & $(29,30,31,32,33,34)$ -- $(35,36,37,38,39,40)$ -- \\
              & $(41,42,43,44,45,46)$ -- $(47,48,49,50,51,52)$ \\
      \hline
      $d = 3$ & $(53,54,55,56)$ \\
      \hline
    \end{tabular}
    \caption{Specifying the nodes for continuous fifth-degree
      Lagrange finite elements on tetrahedra.}
    \label{tab:entityids,P5}
  \end{center}
\end{table}

\begin{table}[htbp]
  \small
  \begin{code}
 void nodemap(int nodes[], const Cell& cell, const Mesh& mesh)
 {
   static unsigned int edge_reordering[2][4] = {{0, 1, 2, 3}, {3, 2, 1, 0}};
   static unsigned int face_reordering[6][6] = {{0, 1, 2, 3, 4, 5},
                                                {0, 3, 5, 1, 4, 2},
                                                {5, 3, 0, 4, 1, 2},
                                                {2, 1, 0, 4, 3, 5},
                                                {2, 4, 5, 1, 3, 0},
                                                {5, 4, 2, 3, 1, 0}};
   nodes[0] = cell.vertexID(0);
   nodes[1] = cell.vertexID(1);
   nodes[2] = cell.vertexID(2);
   nodes[3] = cell.vertexID(3);
   int alignment = cell.edgeAlignment(0);
   int offset = mesh.numVertices();
   nodes[4] = offset + 4*cell.edgeID(0) + edge_reordering[alignment][0];
   nodes[5] = offset + 4*cell.edgeID(0) + edge_reordering[alignment][1];
   nodes[6] = offset + 4*cell.edgeID(0) + edge_reordering[alignment][2];
   nodes[7] = offset + 4*cell.edgeID(0) + edge_reordering[alignment][3];
   alignment = cell.edgeAlignment(1);
   nodes[8] = offset + 4*cell.edgeID(1) + edge_reordering[alignment][0];
   nodes[9] = offset + 4*cell.edgeID(1) + edge_reordering[alignment][1];
   nodes[10] = offset + 4*cell.edgeID(1) + edge_reordering[alignment][2];
   nodes[11] = offset + 4*cell.edgeID(1) + edge_reordering[alignment][3];
   ...
   alignment = cell.faceAlignment(0);
   offset = offset + 4*mesh.numEdges();
   nodes[28] = offset + 6*cell.faceID(0) + face_reordering[alignment][0];
   nodes[29] = offset + 6*cell.faceID(0) + face_reordering[alignment][1];
   nodes[30] = offset + 6*cell.faceID(0) + face_reordering[alignment][2];
   nodes[31] = offset + 6*cell.faceID(0) + face_reordering[alignment][3];
   nodes[32] = offset + 6*cell.faceID(0) + face_reordering[alignment][4];
   nodes[33] = offset + 6*cell.faceID(0) + face_reordering[alignment][5];
   ...
   offset = offset + 6*mesh.numFaces();
   nodes[52] = offset + 4*cell.id() + 0;
   nodes[53] = offset + 4*cell.id() + 1;
   nodes[54] = offset + 4*cell.id() + 2;
   nodes[55] = offset + 4*cell.id() + 3;
 }
  \end{code}
  \normalsize
  \caption{A C++ implementation (excerpt) of the mapping from local to
    global node numbers for continuous fifth-degree Lagrange finite
    elements on tetrahedra. One node is associated with each vertex,
    four nodes with each edge, six nodes with each face and four nodes
    with the tetrahedron itself.}
  \label{tab:mapping,P5}
\end{table}

%------------------------------------------------------------------------------
\section{OPTIMIZATIONS}
\label{sec:optimizations}

As we saw in Section~\ref{sec:automation,forms}, the (affine) tensor
contraction representation of the element tensor for Poisson's
equation may significantly reduce the operation count in the
computation of the element tensor. This is true for a wide range of
multilinear forms, in particular test cases~1--4 presented in
Tables~\ref{tab:testcase1}--\ref{tab:testcase4}.

In some cases however, it may be more efficient to compute the element
tensor by quadrature, either using the direct approach of
Section~\ref{sec:quadrature} or by a tensor contraction representation
of the quadrature evaluation as in Section~\ref{sec:extension}. Which
approach is more efficient depends on the multilinear form and the
function spaces on which it is defined. In particular, the relative
efficiency of a quadrature-based approach increases as the number of
coefficients in the multilinear form increases, since then the rank of
the reference tensor increases. On the other hand, the relative
efficiency of the (affine) tensor contraction representation increases
when the polynomial degree of the basis functions and thus the number
of quadrature points increases. See~\cite{logg:article:10} for a more
detailed account.

\subsection{Tensor Contractions as Matrix--Vector Products}
\label{sec:matrixvector}

As demonstrated above, the representation of the element tensor~$A^K$
as a tensor contraction~$A^K = A^0 : G_K$ may be generated
automatically from a given multilinear form. To evaluate the element
tensor~$A^K$, it thus remains to evaluate the tensor contraction. A
simple approach would be to iterate over the entries
$\{A^K_i\}_{i\in\mathcal{I}_K}$ of $A^K$ and for each entry $A^K_i$
compute the value of the entry by summing over the set of secondary
indices as outlined in Algorithm~\ref{alg:contraction,simple}.

\begin{algorithm}
  \begin{tabbing}
    \textbf{for}  {$i \in \mathcal{I}_K$}\\
    \tab $A^K_i = 0$ \\
    \tab \textbf{for} $\alpha \in \mathcal{A}$ \\
    \tab \tab $A^K_i = A^K_i + A^0_{i\alpha} G_K^{\alpha}$ \\
    \tab \textbf{end for} \\
    \textbf{end for}
  \end{tabbing}
  \caption{$A^K$ = ComputeElementTensor()}
  \label{alg:contraction,simple}
\end{algorithm}

Examining Algorithm~\ref{alg:contraction,simple}, we note that by an
appropriate ordering of the entries in $A^K$, $A^0$ and $G_K$, one may
rephrase the tensor contraction as a
matrix--vector product and call an optimized library
routine\footnotemark{} for the computation of the matrix--vector
product.

\footnotetext{Such a library call is available with the standard
level 2 BLAS~\cite{BlaDem02} routine DGEMV, with optimized
implementations provided for different architectures by
ATLAS~\cite{www:ATLAS,WhaDon97,WhaPet01}.}

To see how to write the tensor
contraction as a matrix--vector product,
we let $\{i^j\}_{j=1}^{|\mathcal{I}_K|}$ be an enumeration of the set of
primary multiindices $\mathcal{I}_K$ and let
$\{\alpha^j\}_{j=1}^{|\mathcal{A}|}$ be an enumeration of the set of
secondary multiindices $\mathcal{A}$. As an example, for the
computation of the $6 \times 6$ element tensor for Poisson's equation
with quadratic elements on triangles, we may enumerate the primary and
secondary multiindices by
\begin{equation}
  \begin{split}
    \{i^j\}_{j=1}^{|\mathcal{I}_K|}    &= \{ (1,1), (1,2), \ldots, (1,6), (2,1), \ldots, (6,6) \}, \\
    \{\alpha^j\}_{j=1}^{|\mathcal{A}|} &= \{ (1,1), (1,2), (2,1), (2,2) \}.
  \end{split}
\end{equation}
By similarly enumerating the $36$~entries of the $6 \times 6$ element
tensor~$A^K$ and the four entries of the $2 \times 2$ geometry
tensor~$G_K$, one may define two vectors~$a^K \in
\R^{36}$ and $g_K \in \R^4$ corresponding to the two tensors~$A^K$
and~$G_K$ respectively.

In general, the element tensor~$A^K$ and the geometry tensor~$G_K$ may
thus be flattened to create the corresponding vectors~$a^K
\leftrightarrow A^K$ and $g_K \leftrightarrow G_K$, defined by
\begin{equation}
  \begin{split}
    a^K &= (A^K_{i^1}, A^K_{i^2}, \ldots, A^K_{i^{|\mathcal{I}_K|}})^{\top}, \\
    g_K &= (G_K^{\alpha^1}, G_K^{\alpha^2}, \ldots, G_K^{\alpha^{|\mathcal{A}|}})^{\top}.
  \end{split}
\end{equation}
Similarly, we define the $|\mathcal{I}_K| \times |\mathcal{A}|$ matrix
$\bar{A}^0$ by
\begin{equation}
  \bar{A}^0_{jk} = A^0_{i^j \alpha^k}, \quad
  j = 1,2,\ldots,|\mathcal{I}_K|, \quad k = 1,2,\ldots,|\mathcal{A}|.
\end{equation}
Since now
\begin{equation}
  a^K_j = A^K_{i^j} = \sum_{\alpha \in \mathcal{A}} A^0_{i^j\alpha}
  G_K^{\alpha}
  = \sum_{k=1}^{|\mathcal{A}|} A^0_{i^j\alpha^k} G_K^{\alpha^k}
  = \sum_{k=1}^{|\mathcal{A}|} \bar{A}^0_{jk} (g_K)_k,
\end{equation}
it follows that the tensor contraction $A^K = A^0 : G_K$
corresponds to the matrix--vector product
\begin{equation}
  a^K = \bar{A}^0 g_K.
\end{equation}

As noted earlier, the element tensor $A^K$ may generally be expressed
as a sum of tensor contractions, rather than as a single tensor
contraction, that is,
\begin{equation} \label{eq:representation,sum,again}
  A^K = \sum_k A^{0,k} : G_{K,k}.
\end{equation}
In that case, we may still compute the (flattened) element
tensor~$A^K$ by a single matrix--vector product,
\begin{equation}
  a^K = \sum_k \bar{A}^{0,k} g_{K,k}
  =
  \left[ \bar{A}^{0,1} \,\, \bar{A}^{0,2} \,\, \cdots \right]
  \left[ \begin{array}{c} g_{K,1} \\ \\ g_{K,2} \\ \vdots \end{array}\right]
  = \bar{A}^0 g_K.
\end{equation}

Having thus phrased the general tensor
contraction~(\ref{eq:representation,sum,again}) as a matrix--vector
product, we note that by grouping the cells of the mesh $\mathcal{T}$
into subsets, one may compute the set of element tensors for all cells
in a subset by one matrix--matrix product (corresponding to a level~3
BLAS call) instead of by a sequence of matrix--vector products (each
corresponding to a level~2 BLAS call), which will typically lead to
improved floating-point performance. This is possible since the
(flattened) reference tensor~$\bar{A}^0$ remains constant over the
mesh. Thus, if $\{K_k\}_k \subset \mathcal{T}$ is a subset of the
cells in the mesh, we have
\begin{equation}
  \left[a^{K_1} \,\, a^{K_2} \,\, \cdots \right]
  = \left[\bar{A}^0 g_{K_1} \,\, \bar{A}^0 g_{K_2} \,\, \ldots\right]
  = \bar{A}^0 \left[g_{K_1} \,\, g_{K_2} \,\, \ldots\right].
\end{equation}
The optimal size of each subset is problem and architecture
dependent. Since the geometry tensor may sometimes contain a large
number of entries, the size of the subset may be limited by the
available memory.

\subsection{Finding an Optimized Computation}

Although the techniques discussed in the previous section may often
lead to good floating-point performance, they do not take full
advantage of the fact that the reference tensor is generated
automatically. In~\cite{logg:article:07} and later
in~\cite{logg:article:09}, it was noted that by knowing the size and
structure of the reference tensor at compile-time, one may generate
very efficient code for the computation of the reference tensor.

Letting $g_K \in \R^{|\mathcal{A}|}$ be the vector obtained by
flattening the geometry tensor~$G_K$ as above, we
note that each entry~$A^K_i$ of the element tensor~$A^K$
is given by the inner product
\begin{equation} \label{eq:innerproduct}
  A^K_i = a^0_{i} \cdot g_K,
\end{equation}
where $a^0_{i}$ is the vector defined by
\begin{equation}
  a^0_{i} = (A^0_{i\alpha^1}, A^0_{i\alpha^2}, \ldots,
  A^0_{i\alpha^{|\mathcal{A}|}})^{\top}.
\end{equation}
To optimize the evaluation of the element tensor, we look for
dependencies between the vectors $\{a^0_{i}\}_{i\in\mathcal{I}_K}$ and
use the dependencies to reduce the operation count.
There are many such dependencies to explore. Below, we consider
collinearity and closeness in Hamming distance between pairs of
vectors~$a^0_i$ and~$a^0_{i'}$.

\subsubsection{Collinearity}

We first consider the case when two vectors $a^0_i$ and
$a^0_{i'}$ are collinear, that is,
\begin{equation}
  a^0_{i'} = \alpha a^0_i,
\end{equation}
for some nonzero $\alpha\in\R$. If $a^0_i$ and $a^0_{i'}$ are collinear,
it follows that
\begin{equation}
  A^K_{i'} = a^0_{i'} \cdot g_K = (\alpha a^0_{i}) \cdot g_K =
  \alpha A^K_i.
\end{equation}
We may thus compute the entry $A^K_{i'}$ in a single multiplication,
if the entry $A^K_i$ has already been computed.

\subsubsection{Closeness in Hamming distance}

Another possibility is to look for closeness between pairs of vectors
$a^0_i$ and $a^0_{i'}$ in Hamming
distance (see~\cite{CorLei01}), which is defined
as the number entries in which two vectors differ.
If the Hamming distance between
$a^0_i$ and
$a^0_{i'}$ is $\rho$, then the entry $A^0_{i'}$ may be computed from
the entry $A^0_{i}$ in at most $\rho$ multiply--add pairs. To see this,
we assume that $a^0_i$ and $a^0_{i'}$ differ only in the first $\rho$
entries. It then follows that
\begin{equation}
  A^K_{i'} = a^0_{i'} \cdot g_K =
  a^0_i \cdot g_K + \sum_{k=1}^{\rho} (A^0_{i'\alpha^k} -
  A^0_{i\alpha^k}) G_K^{\alpha^k}
  =
  A^K_i + \sum_{k=1}^{\rho} (A^0_{i'\alpha^k} - A^0_{i\alpha^k}) G_K^{\alpha^k},
\end{equation}
where we note that the vector $(A^0_{i'\alpha^1} - A^0_{i\alpha^1},
A^0_{i'\alpha^2} - A^0_{i\alpha^2}, \ldots, A^0_{i'\alpha^{\rho}} -
A^0_{i\alpha^{\rho}})^{\top}$ may be precomputed at compile-time. We note
that the maximum Hamming distance between $a^0_i$ and $a^0_{i'}$ is
$\rho = |\mathcal{A}|$, that is, the length of the vectors, which is
also the cost for the direct computation of an entry by the inner
product~(\ref{eq:innerproduct}). We also note that if $a^0_i =
a^0_{i'}$ and consequently $A^K_i = A^K_{i'}$, then the Hamming
distance and the cost of obtaining $A^K_{i'}$ from $A^K_i$ are both
zero.

\subsubsection{Complexity-reducing relations}

In~\cite{logg:article:09}, dependencies between pairs of vectors, such
as collinearity and closeness in Hamming distance, that can be used to
reduce the operation count in computing one entry from another, are
referred to as \emph{complexity-reducing relations}. In general, one
may define for any pair of vectors $a^0_i$ and $a^0_{i'}$ the
complexity-reducing relation $\rho(a^0_i, a^0_{i'}) \leq
|\mathcal{A}|$ as the minimum of all complexity complexity reducing
relations found between $a^0_i$ and $a^0_{i'}$. Thus, if we look for
collinearity and closeness in Hamming distance, we may say that
$\rho(a^0_i, a^0_{i'})$ is in general given by the Hamming distance
between $a^0_i$ and $a^0_{i'}$ unless the two vectors are collinear,
in which case $\rho(a^0_i, a^0_{i'}) \leq 1$.

\subsubsection{Finding a minimum spanning tree}

Given the set of vectors $\{a^0_{i}\}_{i\in\mathcal{I}_K}$ and a
complexity-reducing relation $\rho$, the problem is now to find an
optimized computation of the element tensor~$A^K$ by systematically
exploring the complexity-reducing
relation~$\rho$. In~\cite{logg:article:09}, it was found that this
problem has a simple solution. By constructing a weighted undirected
graph $G = (V, E)$ with vertices given by the vectors
$\{a^0_{i}\}_{i\in\mathcal{I}_K}$ and the weight at each edge given by the
value of the complexity-reducing relation~$\rho$ evaluated at the pair
of end-points, one may find an optimized (but not necessarily optimal)
evaluation of the element
tensor by computing the minimum spanning tree\footnotemark{}
$G' = (V, E')$ for the graph~$G$.

\footnotetext{A spanning tree for a graph $G = (V, E)$ is any
  connected acyclic subgraph $(V, E')$ of $(V, E)$, that is, each
  vertex in $V$ is connected to an edge in $E' \subset E$ and there are no
  cycles. The (generally non-unique) \emph{minimum} spanning tree of a
  weighted graph $G$ is a spanning tree $G' =
  (V, E')$ that minimizes the sum of edge weights for $E'$. The
  minimum spanning tree may be computed using standard algorithms
  such as Kruskal's and Prim's algorithms, see~\cite{CorLei01}.}

The minimum spanning tree directly provides an algorithm for the
evaluation of the element tensor~$A^K$. If one first computes the
entry of $A^K$ corresponding to the root vertex of the minimum
spanning tree, which may be done in $|\mathcal{A}|$ multiply--add
pairs, the remaining entries may then be computed by traversing the
tree (following the edges), either breadth-first or depth-first, and
at each vertex computing the corresponding entry of~$A^K$ from the
parent vertex at a cost given by the weight of the connecting edge.
The total cost of computing the element tensor~$A^K$ is thus given by
\begin{equation}
  |\mathcal{A}| + |E'|,
\end{equation}
where $|E'|$ denotes the weight of the minimum spanning tree. As we
shall see, computing the minimum spanning tree may significantly
reduce the operation count, compared to the straightforward approach
of Algorithm~\ref{alg:contraction,simple} for which the operation
count is given by $|\mathcal{I}_K| \, |\mathcal{A}|$.

\subsubsection{A concrete example}

To demonstrate these ideas, we compute the minimum spanning tree for
the computation of the $36$ entries of the $6 \times 6$ element tensor
for Poisson's equation with quadratic elements on triangles and obtain
a reduction in the operation count from from a total of
$|\mathcal{I}_K| \, |\mathcal{A}| = 36 \times 4 = 144$ multiply--add
pairs to less than $17$ multiply--add pairs. Since there are $36$
entries in the element tensor, this means that we are be able to
compute the element tensor in less than one operation per entry
(ignoring the cost of computing the geometry tensor).

As we saw above in Section~\ref{sec:representation}, the rank~four
reference tensor is $A^0$ is given by
\begin{equation} \label{eq:poisson,referencetensor}
  A^0_{i\alpha} =
  \int_{K_0}
  \frac{\partial \Phi^1_{i_1}}{\partial X_{\alpha_1}}
  \frac{\partial \Phi^2_{i_2}}{\partial X_{\alpha_2}} \dX
  \quad
  \forall i \in \mathcal{I}_K \quad
  \forall \alpha \in \mathcal{A},
\end{equation}
where now $\mathcal{I}_K = [1,6]^2$ and $\mathcal{A} = [1,2]^2$. To
compute the $6^2 \times 2^2 = 144$ entries of the reference tensor, we
evaluate the set of integrals~(\ref{eq:poisson,referencetensor})
for the basis defined in~(\ref{eq:P2}). In Table~\ref{tab:A0,P2}, we
give the corresponding set of (scaled) vectors
$\{a^0_i\}_{i\in\mathcal{I}_K}$ displayed as a $6 \times 6$ matrix of
vectors with rows corresponding to the first component~$i_1$ of the
multiindex~$i$ and columns corresponding to the second component~$i_2$
of the multiindex~$i$. Note that the entries in Table~\ref{tab:A0,P2}
have been scaled with a factor $6$ for ease of notation (corresponding
to the bilinear form $a(v, U) = 6 \int_{\Omega} \nabla v \cdot \nabla
U \dx $). Thus, the entries of the reference tensor are given by
$A^0_{1111} = A^0_{1112} = A^0_{1121} = A^0_{1122} = 3/6 = 1/2$,
$A^0_{1211} = 1/6$, $A^0_{1212} = 0$, etc.

\vspace{0.5cm}

\linespread{2.0}
\begin{table}[htbp]
  \begin{center}
    \small
    \newcommand{\toffset}{1.3cm}
    \newcommand{\tspace}{1.97cm}
    \begin{tabular}{p{\toffset}p{\tspace}p{\tspace}p{\tspace}p{\tspace}p{\tspace}p{\tspace}}
      &1 & 2 & 3 & 4 & 5 & 6
    \end{tabular}
    \\
    \begin{tabular}{r}
      1 \\ 2 \\ 3 \\ 4 \\ 5 \\ 6
    \end{tabular}
    \begin{tabular}{|r|r|r|r|r|r|}
      \hline
      $(3, 3, 3, 3)^{\top}$
      &
      $(1, 0, 1, 0)^{\top}$
      &
      $(0, 1, 0, 1)^{\top}$
      &
      $(0, 0, 0, 0)^{\top}$
      &
      $-(0, 4, 0, 4)^{\top}$
      &
      $-(4, 0, 4, 0)^{\top}$
      \\\hline
      $(1, 1, 0, 0)^{\top}$
      &
      $(3, 0, 0, 0)^{\top}$
      &
      $-(0, 1, 0, 0)^{\top}$,
      &
      $(0, 4, 0, 0)^{\top}$
      &
      $(0, 0, 0, 0)^{\top}$
      &
      $-(4, 4, 0, 0)^{\top}$
      \\\hline
      $(0, 0, 1, 1)^{\top}$
      &
      $-(0, 0, 1, 0)^{\top}$
      &
      $(0, 0, 0, 3)^{\top}$
      &
      $(0, 0, 4, 0)^{\top}$
      &
      $-(0, 0, 4, 4)^{\top}$
      &
      $(0, 0, 0, 0)^{\top}$
      \\\hline
      $(0, 0, 0, 0)^{\top}$
      &
      $(0, 0, 4, 0)^{\top}$
      &
      $(0, 4, 0, 0)^{\top}$
      &
      $(8, 4, 4, 8)^{\top}$
      &
      $-(8, 4, 4, 0)^{\top}$
      &
      $-(0, 4, 4, 8)^{\top}$
      \\\hline
      $-(0, 0, 4, 4)^{\top}$
      &
      $(0, 0, 0, 0)^{\top}$
      &
      $-(0, 4, 0, 4)^{\top}$
      &
      $-(8, 4, 4, 0)^{\top}$
      &
      $(8, 4, 4, 8)^{\top}$
      &
      $(0, 4, 4, 0)^{\top}$
      \\\hline
      $-(4, 4, 0, 0)^{\top}$
      &
      $-(4, 0, 4, 0)^{\top}$
      &
      $(0, 0, 0, 0)^{\top}$
      &
      $-(0, 4, 4, 8)^{\top}$
      &
      $(0, 4, 4, 0)^{\top}$
      &
      $(8, 4, 4, 8)^{\top}$
      \\\hline
    \end{tabular}
    \linespread{1.0}
    \caption{The $6\times 6\times 2\times 2$ reference tensor~$A^0$
    for Poisson's equation with quadratic elements on triangles,
    displayed here as the set of vectors~$\{a^0_i\}_{i \in \mathcal{I}_K}$.}
    \label{tab:A0,P2}
  \end{center}
\end{table}
\linespread{1.0}

\vspace{0.5cm}

Before proceeding to compute the minimum spanning tree for the $36$
vectors in Table~\ref{tab:A0,P2}, we note that the element tensor $A^K$
for Poisson's equation is symmetric, and as a consequence we only need
to compute $21$ of the $36$ entries of the element tensor. The
remaining $15$ entries are given by symmetry. Furthermore, since the
geometry tensor $G_K$ is symmetric (see Table~\ref{tab:testcase2}), it
follows that
\begin{equation}
  \begin{split}
    A^K_i &= a^0_i \cdot g_K =
    A^0_{i11} G_K^{11} +
    A^0_{i12} G_K^{12} +
    A^0_{i21} G_K^{21} +
    A^0_{i22} G_K^{22} \\
    &=
    A^0_{i11} G_K^{11} +
    (A^0_{i12} + A^0_{i21}) G_K^{12} +
    A^0_{i22} G_K^{22}
    =
    \bar{a}^0_i \cdot \bar{g}_K,
  \end{split}
\end{equation}
where
\begin{equation}
  \begin{split}
    \bar{a}^0_i &= (A^0_{i11}, A^0_{i12} + A^0_{i21}, A^0_{i22})^{\top}, \\
    \bar{g}_K &= (G_K^{11}, G_K^{12}, G_K^{22})^{\top}.
  \end{split}
\end{equation}
As a consequence, each of the $36$ entries of
the element tensor $A^K$ may be obtained in at most $3$~multiply--add pairs,
and since only $21$ of the entries need to be computed, the total
operation count is directly reduced from $144$ to
$21 \times 3 = 63$.

The set of symmetry-reduced vectors
$\{\bar{a}^0_{11},\bar{a}^0_{12},\ldots,\bar{a}^0_{66}\}$ are given in
Table~\ref{tab:A0,P2,symmetric}. We immediately note a number of
complexity-reducing relations between the vectors. Entries
$\bar{a}^0_{12} = (1,1,0)^{\top}$, $\bar{a}^0_{16} =
(-4,-4,0)^{\top}$, $\bar{a}^0_{26} = (-4,-4,0)^{\top}$ and
$\bar{a}^0_{45} = (-8,-8,0)^{\top}$ are collinear, entries
$\bar{a}^0_{44} = (8,8,8)^{\top}$ and $\bar{a}^0_{45} =
(-8,-8,0)^{\top}$ are close in Hamming distance\footnotemark{} etc.

\footnotetext{We use an extended concept of Hamming distance by
  allowing an optional negation of vectors (which is cheap to compute).}

\linespread{2.0}
\begin{table}[htbp]
  \begin{center}
    \small
    \newcommand{\toffset}{1.3cm}
    \newcommand{\tspace}{1.5cm}
    \begin{tabular}{p{\toffset}p{\tspace}p{\tspace}p{\tspace}p{\tspace}p{\tspace}p{\tspace}}
      &1 & 2 & 3 & 4 & 5 & 6
    \end{tabular}
    \\
    \begin{tabular}{r}
      1 \\ 2 \\ 3 \\ 4 \\ 5 \\ 6
    \end{tabular}
    \begin{tabular}{|r|r|r|r|r|r|}
      \hline
      $(3, 6, 3)^{\top}$
      &
      $(1, 1, 0)^{\top}$
      &
      $(0, 1, 1)^{\top}$
      &
      $(0, 0, 0)^{\top}$
      &
      $-(0, 4, 4)^{\top}$
      &
      $-(4, 4, 0)^{\top}$
      \\\hline
      &
      $(3, 0, 0)^{\top}$
      &
      $-(0, 1, 0)^{\top}$,
      &
      $(0, 4, 0)^{\top}$
      &
      $(0, 0, 0)^{\top}$
      &
      $-(4, 4, 0)^{\top}$
      \\\hline
      &
      &
      $(0, 0, 3)^{\top}$
      &
      $(0, 4, 0)^{\top}$
      &
      $-(0, 4, 4)^{\top}$
      &
      $(0, 0, 0)^{\top}$
      \\\hline
      &
      &
      &
      $(8, 8, 8)^{\top}$
      &
      $-(8, 8, 0)^{\top}$
      &
      $-(0, 8, 8)^{\top}$
      \\\hline
      &
      &
      &
      &
      $(8, 8, 8)^{\top}$
      &
      $(0, 8, 0)^{\top}$
      \\\hline
      &
      &
      &
      &
      &
      $(8, 8, 8)^{\top}$
      \\\hline
    \end{tabular}
    \linespread{1.0}
    \caption{The upper triangular part of the symmetry-reduced
    reference tensor~$A^0$ for Poisson's equation with quadratic elements on
    triangles.}.
    \label{tab:A0,P2,symmetric}
  \end{center}
\end{table}
\linespread{1.0}

\begin{figure}[htbp]
  \begin{center}
    \small
    \psfrag{A00 = (1, 1, 1)}{$\,\, \bar{a}^0_{11} =  (3, 6, 3)^{\top}$}
    \psfrag{A01 = (1, 1, 1)}{$\,\, \bar{a}^0_{12} =  (1, 1, 0)^{\top}$}
    \psfrag{A02 = (1, 1, 1)}{$\,\, \bar{a}^0_{13} =  (0, 1, 1)^{\top}$}
    \psfrag{A03 = (1, 1, 1)}{$\,\, \bar{a}^0_{14} =  (0, 0, 0)^{\top}$}
    \psfrag{A04 = (1, 1, 1)}{$     \bar{a}^0_{15} = -(0, 4, 4)^{\top}$}
    \psfrag{A05 = (1, 1, 1)}{$     \bar{a}^0_{16} = -(4, 4, 0)^{\top}$}
    \psfrag{A11 = (1, 1, 1)}{$\,\, \bar{a}^0_{22} =  (3, 0, 0)^{\top}$}
    \psfrag{A12 = (1, 1, 1)}{$     \bar{a}^0_{23} = -(0, 1, 0)^{\top}$}
    \psfrag{A13 = (1, 1, 1)}{$\,\, \bar{a}^0_{24} =  (0, 4, 0)^{\top}$}
    \psfrag{A14 = (1, 1, 1)}{$\,\, \bar{a}^0_{25} =  (0, 0, 0)^{\top}$}
    \psfrag{A15 = (1, 1, 1)}{$     \bar{a}^0_{26} = -(4, 4, 0)^{\top}$}
    \psfrag{A22 = (1, 1, 1)}{$\,\, \bar{a}^0_{33} =  (0, 0, 3)^{\top}$}
    \psfrag{A23 = (1, 1, 1)}{$\,\, \bar{a}^0_{34} =  (0, 4, 0)^{\top}$}
    \psfrag{A24 = (1, 1, 1)}{$     \bar{a}^0_{35} = -(0, 4, 4)^{\top}$}
    \psfrag{A25 = (1, 1, 1)}{$\,\, \bar{a}^0_{36} =  (0, 0, 0)^{\top}$}
    \psfrag{A33 = (1, 1, 1)}{$\,\, \bar{a}^0_{44} =  (8, 8, 8)^{\top}$}
    \psfrag{A34 = (1, 1, 1)}{$     \bar{a}^0_{45} = -(8, 8, 0)^{\top}$}
    \psfrag{A35 = (1, 1, 1)}{$     \bar{a}^0_{46} = -(0, 8, 8)^{\top}$}
    \psfrag{A44 = (1, 1, 1)}{$\,\, \bar{a}^0_{55} =  (8, 8, 8)^{\top}$}
    \psfrag{A45 = (1, 1, 1)}{$\,\, \bar{a}^0_{56} =  (0, 8, 0)^{\top}$}
    \psfrag{A55 = (1, 1, 1)}{$\,\, \bar{a}^0_{66} =  (8, 8, 8)^{\top}$}
    \includegraphics[width=15cm]{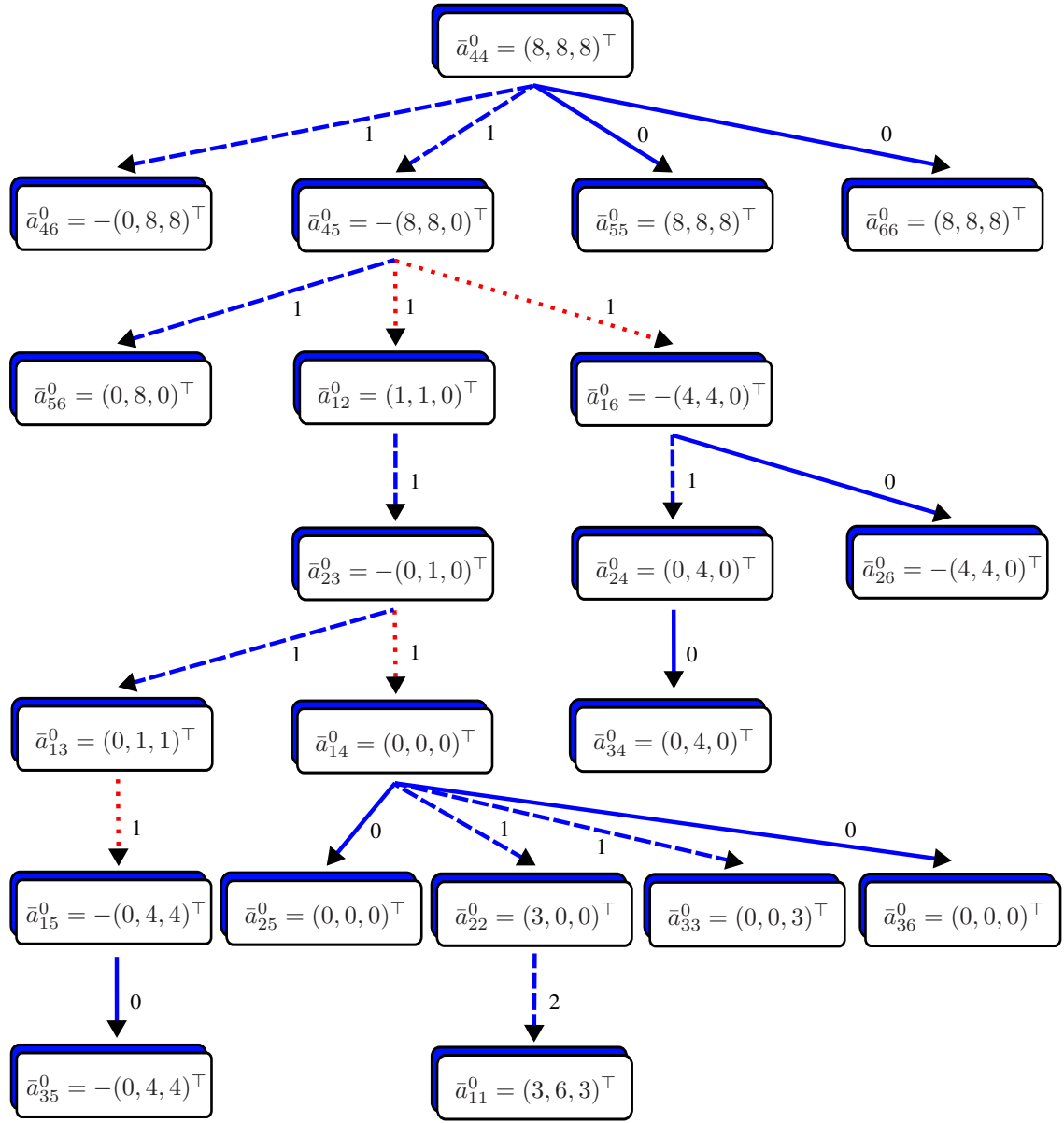}
    \vspace{1cm}
    \caption{The minimum spanning tree for the optimized computation
      of the upper triangular part (Table~\ref{tab:A0,P2,symmetric}) of the
      element tensor  for Poisson's equation with quadratic elements
      on triangles. Solid (blue) lines indicate zero Hamming distance
      (equality), dashed (blue) lines indicate a small but nonzero
      Hamming distance and dotted (red) lines indicate collinearity.}
    \label{fig:mst}
  \end{center}
\end{figure}

To systematically explore these dependencies, we form a weighted graph
$G = (V,E)$ and compute a minimum spanning tree. We let the vertices
$V$ be the set of symmetry-reduced vectors, $V =
\{\bar{a}^0_{11},\bar{a}^0_{12},\ldots,\bar{a}^0_{66}\}$, and form the
set of edges $E$ by adding between each pair of vertices an edge with
weight given by the minimum of all complexity-reducing relations
between the two vertices. The resulting minimum spanning tree is shown
in Figure~\ref{fig:mst}. We note that the total edge weight of the
minimum spanning tree is $14$. This means that once the value of the
entry in the element tensor corresponding to the root vertex is known,
the remaining entries may be computed in at most $14$ multiply--add
pairs. Adding the $3$ multiply--add pairs needed to compute the root
entry, we thus find that all $36$ entries of the element tensor $A^K$
may be computed in at most $17$ multiply--add pairs.

An optimized algorithm for the computation of the element tensor~$A^K$
may then be found by starting at the root vertex and computing the
remaining entries by traversing the minimum spanning tree, as
demonstrated in Algorithm~\ref{alg:traversing}. Note that there are
several ways to traverse the tree. In particular, it is possible to
pick any vertex as the root vertex and start from there. Furthermore,
there are many ways to traverse the tree given the root
vertex. Algorithm~\ref{alg:traversing} is generated by traversing the
tree breadth-first, starting at the root vertex $\bar{a}^0_{44} =
(8,8,8)^{\top}$. Finally, we note that the operation count may be
further reduced by not counting multiplications with zeros and ones.

\begin{algorithm}
  \small
  \begin{tabular}{l|l}
    \begin{minipage}{8cm}
    $A^K_{44} = A^0_{4411} G_K^{11} + (A^0_{4412} + A^0_{4421})
    G_K^{12} +  A^0_{4422} G_K^{22}$ \\
    $A^K_{46} = -A^K_{44} + 8 G_K^{11}$ \\
    $A^K_{45} = -A^K_{44} + 8 G_K^{22}$ \\
    $A^K_{55} = A^K_{44}$ \\
    $A^K_{66} = A^K_{44}$ \\
    $A^K_{56} = -A^K_{45} - 8 G_K^{11}$ \\
    $A^K_{12} = -\frac{1}{8} A^K_{45}$ \\
    $A^K_{16} = \frac{1}{2} A^K_{45}$ \\
    $A^K_{23} = -A^K_{12} + 1 G_K^{11}$ \\
    $A^K_{24} = -A^K_{16} - 4 G_K^{11}$ \\
    $A^K_{26} = A^K_{16}$
    \end{minipage}
    &
    \begin{minipage}{8cm}
    $A^K_{13} = -A^K_{23} + 1 G^K_{22}$ \\
    $A^K_{14} = 0 A^K_{23} $ \\
    $A^K_{34} = A^K_{24}$ \\
    $A^K_{15} = -4 A^K_{13}$ \\
    $A^K_{25} = A^K_{14}$ \\
    $A^K_{22} = A^K_{14} + 3 G^K_{11}$ \\
    $A^K_{33} = A^K_{14} + 3 G^K_{22}$ \\
    $A^K_{36} = A^K_{14}$ \\
    $A^K_{35} = A^K_{15}$ \\
    $A^K_{11} = A^K_{22} + 6 G^K_{12} + 3 G^K_{22}$ \\
    \end{minipage}
  \end{tabular}
  \caption{An optimized (but not optimal) algorithm for computing the
    upper triangular part of the element tensor~$A^K$ for Poisson's
    equation with quadratic elements on triangles in $17$ multiply--add
    pairs.}
  \label{alg:traversing}
\end{algorithm}

% Exercise: Find something more to optimize in Algorithm~\ref{alg:traversing}

\subsubsection{Extensions}

By use of symmetry and relations between subsets of the reference
tensor $A^0$ we have seen that it is possible to significantly reduce
the operation count for the computation of the tensor contraction $A^K
= A^0 : G_K$. We have here only discussed the use of binary relations
(collinearity and Hamming distance) but further reductions may be made
by considering ternary relations, such as coplanarity, and
higher-arity relations between the vectors.

%------------------------------------------------------------------------------
\section{AUTOMATION AND SOFTWARE ENGINEERING}
\label{sec:se}

In this section, we comment briefly on some topics of software
engineering relevant to the automation of the finite element method.
A number of books and papers have been written on the subject of
software engineering for the implementation of finite element methods,
see for
example~\cite{Lan99,Lan05,ArgBru97,Mac92,MasUsm97}.
In particular these works point out the importance of
object-oriented, or concept-oriented, design in developing mathematical
software; since the mathematical concepts have already been hammered
out, it may be advantageous to reuse these concepts in the system
design, thus providing abstractions for important concepts, including
\texttt{Vector}, \texttt{Matrix}, \texttt{Mesh}, \texttt{Function},
\texttt{BilinearForm}, \texttt{LinearForm}, \texttt{FiniteElement}
etc.

\begin{figure}[htbp]
  \begin{center}
    \includegraphics[width=12cm]{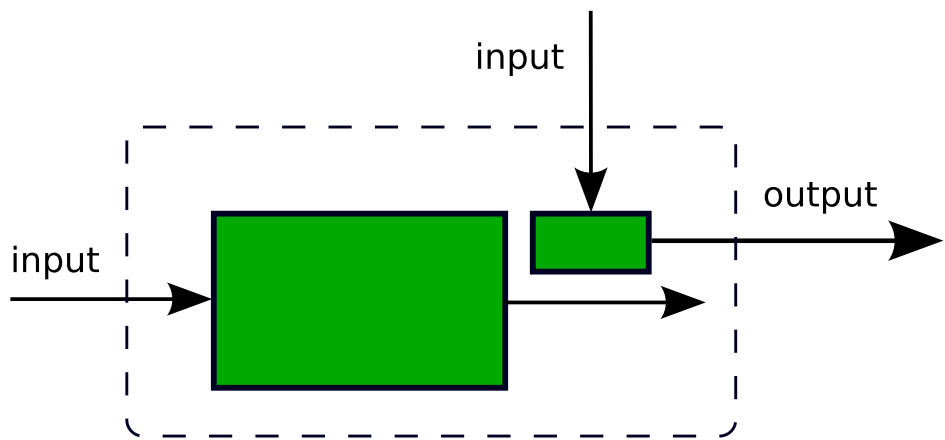}
    \caption{A machine (computer program) that automates the finite
      element method by automatically generating a particular machine
      (computer program) for a suitable subset of the given input
      data.}
    \label{fig:machines}
  \end{center}
\end{figure}

We shall not repeat these arguments, but instead point out a
couple of issues that might be less obvious. In particular, a
straightforward implementation of all the mathematical concepts
discussed in the previous sections may be difficult or even
impossible to attain. Therefore, we will argue that a level of
automation is needed also in the implementation or \emph{realization}
of an automation of the finite element method, that is, the automatic
generation of computer code for the specific mathematical concepts
involved in the specification of any particular finite element method
and differential equation, as illustrated in
Figure~\ref{fig:machines}.

We also point out that the automation of the finite element method is
not only a software engineering problem. In addition to identifying
and implementing the proper mathematical concepts, one must develop
new mathematical tools and ideas that make it possible for the
automating system to realize the full generality of the finite element
method. In addition, new insights are needed to build an
\emph{efficient} automating system that can compete with or outperform
hand-coded specialized systems for any given input.

\subsection{Code Generation}

As in all types of engineering, software for scientific computing must
try to find a suitable trade-off between generality and efficiency; a
software system that is general in nature, that is, it accepts a wide
range of inputs, is often less efficient than another software
system that performs the same job on a more limited set of
inputs. As a result, most codes used by practitioners for the solution
of differential equations are very specific, often specialized to a
specific method for a specific differential equation.

However, by using a compiler approach, it is possible to combine
generality and efficiency without loss of generality and without loss
of efficiency. Instead of developing a potentially inefficient
general-purpose program that accepts a wide range of inputs, a
suitable subset of the input is given to an optimizing compiler that
generates a specialized program that takes a more limited set of
inputs. In particular, one may automatically generate a specialized
simulation code for any given method and differential equation.

An important task is to identify a suitable subset of the input to
move to a precompilation phase. In the case of a system automating the
solution of differential equations by the finite element method, a
suitable subset of input includes the variational
problem~(\ref{eq:varproblem}) and the choice of approximating finite
element spaces. We thus develop a domain-specific compiler that
accepts as input a variational problem and a set of finite elements
and generates optimized low-level code (in some general-purpose
language such as C or C++). Since the compiler may thus work on only a
small family of inputs (multilinear forms), domain-specific knowledge
allows the compiler to generate very efficient code, using the
optimizations discussed in the previous section. We return to this in
more detail below in Section~\ref{sec:ffc} when we discuss the FEniCS
Form Compiler FFC.

We note that to limit the complexity of the automating system, it is
important to identify a minimal set of code to be generated
at a precompilation stage, and implement the remaining code in a
general-purpose language. It makes less sense to generate the code for
administrative tasks such as reading and writing data to file, special
algorithms like adaptive mesh refinement etc. These tasks can be
implemented as a library in a general-purpose language.

\subsection{Just-In-Time Compilation}

To make an automating system for the solution of differential
equations truly useful, the generation and precompilation of code
according to the above discussion must also be automated. Thus, a user
should ultimately be presented with a single user-interface and the
code should automatically and transparently be generated and compiled
\emph{just-in-time} for a given problem specification.

Achieving just-in-time compilation of variational problems is
challenging, not only to construct the exact mechanism by which code
is generated, compiled and linked back in at run-time, but also to
reduce the precompilation phase to a minimum so that the overhead of
code generation and compilation is acceptable. To compile and generate
the code for the evaluation of a multilinear form as discussed in
Section~\ref{sec:automation,forms}, we need to compute the tensor
representation~(\ref{eq:representation}), including the evaluation of
the reference tensor. Even with an optimized algorithm for the
computation of the reference tensor as discussed
in~\cite{logg:submit:toms-ffc-monomial-2006}, the computation of the
reference tensor may be very costly, especially for high-order
elements and complicated forms. To improve the situation, one may
consider caching previously computed reference tensors (similarly
to how~\LaTeX{} generates and caches fonts in different resolutions)
and reuse previously computed reference tensors. As discussed
in~\cite{logg:submit:toms-ffc-monomial-2006}, a reference tensor may
be uniquely identified by a (short) string referred to as a
\emph{signature}. Thus, one may store reference tensors along with
their signatures to speed up the precomputation and allow run-time
just-in-time compilation of variational problems with little overhead.

%------------------------------------------------------------------------------
\section{A PROTOTYPE IMPLEMENTATION (FEniCS)}
\label{sec:prototype}

\begin{center}
  \begin{minipage}{10cm}
  \emph{An algorithm must be seen to be believed, and the best way
  to learn what an algorithm is all about is to try it.}

  \vspace{0.5cm}

  \raggedleft
  Donald E. Knuth \\
  \emph{The Art of Computer Programming} (1968) \\
  \end{minipage}
\end{center}

The automation of the finite element method includes its own
realization, that is, a software system that implements the algorithms
discussed in Sections~\ref{sec:fem}--\ref{sec:automation,assembly}.
Such a system is provided by the FEniCS
project~\cite{logg:www:03,logg:preprint:10}. We present below some of
the key components of FEniCS, including FIAT, FFC and DOLFIN, and
point out how they relate to the various aspects of the automation of
the finite element method. In particular, the automatic tabulation of
finite element basis functions discussed in
Section~\ref{sec:automation,tabulation} is provided by
FIAT~\cite{www:FIAT,Kir04,Kir06}, the automatic evaluation of the
element tensor as discussed in Section~\ref{sec:automation,forms} is
provided by
FFC~\cite{logg:www:04,logg:article:10,logg:submit:toms-ffc-monomial-2006,logg:manual:02}
and the automatic assembly
of the discrete system
as discussed in Section~\ref{sec:automation,assembly}
is provided by
DOLFIN~\cite{logg:www:01,logg:preprint:06,logg:manual:01}. The FEniCS
project thus serves as a testbed for development of new ideas for
the automatic and efficient implementation of finite element methods. At
the same, it provides a reference implementation of these ideas.

FEniCS software is free software~\cite{www:freesoftware}. In
particular, the components of FEniCS are licensed under the GNU
General Public License \cite{www:GPL}.\footnotemark{} The source code
is freely available on the FEniCS web site \cite{logg:www:03} and the
development is discussed openly on public mailing lists.

\footnotetext{FIAT is licensed under the Lesser General Public
License~\cite{www:LGPL}.}

\subsection{FIAT}

The FInite element Automatic Tabulator FIAT~\cite{www:FIAT}
was first introduced in~\cite{Kir04} and implements the ideas
discussed above in Section~\ref{sec:automation,tabulation} for
the automatic tabulation of finite element basis functions based on a
linear algebraic representation of function spaces and constraints.

FIAT provides functionality for defining finite element function
spaces as constrained subsets of polynomials on the simplices in one,
two and three space dimensions, as well as a library of predefined
finite elements, including
arbitrary degree
Lagrange~\cite{Cia78,BreSco94},
Hermite~\cite{Cia78,BreSco94},
Raviart--Thomas~\cite{RavTho77a},
Brezzi--Douglas--Marini~\cite{BreDou85}
and
Nedelec~\cite{Ned80}
elements, as well as the (first degree)
Crouzeix--Raviart element~\cite{CroRav73}. Furthermore, the plan is to
support Brezzi--Douglas--Fortin--Marini \cite{BreFor91} and
Arnold--Winther \cite{ArnWin02} elements in future
versions.

In addition to tabulating finite element nodal basis functions (as
linear combinations of a reference basis), FIAT generates quadrature
points of any given order on the reference simplex and provides
functionality for efficient tabulation of the basis functions and
their derivatives at any given set of points. In
Figure~\ref{fig:fiatexample,1} and Figure~\ref{fig:fiatexample,2}, we
present some examples of basis functions generated by FIAT.

Although FIAT is implemented in Python, the interpretive overhead of
Python compared to compiled languages is small, since the operations
involved may be phrased in terms of standard linear algebra
operations, such as the solution of linear systems and singular value
decomposition, see~\cite{Kir06}. FIAT may thus make use of optimized
Python linear algebra libraries such Python
Numeric~\cite{www:Numeric}. Recently, a C++ version of FIAT
called FIAT++ has also been developed with run-time bindings for
Sundance~\cite{www:Sundance,Lon03,Lon04}.

\begin{figure}[htbp]
  \begin{center}
    % These are basis functions number 1, 2, 3
    \includegraphics[width=14cm]{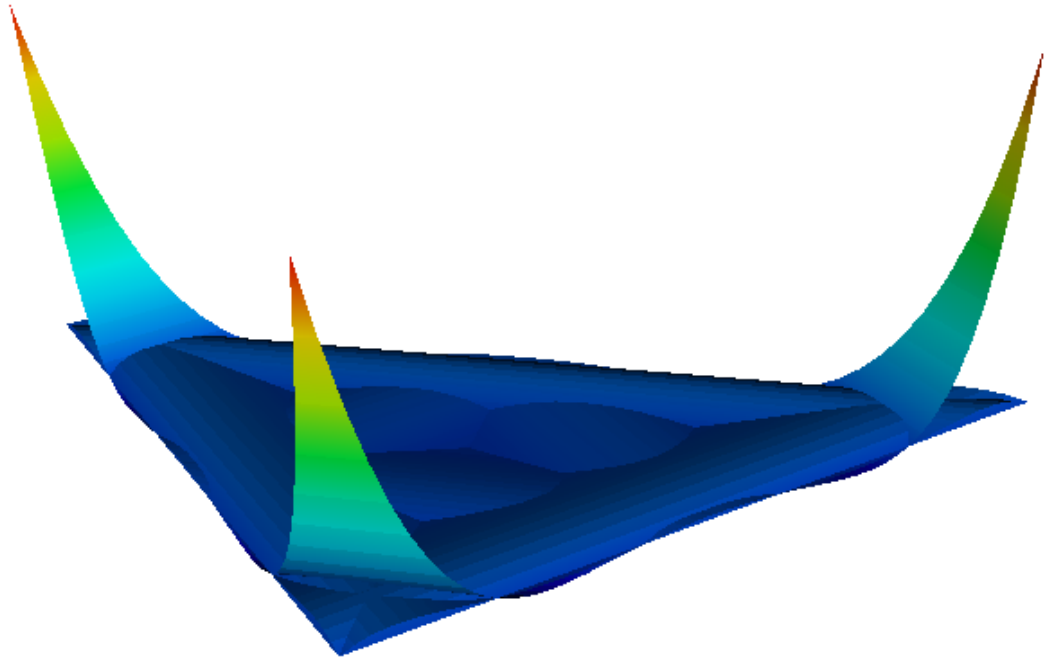}
    \caption{The first three basis functions for a fifth-degree Lagrange
      finite element on a triangle, associated with the three vertices
      of the triangle. (Courtesy of Robert C. Kirby.)}
    \label{fig:fiatexample,1}
  \end{center}
\end{figure}

\begin{figure}[htbp]
  \begin{center}
    % This is basis function number 21 (the last one)
    \includegraphics[width=14cm]{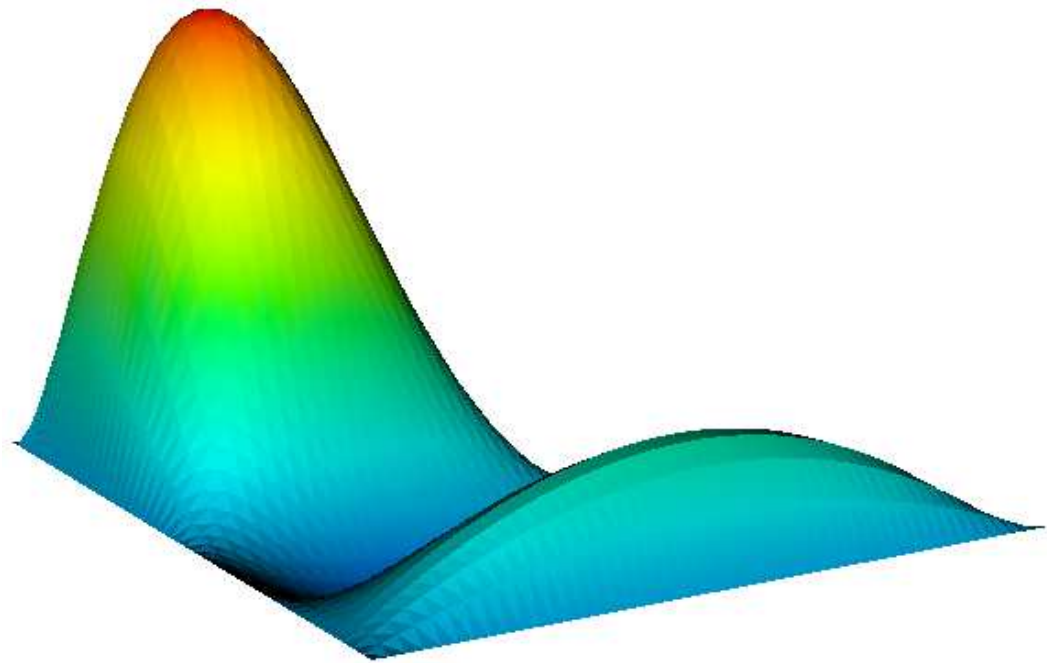}
    \caption{A basis function associated with an interior point for a
      fifth-degree Lagrange finite element on a triangle. (Courtesy of
      Robert C. Kirby.)}
    \label{fig:fiatexample,2}
  \end{center}
\end{figure}

\subsection{FFC}
\label{sec:ffc}

\begin{figure}[htbp]
  \begin{center}
    \includegraphics[width=12cm]{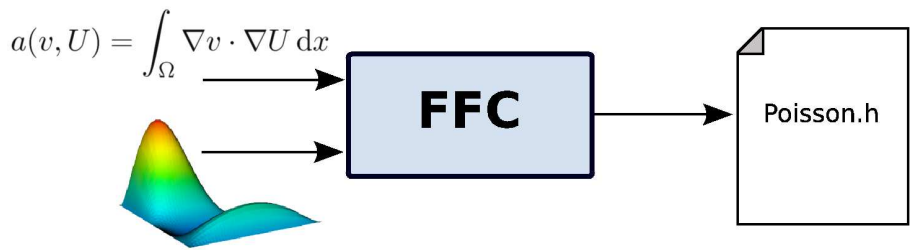}
    \caption{The form compiler FFC takes as input a multilinear form
      together with a set of function spaces and generates optimized
      low-level (C++) code for the evaluation of the associated
      element tensor.}
    \label{fig:ffc}
  \end{center}
\end{figure}

The FEniCS Form Compiler FFC~\cite{logg:www:04}, first introduced
in~\cite{logg:article:10}, automates the evaluation of multilinear forms as
outlined in Section~\ref{sec:automation,forms} by automatically
generating code for the efficient computation of the element tensor
corresponding to a given multilinear form. FFC thus functions as
domain-specific compiler for multilinear forms, taking as input a set
of discrete function spaces together with a multilinear form defined
on these function spaces, and produces as output optimized
low-level code, as illustrated in Figure~\ref{fig:ffc}. In its
simplest form, FFC generates code in the form of a single C++ header
file that can be included in a C++ program, but FFC can also be used
as a just-in-time compiler within a scripting
environment like Python, for seamless definition and evaluation of
multilinear forms.

\subsubsection{Form language}

The FFC form language is generated from a small set of basic data
types and operators that allow a user to define a wide range of
multilinear forms, in accordance with the discussion of
Section~\ref{sec:operatoroverloading}. As an illustration, we include
below the complete definition in the FFC form language of the bilinear
forms for the test cases considered above in
Tables~\ref{tab:testcase1}--\ref{tab:testcase4}. We refer to the FFC
user manual~\cite{logg:manual:02} for a detailed discussion of the
form language, but note here that in addition to a set of standard
operators, including the inner product~\texttt{dot}, the partial
derivative~\texttt{D}, the gradient~\texttt{grad}, the
divergence~\texttt{div} and the rotation~\texttt{rot}, FFC
supports Einstein tensor-notation (Table~\ref{tab:ffc,testcase3}) and
user-defined operators (operator \texttt{epsilon} in
Table~\ref{tab:ffc,testcase4}).

\begin{table}[htbp]
  \begin{center}
    \begin{code}
 element = FiniteElement("Lagrange", "tetrahedron", 1)

 v = BasisFunction(element)
 U = BasisFunction(element)

 a = v*U*dx
    \end{code}
    \caption{The complete definition of the bilinear form
      $a(v,U) = \int_{\Omega} v \, U \dx$
      in the FFC form language (test case 1).}
    \label{tab:ffc,testcase1}
  \end{center}
\end{table}

\begin{table}[htbp]
  \begin{center}
    \begin{code}
 element = FiniteElement("Lagrange", "tetrahedron", 1)

 v = BasisFunction(element)
 U = BasisFunction(element)

 a = dot(grad(v), grad(U))*dx
    \end{code}
    \caption{The complete definition of the bilinear form
      $a(v,U) = \int_{\Omega} \nabla v \cdot \nabla U \dx$
      in the FFC form language (test case 2).}
    \label{tab:ffc,testcase2}
  \end{center}
\end{table}

\begin{table}[htbp]
  \begin{center}
    \begin{code}
 element = FiniteElement("Vector Lagrange", "tetrahedron", 1)

 v = BasisFunction(element)
 U = BasisFunction(element)
 w = Function(element)

 a = v[i]*w[j]*D(U[i],j)*dx
    \end{code}
    \caption{The complete definition of the bilinear form
      $a(v,U) = \int_{\Omega} v \cdot (w \cdot \nabla) U \dx$
      in the FFC form language (test case 3).}
    \label{tab:ffc,testcase3}
  \end{center}
\end{table}

\begin{table}[htbp]
  \begin{center}
    \begin{code}
 element = FiniteElement("Vector Lagrange", "tetrahedron", 1)

 v = BasisFunction(element)
 U = BasisFunction(element)

 def epsilon(v):
     return 0.5*(grad(v) + transp(grad(v)))

 a = dot(epsilon(v), epsilon(U))*dx
    \end{code}
    \caption{The complete definition of the bilinear form
      $a(v,U) = \int_{\Omega} \epsilon(v) : \epsilon(U) \dx$
      in the FFC form language (test case 4).}
    \label{tab:ffc,testcase4}
  \end{center}
\end{table}

\subsubsection{Implementation}

The FFC form language is implemented in Python as a collection of
Python classes (including \texttt{BasisFunction}, \texttt{Function},
\texttt{FiniteElement} etc.) and operators on theses classes. Although
FFC is implemented in Python, the interpretive overhead of Python has
been minimized by judicious use of optimized numerical libraries such
as Python Numeric~\cite{www:Numeric}. The computationally most expensive
part of the compilation of a multilinear form is the precomputation of
the reference tensor. As demonstrated
in~\cite{logg:submit:toms-ffc-monomial-2006}, by suitably
pretabulating basis functions and their derivatives at a set of
quadrature points (using FIAT), the reference
tensor can be computed by assembling a set of outer products, which may each
be efficiently computed by a call to Python Numeric.

Currently, the only optimization FFC makes is to avoid multiplications
with any zeros of the reference tensor~$A^0$ when generating code for the
tensor contraction~$A^K = A^0 : G_K$. As part of the FEniCS project,
an optimizing backend, FErari (Finite Element Re-arrangement
Algorithm to Reduce Instructions), is currently being
developed. Ultimately, FFC will call FErari at compile-time to find an
optimized computation of the tensor contraction, according to the
discussion in Section~\ref{sec:optimizations}.

\subsubsection{Benchmark results}
\label{sec:benchmarks}

As a demonstration of the efficiency of the code generated by FFC, we
include in Table~\ref{tab:speedup} a comparison taken
from~\cite{logg:article:10} between a standard implementation, based
on computing the element tensor~$A^K$ on each cell~$K$ by a loop
over quadrature points, with the code automatically generated by FFC,
based on precomputing the reference tensor~$A^0$ and computing the
element tensor~$A^K$ by the tensor contraction~$A^K = A^0 : G_K$ on
each cell.

\begin{table}[htbp]
  \begin{center}
    \small
    \begin{tabular}{|l|r|r|r|r|r|r|r|r|}
      \hline
      Form & $q = 1$ & $q = 2$ & $q = 3$ & $q = 4$ & $q = 5$ & $q = 6$
      & $q = 7$ &$q = 8$ \\
      \hline
      \hline
      Mass 2D           & 12 & 31  & 50  & 78  & 108 & 147 & 183  &
      232  \\
      Mass 3D           & 21 & 81  & 189 & 355 & 616 & 881 & 1442 &
      1475 \\
      Poisson 2D        & 8  & 29  & 56  & 86  & 129 & 144 & 189  &
      236  \\
      Poisson 3D        & 9  & 56  & 143 & 259 & 427 & 341 & 285  &
      356  \\
      Navier--Stokes 2D & 32 & 33  & 53  & 37  & --- & --- & --- & ---
      \\
      Navier--Stokes 3D & 77 & 100 & 61  & 42  & --- & --- & --- & ---
      \\
      Elasticity 2D     & 10 & 43  & 67  & 97  & --- & --- & --- & ---
      \\
      Elasticity 3D     & 14 & 87  & 103 & 134 & --- & --- & --- & ---
      \\
      \hline
    \end{tabular}
    \normalsize
    \caption{Speedups for test cases 1--4
      (Tables~\ref{tab:testcase1}--\ref{tab:testcase4} and Tables~\ref{tab:ffc,testcase1}--\ref{tab:ffc,testcase4}) in two and
      three space dimensions.}
    \label{tab:speedup}
  \end{center}
\end{table}

As seen in Table~\ref{tab:speedup}, the speedup ranges between one and
three orders of magnitude, with larger speedups for higher degree
elements. In Figure~\ref{fig:result,1} and Figure~\ref{fig:result,2},
we also plot the dependence of the speedup on the polynomial degree
for test cases~1 and~2 respectively.

 It should be noted that the total work in a
simulation also includes the assembly of the local element tensors
$\{A^K\}_{K\in\mathcal{T}}$ into the global tensor~$A$, solving the
linear system, iterating on the nonlinear problem etc. Therefore, the
overall speedup may be significantly less than the speedups reported
in Table~\ref{tab:speedup}. We note that if the computation of the
local element tensors normally accounts for a fraction~$\theta \in
(0,1)$ of the total run-time, then the overall speedup gained by
a speedup of size~$s > 1$ for the computation of the element tensors
will be
\begin{equation}
  1 < \frac{1}{1-\theta + \theta/s} \leq \frac{1}{1 - \theta},
\end{equation}
which is significant only if $\theta$ is significant. As noted
in~\cite{logg:article:07}, $\theta$ may be significant in many cases,
in particular for nonlinear problems where a nonlinear system (or the
action of a linear operator) needs to be repeatedly reassembled as
part of an iterative method.

\begin{figure}[htbp]
  \begin{center}
    \psfrag{q}{$q$}
    \includegraphics[width=12cm]{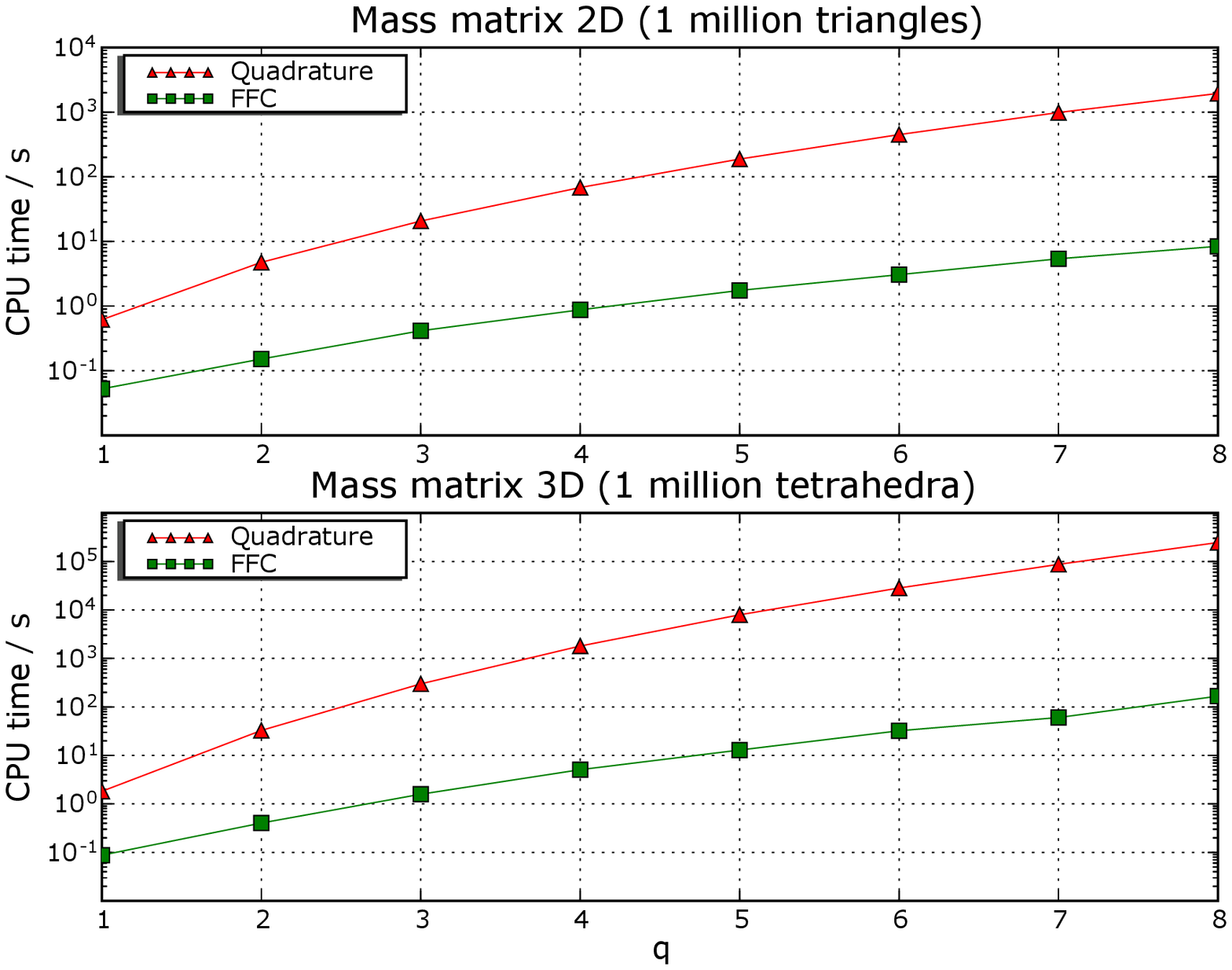}
    \caption{Benchmark results for test case~1, the mass matrix,
      specified in FFC by \texttt{a = v*U*dx}.}
    \label{fig:result,1}
  \end{center}
\end{figure}

\begin{figure}[htbp]
  \begin{center}
    \psfrag{q}{$q$}
    \includegraphics[width=12cm]{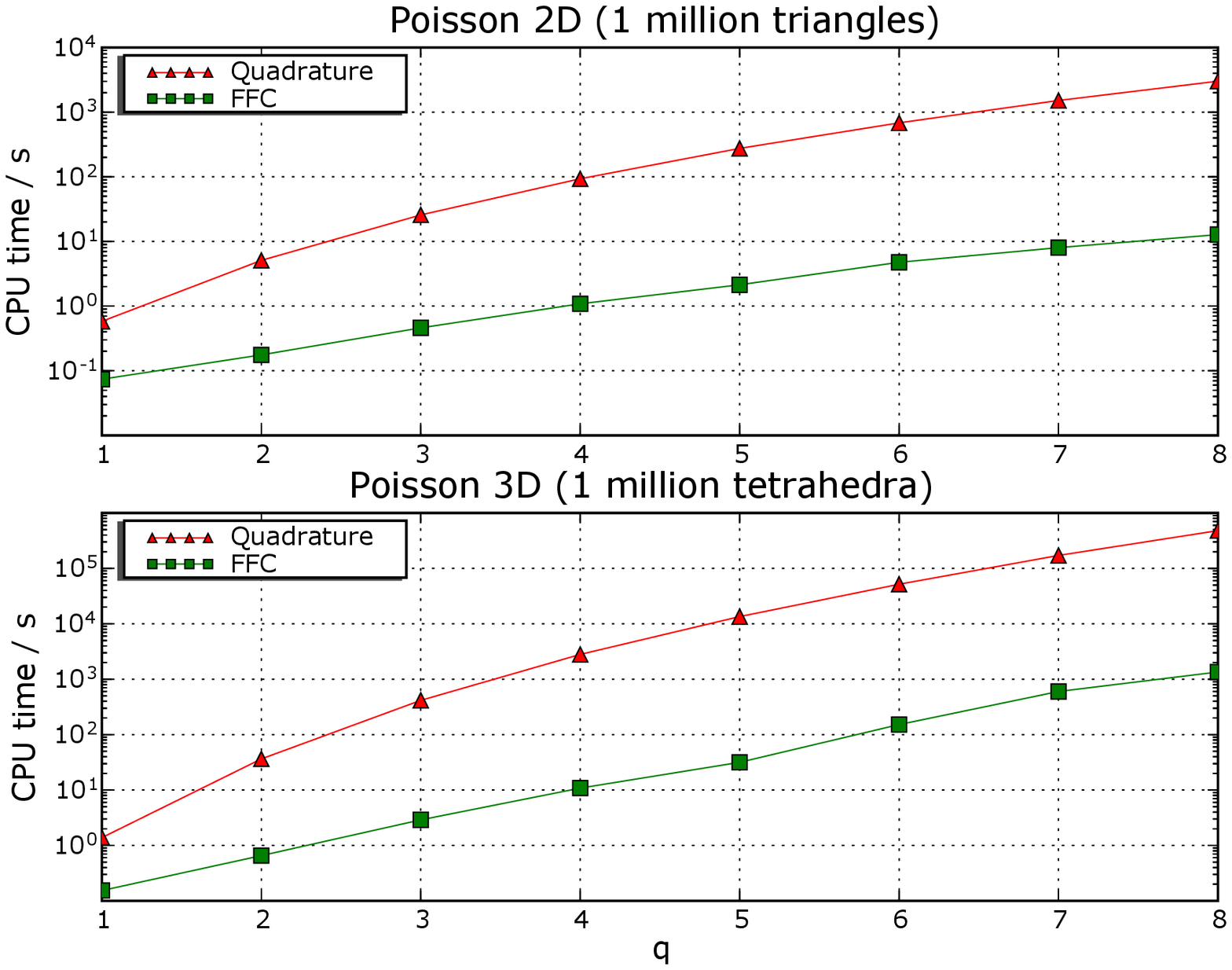}
    \caption{Benchmark results for test case~2, Poisson's equation,
      specified in FFC by \texttt{a = dot(grad(v), grad(U))*dx}.}
    \label{fig:result,2}
  \end{center}
\end{figure}

\subsubsection{User interfaces}

FFC can be used either as a stand-alone compiler on the command-line,
or as a Python module from within a Python script. In the first case,
a multilinear form (or a pair of bilinear and linear forms) is entered
in a text file with suffix~\texttt{.form} and then compiled by calling
the command~\texttt{ffc} with the form file on the command-line.

\begin{figure}[htbp]
  \begin{center}
    \includegraphics[width=14cm]{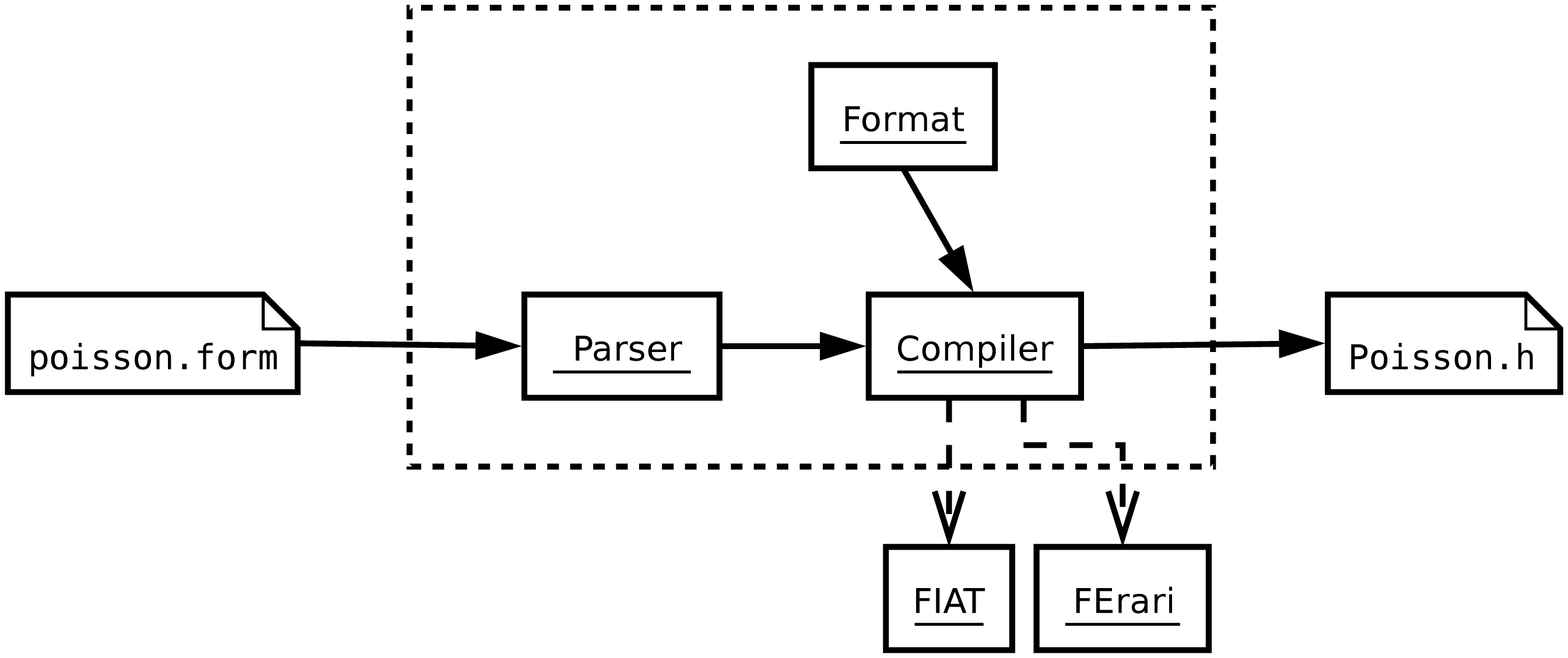}
    \caption{Component diagram for FFC.}
    \label{fig:ffccomponents}
  \end{center}
\end{figure}

By default, FFC generates C++ code for inclusion in a DOLFIN C++ program
(see Section~\ref{sec:dolfin} below) but FFC can also compile code for
other backends (by an appropriate compiler flag), including the ASE
(ANL SIDL Environment) format~\cite{www:ASE}, XML format, and \LaTeX{}
format (for inclusion of the tensor representation in reports and
presentations). The format of the generated code is separated from the
parsing of forms and the generation of the tensor contraction, and new
formats for alternative backends may be added with little effort,
see Figure~\ref{fig:ffccomponents}.

Alternatively, FFC can be used directly from within Python as a Python
module, allowing definition and compilation of multilinear forms from
within a Python script. If used together with the recently developed
Python interface of DOLFIN (PyDOLFIN), FFC functions as a just-in-time
compiler for multilinear forms, allowing forms to be defined and
evaluated from within Python.

\subsection{DOLFIN}
\label{sec:dolfin}

DOLFIN~\cite{logg:www:01,logg:preprint:06,logg:manual:01}, Dynamic
Object-oriented Library for FINite element computation, functions as a
general programming interface to DOLFIN and provides a problem-solving
environment (PSE) for differential equations in the form of a C++/Python
class library.

Initially, DOLFIN was developed as a self-contained (but modularized)
C++ code for finite element simulation, providing basic functionality
for the definition and automatic evaluation of multilinear forms,
assembly, linear algebra, mesh data structures and adaptive mesh
refinement, but as a consequence of the development of focused
components for each of these tasks as part of the FEniCS project, a
large part (but not all) of the functionality of DOLFIN has been
delegated to these other components while maintaining a consistent
programming interface. Thus, DOLFIN relies on FIAT for the automatic
tabulation of finite element basis functions and on FFC for the
automatic evaluation of multilinear forms. We discuss below some of
the key aspects of DOLFIN and its role as a component of the FEniCS
project.

\subsubsection{Automatic assembly of the discrete system}

DOLFIN implements the automatic assembly of the discrete system
associated with a given variational problem as outlined in
Section~\ref{sec:automation,assembly}. DOLFIN iterates over the
cells $\{K\}_{K\in\mathcal{T}}$ of a given mesh~$\mathcal{T}$ and
calls the code generated by FFC on each cell~$K$ to evaluate the
element tensor~$A^K$. FFC also generates the code for the
local-to-global mapping which DOLFIN calls to obtain a rule for the
addition of each element tensor~$A^K$ to the global tensor~$A$.

Since FFC generates the code for both the evaluation of the element
tensor and for the local-to-global mapping, DOLFIN needs to know very
little about the finite element method. It only follows the
instructions generated by FFC and operates abstractly on the level of
Algorithm~\ref{alg:assembly,2}.

\subsubsection{Meshes}

DOLFIN provides basic data structures and algorithms for simplicial
meshes in two and three space dimensions (triangular and tetrahedral
meshes) in the form of a class \texttt{Mesh}, including adaptive mesh
refinement. As part of PETSc~\cite{www:PETSc,BalBus04,BalEij97} and the FEniCS
project, the new component Sieve~\cite{www:Sieve,KneKar05} is currently
being developed. Sieve generalizes the mesh concept and provides
powerful abstractions for dimension-independent operations on mesh
entities and will function as a backend for the mesh data structures
in DOLFIN.

\subsubsection{Linear algebra}

Previously, DOLFIN provided a stand-alone basic linear algebra library
in the form of a class \texttt{Matrix}, a class \texttt{Vector} and a
collection of iterative and direct solvers. This implementation has
recently been replaced by a set of simple wrappers for the sparse
linear algebra library provided by PETSc~\cite{www:PETSc,BalBus04,BalEij97}. As a
consequence, DOLFIN is able to provide sophisticated
high-performance parallel linear algebra with an easy-to-use
object-oriented interface suitable for finite element computation.

\subsubsection{ODE solvers}

DOLFIN also provides a set of general order mono-adaptive and
multi-adaptive~\cite{logg:article:01,logg:article:02,logg:article:03,logg:article:05,logg:article:08}
ODE-solvers, automating the solution of ordinary differential
equations. Although the ODE-solvers may be used in connection with the
automated assembly of discrete systems, DOLFIN does currently not
provide any level of automation for the discretization of
time-dependent PDEs. Future versions of DOLFIN (and FFC) will allow
time-dependent PDEs to be defined directly in the FFC form language
with automatic discretization and adaptive time-integration.

\subsubsection{PDE solvers}
\label{sec:solvers}

In addition to providing a class library of basic tools that automate
the implementation of adaptive finite element methods, DOLFIN provides
a collection of ready-made solvers for a number of standard
equations. The current version of DOLFIN provides solvers for
Poisson's equation, the heat equation, the convection--diffusion
equation, linear elasticity, updated large-deformation elasticity, the
Stokes equations and the incompressible Navier--Stokes equations.

\subsubsection{Pre- and post-processing}

DOLFIN relies on interaction with external tools for pre-processing
(mesh generation) and post-processing (visualization). A number of
output formats are provided for visualization, including DOLFIN
XML~\cite{logg:manual:01}, VTK~\cite{www:VTK} (for use in
ParaView~\cite{www:ParaView} or MayaVi~\cite{www:MayaVi}),
Octave~\cite{www:Octave}, MATLAB~\cite{www:MATLAB}, OpenDX~\cite{www:OpenDX},
GiD~\cite{www:GiD} and Tecplot~\cite{www:Tecplot}. DOLFIN may also be easily
extended with new output formats.

\subsubsection{User interfaces}

DOLFIN can be accessed either as a C++ class library or as a Python
module, with the Python interface generated semi-automatically from
the C++ class library using SWIG~\cite{www:SWIG,Bea96}. In both cases,
the user is presented with a simple and consistent but powerful
programming interface.

As discussed in Section~\ref{sec:solvers}, DOLFIN provides a set of
ready-made solvers for standard differential equations. In the
simplest case, a user thus only needs to supply a mesh, a set of
boundary conditions and any parameters and variable coefficients to
solve a differential equation, by calling one of the existing solvers.
For other differential equations, a solver may be implemented with
minimal effort using the set of tools provided by the DOLFIN class
library, including variational problems, meshes and linear algebra as
discussed above.

\subsection{Related Components}

We also mention two other projects developed as part of FEniCS. One of
these is Puffin\cite{logg:www:05,logg:manual:03}, a light-weight
educational implementation of the basic functionality of FEniCS for
Octave/MATLAB, including automatic assembly of the linear system from
a given variational problem. Puffin has been used with great success
in introductory undergraduate mathematics courses and is accompanied
by a set of exercises~\cite{logg:www:06} developed as part of the
Body~and~Soul reform
project~\cite{logg:www:07,EriEst03a,EriEst03b,EriEst03c} for applied
mathematics education.

The other project is the Ko mechanical simulator~\cite{www:Ko}. Ko
uses DOLFIN as the computational backend and provides a specialized
interface to the simulation of mechanical systems, including
large-deformation elasticity and collision detection. Ko provides two
different modes of simulation: either a simple mass--spring model solved as a
system of ODEs, or a large-deformation updated elasticity
model~\cite{logg:inprep:updatedelasticity} solved as a system of
time-dependent PDEs. As a consequence of the efficient assembly
provided by DOLFIN, based on efficient code being generated by
FFC, the overhead of the more complex PDE model compared to the simple
ODE model is relatively small.

%------------------------------------------------------------------------------
\section{EXAMPLES}
\label{sec:examples}

In this section, we present a number of examples chosen to illustrate
various aspects of the implementation of finite element methods for a
number of standard partial differential equations with the FEniCS
framework. We already saw in Section~\ref{sec:survey} the
specification of the variational problem for Poisson's equation in the
FFC form language. The examples below include static linear
elasticity, two different formulations for the Stokes equations and
the time-dependent convection--diffusion equations with the velocity
field given by the solution of the Stokes equations. For simplicity,
we consider only linear problems but note that the framework allows
for implementation of methods for general nonlinear problems. See in
particular~\cite{logg:inprep:updatedelasticity} and
\cite{HofJoh04a,HofJoh04b,HofJoh06}.

\subsection{Static Linear Elasticity}

As a first example, consider the equation of static linear
elasticity~\cite{BreSco94} for the displacement $u = u(x)$ of an
elastic shape~$\Omega \in \R^d$,
\begin{equation} \label{eq:elasticity}
  \begin{array}{rcll}
  - \nabla \cdot \sigma(u) &=& f \quad & \mbox{in } \Omega, \\
  u &=& u_0 \quad & \mbox{on } \Gamma_0 \subset \partial \Omega, \\
  \sigma(u) \hat{n} &=& 0 \quad & \mbox{on } \partial\Omega \setminus \Gamma_0,
  \end{array}
\end{equation}
where $\hat{n}$ denotes a unit vector normal to the boundary
$\partial\Omega$. The stress tensor~$\sigma$ is given by
\begin{equation}
  \sigma(v) = 2\mu \, \epsilon(v) + \lambda \, \mathrm{trace}(\epsilon(v)) I,
\end{equation}
where $I$ is the $d \times d$ identity matrix and where the strain
tensor~$\epsilon$ is given by
\begin{equation}
  \epsilon(v) = \frac{1}{2}
  \left( \nabla v + (\nabla v)^{\top} \right),
\end{equation}
that is, $\epsilon_{ij}(v) =
\frac{1}{2}
(\frac{\partial v_i}{\partial x_j} +
 \frac{\partial v_j}{\partial x_i})$ for $i,j = 1,\ldots,d$.
The Lam\'e constants $\mu$ and $\lambda$ are given by
\begin{equation}
  \mu = \frac{E}{2(1 + \nu)}, \quad
  \lambda = \frac{E \nu}{(1 + \nu)(1 - 2\nu)},
\end{equation}
with $E$ the Young's modulus of elasticity and $\nu$ the Poisson
ratio, see~\cite{ZieTay67}. In the example below, we take $E = 10$ and $\nu = 0.3$.

To obtain the discrete variational problem corresponding
to~(\ref{eq:elasticity}), we multiply with a test function $v$ in a
suitable discrete test space~$\hat{V}_h$ and integrate by parts to
obtain
\begin{equation} \label{eq:elasticity,varproblem}
  \int_{\Omega} \nabla v : \sigma(U) \dx = \int_{\Omega} v \cdot f \dx
  \quad \forall v \in \hat{V}_h.
\end{equation}
The corresponding formulation in the FFC form language is shown in
Table~\ref{tab:form,elasticity} for an approximation with linear
Lagrange elements on tetrahedra. Note that by defining the operators
$\sigma$ and $\epsilon$, it is possible to obtain a very compact notation
that corresponds well with the mathematical notation
of~(\ref{eq:elasticity,varproblem}).

\begin{table}[htbp]
  \begin{center}
    \small
    \begin{code}
 element = FiniteElement("Vector Lagrange", "tetrahedron", 1)

 v = BasisFunction(element)
 U = BasisFunction(element)
 f = Function(element)

 E  = 10.0
 nu = 0.3

 mu = E / (2*(1 + nu))
 lmbda = E*nu / ((1 + nu)*(1 - 2*nu))

 def epsilon(v):
     return 0.5*(grad(v) + transp(grad(v)))

 def sigma(v):
     return 2*mu*epsilon(v) + lmbda*mult(trace(epsilon(v)), Identity(len(v)))

 a = dot(grad(v), sigma(U))*dx
 L = dot(v, f)*dx
    \end{code}
    \normalsize
    \caption{The complete specification of the variational
      problem~(\ref{eq:elasticity,varproblem}) for static linear
      elasticity in the FFC form language.}
    \label{tab:form,elasticity}
  \end{center}
\end{table}

\begin{figure}[htbp]
  \begin{center}
    \includegraphics[width=11cm]{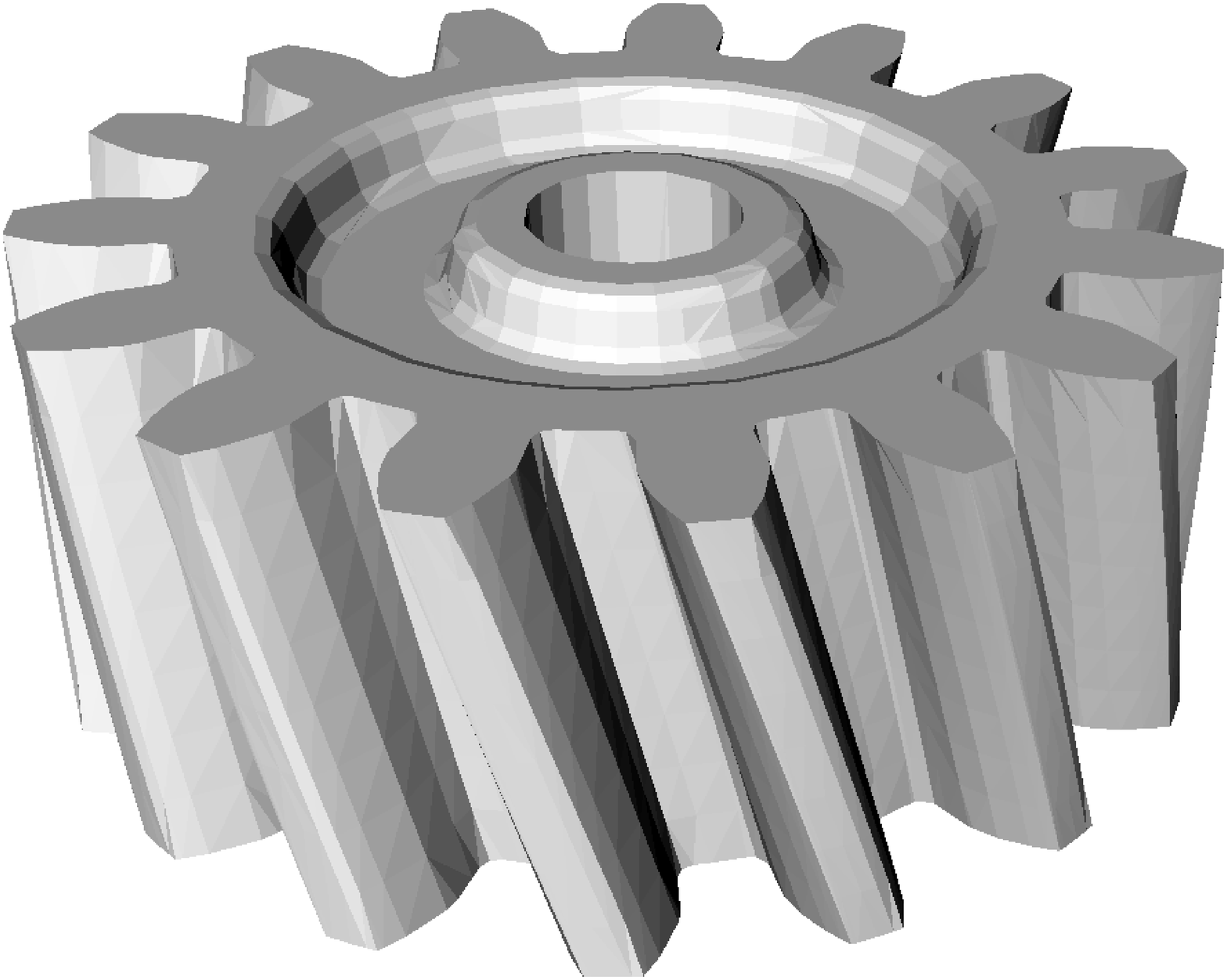} \\
    \includegraphics[width=11cm]{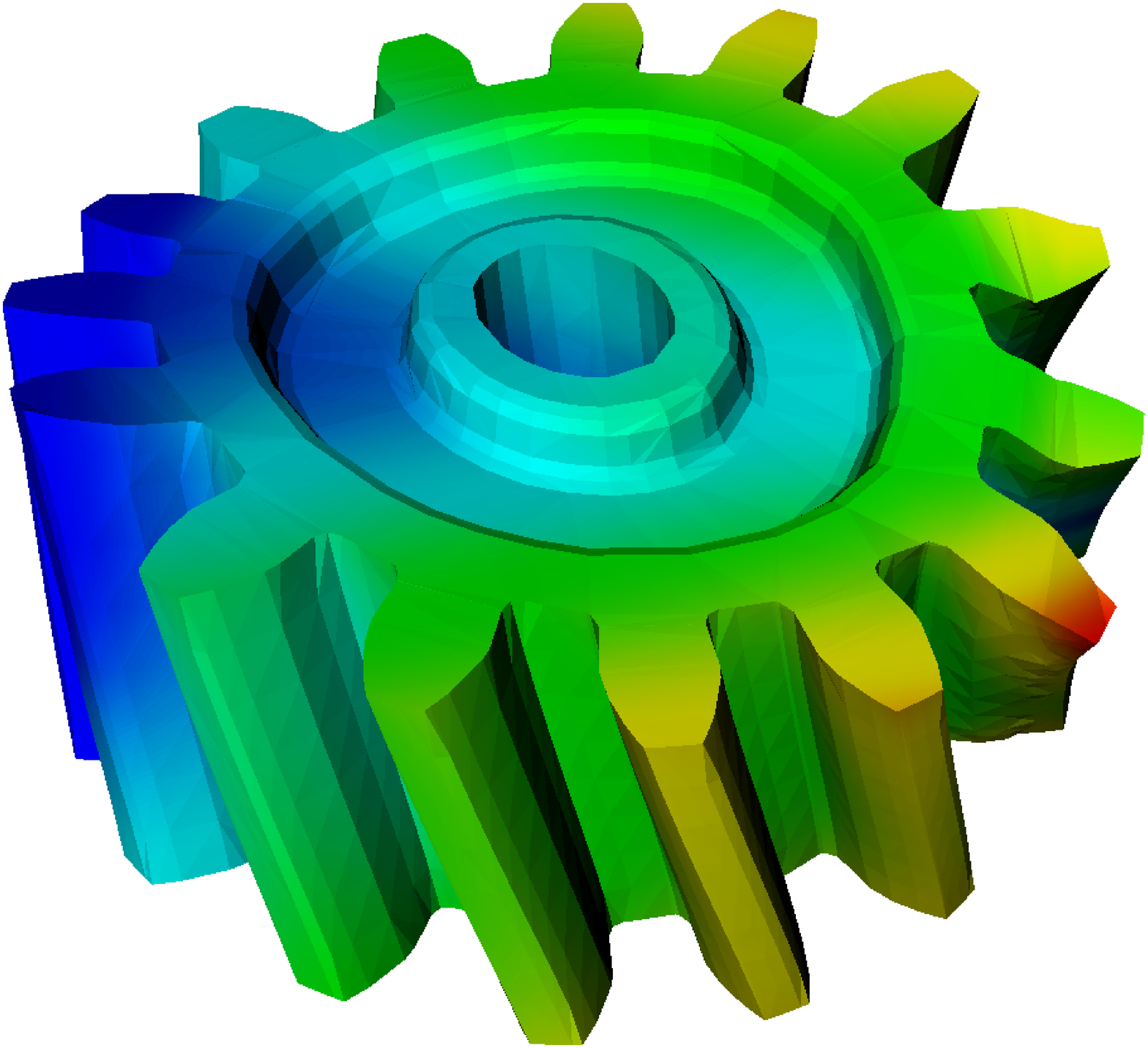}
    \caption{The original domain~$\Omega$ of the gear (above)
      and the twisted gear (below), obtained by displacing $\Omega$
      at each point $x \in \Omega$ by the value of the solution $u$ of
      (\ref{eq:elasticity}) at the point~$x$.}
    \label{fig:gear}
  \end{center}
\end{figure}

Computing the solution of the variational problem for a
domain~$\Omega$ given by a gear, we obtain the solution in
Figure~\ref{fig:gear}. The gear is clamped at two of its ends and
twisted $30$~degrees, as specified by a suitable choice of Dirichlet
boundary conditions on~$\Gamma_0$.

\subsection{The Stokes Equations}

Next, we consider the Stokes equations,
\begin{equation} \label{eq:stokes}
  \begin{array}{rcll}
    - \Delta u + \nabla p &=& f \quad & \mbox{in } \Omega, \\
    \nabla \cdot u &=& 0 \quad & \mbox{in } \Omega, \\
    u &=& u_0 \quad & \mbox{on } \partial \Omega,
  \end{array}
\end{equation}
for the velocity field $u = u(x)$ and the pressure $p = p(x)$ in a
highly viscous medium. By multiplying the two equations with a pair of
test functions $(v, q)$ chosen from a suitable discrete test
space~$\hat{V}_h = \hat{V}_h^u \times \hat{V}_h^p$, we obtain the discrete
variational problem
\begin{equation} \label{eq:stokes,varproblem}
  \int_{\Omega} \nabla v : \nabla U - (\nabla \cdot v) P +
  q \nabla \cdot U \dx =
  \int_{\Omega} v \cdot f \dx
  \quad \forall (v,q) \in \hat{V}_h.
\end{equation}
for the discrete approximate solution~$(U,P) \in V_h = V^u_h \times
V^p_h$. To guarantee the existence of a unique solution of the discrete
variational problem (\ref{eq:stokes,varproblem}), the
discrete function spaces $\hat{V}_h$ and $V_h$ must be chosen
appropriately. The Babu\v{s}ka--Brezzi~\cite{Bab71,Bre74}
inf--sup condition gives a precise condition for the selection of the
approximating spaces.

\subsubsection{Taylor--Hood elements}

One way to fulfill the Babu\v{s}ka--Brezzi condition is to use
different order approximations for the velocity and the pressure, such
as degree~$q$ polynomials for the velocity and degree~$q-1$ for the
pressure, commonly referred to as Taylor--Hood elements, see
\cite{Bof97,BreSco94}. The resulting mixed formulation may be
specified in the FFC form language by defining a Taylor--Hood element
as the direct sum of a degree~$q$ vector-valued Lagrange element and a
degree~$q-1$ scalar Lagrange element, as shown in
Table~\ref{tab:form,stokes,taylorhood}. Figure~\ref{fig:stokes,taylorhood}
shows the velocity field for the flow around a two-dimensional dolphin
computed with a $P_2$--$P_1$ Taylor-Hood approximation.

\begin{table}[htbp]
  \begin{center}
    \begin{code}
 P2 = FiniteElement("Vector Lagrange", "triangle", 2)
 P1 = FiniteElement("Lagrange", "triangle", 1)
 TH = P2 + P1

 (v, q) = BasisFunctions(TH)
 (U, P) = BasisFunctions(TH)

 f = Function(P2)

 a = (dot(grad(v), grad(U)) - div(v)*P + q*div(U))*dx
 L = dot(v, f)*dx
    \end{code}
    \caption{The complete specification of the variational
      problem~(\ref{eq:stokes,varproblem}) for the Stokes
      equations with $P_2$--$P_1$ Taylor--Hood elements.}
    \label{tab:form,stokes,taylorhood}
  \end{center}
\end{table}

\subsubsection{A stabilized equal-order formulation}

Alternatively, the Babu\v{s}ka--Brezzi condition may be circumvented
by an appropriate modification (stabilization) of the variational
problem~(\ref{eq:stokes,varproblem}). In general, an appropriate
modification may be obtained by a Galerkin/least-squares (GLS)
stabilization, that is, by modifying the test function $w = (v, q)$
according to $w \rightarrow w + \delta A w$, where $A$ is
the operator of the differential equation and $\delta = \delta(x)$ is
suitable stabilization parameter.  Here, we a choose simple
pressure-stabilization obtained by modifying the test function $w =
(v, q)$ according to
\begin{equation} \label{eq:pspg}
  (v, q) \rightarrow
  (v, q) + (\delta \nabla q, 0).
\end{equation}
The stabilization~(\ref{eq:pspg}) is sometimes referred to as a
pressure-stabilizing/Petrov-Galerkin (PSPG) method,
see~\cite{HugFra86,FriMat04}. Note that the
stabilization~(\ref{eq:pspg}) may also be viewed as a reduced GLS
stabilization.

We thus obtain the following modified variational problem:
Find $(U,P) \in V_h$ such that
\begin{equation} \label{eq:stokes,varproblem,stabilized}
  \int_{\Omega} \nabla v : \nabla U - (\nabla \cdot v) P +
  q \nabla \cdot U  + \delta \nabla q \cdot \nabla P \dx =
  \int_{\Omega} (v + \delta \nabla q) \cdot f \dx
  \quad \forall (v,q) \in \hat{V}_h.
\end{equation}
Table~\ref{tab:form,stokes,stabilized} shows the stabilized
equal-order method in the FFC form language, with the stabilization
parameter given by
\begin{equation}
  \delta = \beta h^2,
\end{equation}
where $\beta = 0.2$ and $h = h(x)$ is the local mesh size (cell diameter).

\begin{table}[htbp]
  \begin{center}
    \small
    \begin{code}
 vector = FiniteElement("Vector Lagrange", "triangle", 1)
 scalar = FiniteElement("Lagrange", "triangle", 1)
 system = vector + scalar

 (v, q) = BasisFunctions(system)
 (U, P) = BasisFunctions(system)

 f = Function(vector)
 h = Function(scalar)

 d = 0.2*h*h

 a = (dot(grad(v), grad(U)) - div(v)*P + q*div(U) + d*dot(grad(q), grad(P)))*dx
 L = dot(v + mult(d, grad(q)), f)*dx
    \end{code}
    \normalsize
    \caption{The complete specification of the variational
      problem~(\ref{eq:stokes,varproblem,stabilized}) for the Stokes
      equations with an equal-order $P_1$--$P_1$
      stabilized method.}
    \label{tab:form,stokes,stabilized}
  \end{center}
\end{table}

In Figure \ref{fig:stokes,stabilized}, we illustrate the importance of
stabilizing the equal-order method by plotting the solution for the
pressure with and without stabilization. Without stabilization, the
solution oscillates heavily. Note that the scaling is chosen
differently in the two images, with the oscillations scaled down by a
factor two in the unstabilized solution. The situation without
stabilization is thus even worse than what the figure indicates.

\begin{figure}[htbp]
  \begin{center}
    \includegraphics[width=12cm]{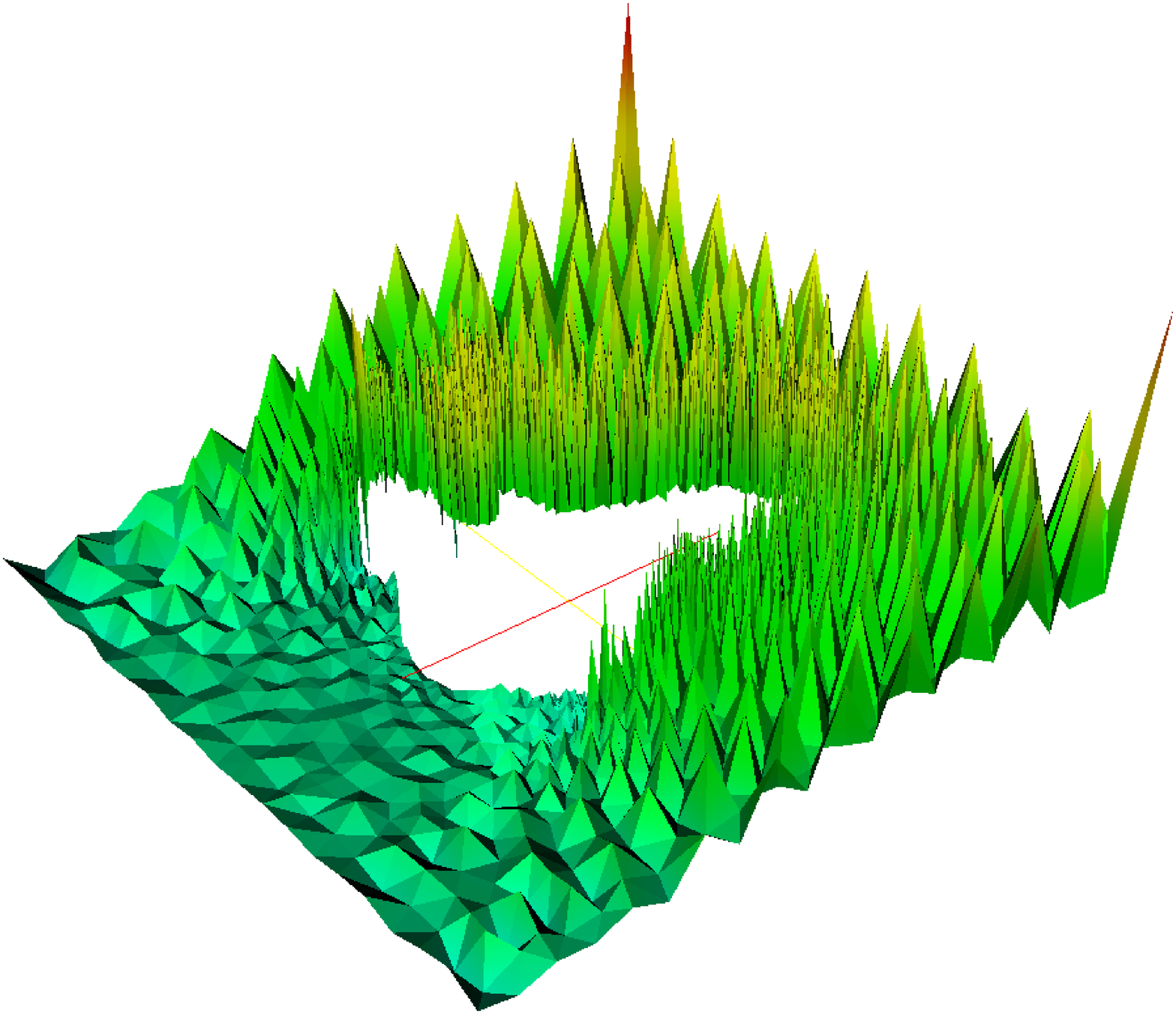} \\
    \includegraphics[width=12cm]{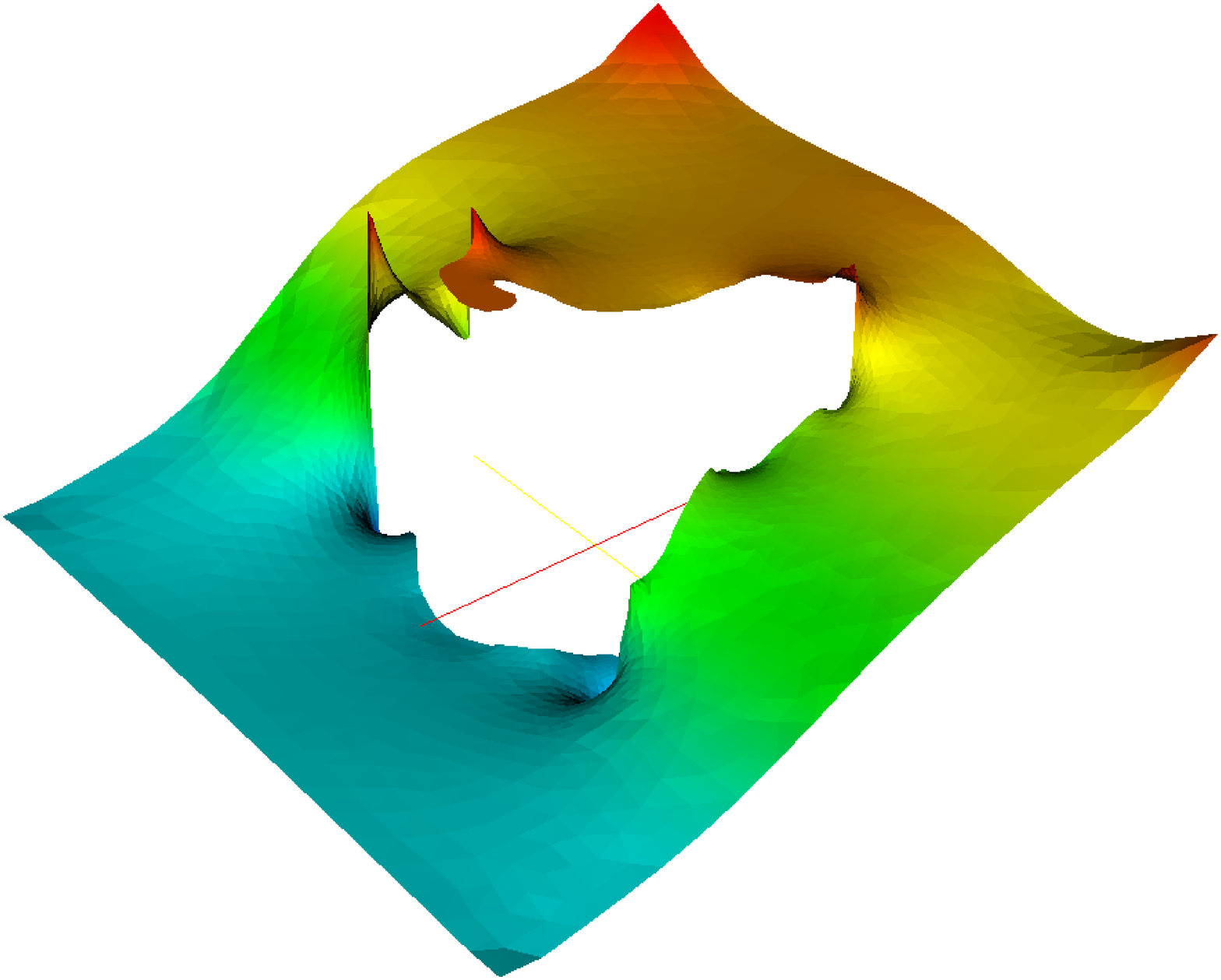}
    \caption{The pressure for the flow around a two-dimensional
      dolphin, obtained by solving the Stokes
      equations~(\ref{eq:stokes}) by an unstabilized
      $P_1$--$P_1$ approximation (above) and a
      stabilized $P_1$--$P_1$ approximation (below).}
    \label{fig:stokes,stabilized}
  \end{center}
\end{figure}

\subsection{Convection--Diffusion}

As a final example, we compute the temperature $u = u(x,t)$
around the dolphin (Figure~\ref{fig:stokes,convectiondiffusion})
from the previous example by solving the
time-dependent convection--diffusion equations,
\begin{equation} \label{eq:convectiondiffusion}
  \begin{array}{rcll}
    \dot{u} + b \cdot \nabla u - \nabla \cdot (c \nabla u) &=& f
    \quad &\mbox{in } \Omega \times (0,T], \\
    u &=& u_{\partial} \quad &\mbox{on } \partial \Omega \times (0,T], \\
    u &=& u_0 \quad &\mbox{at } \Omega \times \{0\},
  \end{array}
\end{equation}
with velocity field $b = b(x)$ obtained by solving the Stokes
equations.

\begin{figure}[htbp]
  \begin{center}
    \includegraphics[width=8.5cm]{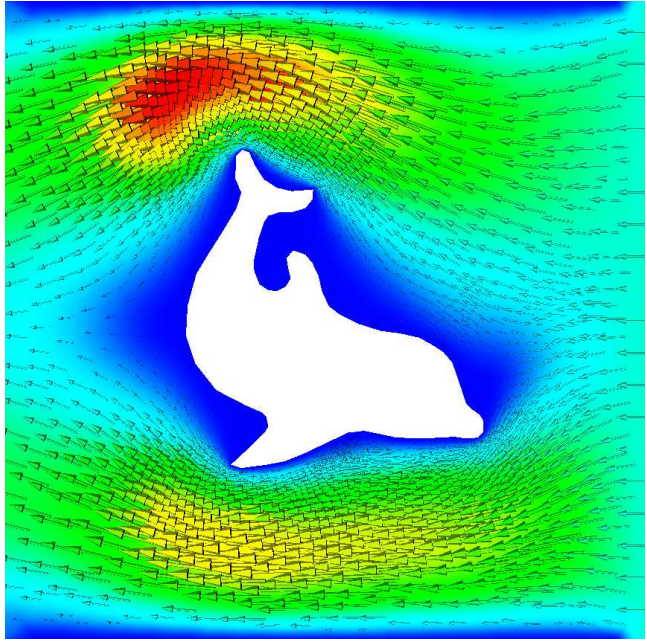}
    \caption{The velocity field for the flow around a two-dimensional
      dolphin, obtained by solving the Stokes equations~(\ref{eq:stokes}) by a
      $P_2$--$P_1$ Taylor-Hood approximation.}
    \label{fig:stokes,taylorhood}
  \end{center}
\end{figure}

\begin{figure}[htbp]
  \begin{center}
    \includegraphics[width=8.5cm]{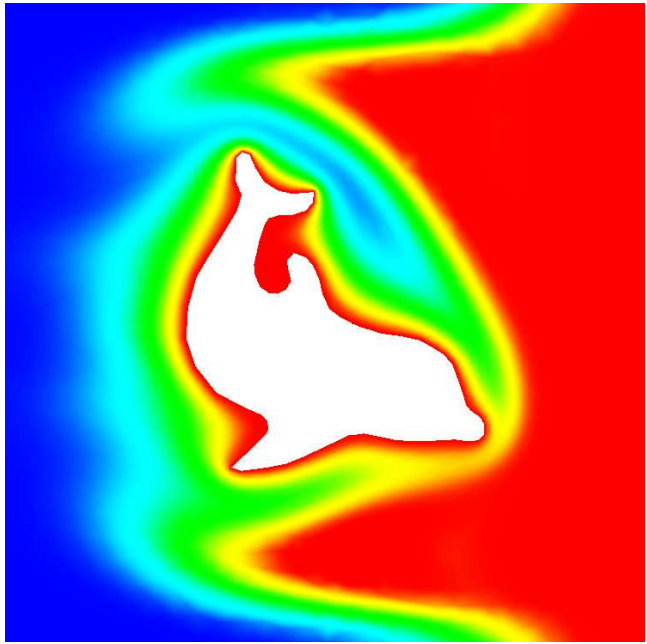}
    \caption{The temperature around a hot dolphin in surrounding cold
      water with a hot inflow, obtained by solving the
      convection--diffusion equation with the velocity field obtained
      from a solution of the Stokes equations with a $P_2$--$P_1$
      Taylor--Hood approximation.}
    \label{fig:stokes,convectiondiffusion}
  \end{center}
\end{figure}

We discretize~(\ref{eq:convectiondiffusion}) with the
$\mathrm{cG}(1)\mathrm{cG}(1)$ method, that is, with continuous piecewise
linear functions in space and time (omitting stabilization for
simplicity). The interval $[0,T]$ is partitioned into a set of
time intervals $0 = t_0 < t_1 < \cdots < t_{n-1} < t_n < \cdots < t_M =
T$ and on each time interval $(t_{n-1}, t_n]$, we pose the variational
problem
\begin{equation}
  \int_{t_{n-1}}^{t_n} \int_{\Omega}
  (v, \dot{U}) + v \, b \cdot \nabla U + c \nabla v \cdot \nabla U \dx
  \dt
  = \int_{t_{n-1}}^{t_n} \int_{\Omega} v \, f \dx \dt
  \quad \forall v \in \hat{V}_h,
\end{equation}
with $\hat{V}_h$ the space of all continuous piecewise linear functions
in space. Note that the $\mathrm{cG}(1)$ method in time uses piecewise constant
test functions, see
\cite{Hul72a,Hul72b,EstFre94,EriEst96,logg:article:01}. As a
consequence, we obtain the following variational problem for $U^n \in
V_h = \hat{V}_h$, the piecewise linear in space solution at time $t = t_n$,
\begin{equation}
  \begin{split}
    \int_{\Omega} v \, \frac{U^n - U^{n-1}}{k_n} +
    v \, b \cdot \nabla (U^n + U^{n-1})/2 +
    c \, \nabla v \cdot \nabla (U^n + U^{n-1})/2 \dx \\
    = \int_{t_{n-1}}^{t_n} \int_{\Omega} v \, f \dx \dt
    \quad \forall v \in \hat{V}_h,
  \end{split}
\end{equation}
where $k_n = t_n - t_{n-1}$ is the size of the time step.
We thus obtain a variational problem of the form
\begin{equation} \label{eq:convectiondiffusion,varproblem}
  a(v, U^n) = L(v) \quad \forall v \in \hat{V}_h,
\end{equation}
where
\begin{equation}
  \begin{split}
    a(v, U^n) &= \int_{\Omega} v \, U^n \dx +
    \frac{k_n}{2} \left(
    v \, b \cdot \nabla U^n +
    c \, \nabla v \cdot \nabla U^n \right) \dx, \\
    L(v) &=
    \int_{\Omega} v \, U^{n-1} \dx -
    \frac{k_n}{2} \left(
    v \, b \cdot \nabla U^{n-1} +
    c \, \nabla v \cdot \nabla U^{n-1} \right) \dx
    + \int_{t_{n-1}}^{t_n} \int_{\Omega} v \, f \dx \dt.
  \end{split}
\end{equation}

The corresponding specification in the FFC form language is presented
in Table~\ref{tab:form,convectiondiffusion}, where for simplicity we
approximate the right-hand side with its value at the right end-point.

\begin{table}[htbp]
  \begin{center}
    \small
    \begin{code}
scalar = FiniteElement("Lagrange", "triangle", 1)
vector = FiniteElement("Vector Lagrange", "triangle", 2)

v  = BasisFunction(scalar)
U1 = BasisFunction(scalar)
U0 = Function(scalar)
b  = Function(vector)
f  = Function(scalar)

c = 0.005
k = 0.05

a = v*U1*dx + 0.5*k*(v*dot(b, grad(U1)) + c*dot(grad(v), grad(U1)))*dx
L = v*U0*dx - 0.5*k*(v*dot(b, grad(U0)) + c*dot(grad(v), grad(U0)))*dx + k*v*f*dx
    \end{code}
    \normalsize
    \caption{The complete specification of the variational
      problem~(\ref{eq:convectiondiffusion,varproblem}) for
      $\mathrm{cG}(1)$ time-stepping of the convection--diffusion
      equation.}
    \label{tab:form,convectiondiffusion}
  \end{center}
\end{table}

%------------------------------------------------------------------------------
\section{OUTLOOK: THE AUTOMATION OF CMM}
\label{sec:outlook}

The automation of the finite element method, as described above,
constitutes an important step towards the Automation of Computational
Mathematical Modeling (ACMM), as outlined in~\cite{logg:thesis:03}.
In this context, the automation of the finite element method amounts
to the \emph{automation of discretization}, that is, the automatic
translation of a given continuous model to a system of discrete
equations. Other key steps include the automation of discrete
solution, the automation of error control, the automation of modeling
and the automation of optimization. We discuss these steps below and
also make some comments concerning automation in general.

\subsection{The Principles of Automation}

An \emph{automatic} system carries out a well-defined task without
intervention from the person or system actuating the automatic
process. The task of the automating system may be formulated as
follows: For given input satisfying a fixed set of conditions (the
\emph{input conditions}), produce output satisfying a given set
of conditions (the \emph{output conditions}).

An automatic process is defined by an \emph{algorithm},
consisting of a sequential list of instructions (like a computer program).
In automated manufacturing, each step of the algorithm operates on and transforms
physical material. Correspondingly, an algorithm for the Automation of CMM
operates on digits and consists of the automated transformation of
digital information.

A key problem of automation is the design of a \emph{feed-back control},
allowing the given output conditions to be satisfied under
variable input and external conditions, ideally at a minimal cost.
Feed-back control is realized through
\emph{measurement}, \emph{evaluation} and \emph{action};
a quantity relating to the given set of conditions to be satisfied by
the output is measured, the measured quantity is evaluated to determine
if the output conditions are satisfied or if an adjustment is
necessary, in which case some action is taken to make the necessary adjustments.
In the context of an algorithm for feed-back control, we refer to the
evaluation of the set of output conditions as the \emph{stopping criterion},
and to the action as the \emph{modification strategy}.

A key step in the automation of a complex process
is \emph{modularization}, that is,
the hierarchical organization of the complex process into components or
sub processes. Each sub process may then itself be automated,
including feed-back control. We may also express this as \emph{abstraction},
that is, the distinction between the properties of a component (its purpose)
and the internal workings of the component (its realization).

Modularization (or abstraction) is central in all engineering and
makes it possible to build complex systems by connecting together
components or subsystems without concern for the internal workings of
each subsystem. The exact partition of a system into components is not
unique. Thus, there are many ways to partition a system into
components. In particular, there are many ways to design a system for
the Automation of Computational Mathematical Modeling.

We thus identify the following basic principles of automation:
algorithms, feed-back control, and modularization.

\subsection{Computational Mathematical Modeling}

In automated manufacturing, the task of the automating system
is to produce a certain product (the output)
from a given piece of material (the input), with the product
satisfying some measure of quality (the output conditions).

\begin{figure}[htbp]
  \begin{center}
    \psfrag{in1}{$A(u) = f$}
    \psfrag{in2}{$\mathrm{TOL} > 0$}
    \psfrag{ut}{\hspace{0.4cm} $U \approx u$}
    \includegraphics[width=12cm]{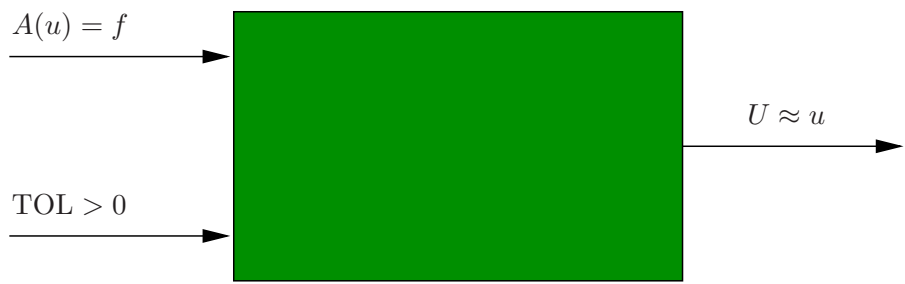}
    \caption{The Automation of Computational Mathematical Modeling.}
    \label{fig:cmm}
  \end{center}
\end{figure}

For the Automation of CMM, the input is a given model of
the form
\begin{equation} \label{eq:model}
  A(u) = f,
\end{equation}
for the solution~$u$ on a given domain $\Omega \times (0,T]$ in
space-time, where~$A$ is a given differential operator and where $f$
is a given source term. The output is a discrete solution $U \approx
u$ satisfying some measure of quality. Typically, the measure of
quality is given in the form of a \emph{tolerance} $\mathrm{TOL} > 0$
for the size of the \emph{error} $e = U - u$ in a suitable norm,
$\|e\| \leq \mathrm{TOL}$, or alternatively, the error in some given
functional~$M$,
\begin{equation} \label{eq:tol}
  \left| M(U) - M(u) \right| \leq \mathrm{TOL}.
\end{equation}
In addition to controlling the quality of the computed solution, one
may also want to determine a parameter that optimizes some given cost
functional depending on the computed solution (optimization). We refer
to the overall process, including optimization, as the Automation of
CMM.

The key problem for the Automation of CMM is thus the design of a
feed-back control for the automatic construction of a discrete
solution $U$, satisfying the output condition (\ref{eq:tol}) at
minimal cost. The design of this feed-back control is based on the
solution of an associated \emph{dual problem}, connecting the size of
the \emph{residual} $R(U) = A(U) - f$ of the computed discrete
solution to the size of the error $e$, and thus to the output
condition (\ref{eq:tol}).

\subsection{An Agenda for the Automation of CMM}

Following our previous discussion on modularization as a basic
principle of automation, we identify the following key steps in
the Automation of CMM:

\begin{itemize}
\item[(i)]
  \emph{the automation of discretization}, that is,
  the automatic translation of a continuous model of the
  form~(\ref{eq:model}) to a system of discrete equations;
\item[(ii)]
  \emph{the automation of discrete solution}, that is, the automatic
  solution of the system of discrete equations obtained from the
  automatic discretization of~(\ref{eq:model}),
\item[(iii)]
  \emph{the automation of error control}, that is,
  the automatic selection of an appropriate resolution of the discrete
  model to produce a discrete solution satisfying the given accuracy
  requirement with minimal work;
\item[(iv)]
  \emph{the automation of modeling}, that is,
  the automatic selection of the model~(\ref{eq:model}),
  either by constructing a model from a given set of data, or by
  constructing from a given model a~\emph{reduced model} for the
  variation of the solution on resolvable scales;
\item[(v)]
  \emph{the automation of optimization}, that is,
  the automatic selection of a parameter in the
  model~(\ref{eq:model}) to optimize a given goal functional.
\end{itemize}

In Figure~\ref{fig:automation,1234}, we demonstrate how (i)--(iv)
connect to solve the overall task of the Automation of CMM (excluding
optimization) in accordance with Figure~\ref{fig:cmm}.
We discuss (i)--(v) in some detail below. In all cases, feed-back control,
or \emph{adaptivity}, plays a key role.

\begin{figure}[htbp]
  \begin{center}
    \psfrag{eq}{$A(u) = f$}
    \psfrag{TOL}{$\mathrm{TOL}$}
    \psfrag{U}{$U \approx u$}
    \psfrag{tol}{$\mathrm{tol} > 0$}
    \psfrag{eq2}{$A(\bar{u}) = \bar{f} + \bar{g}(u)$}
    \psfrag{V}{$(\hat{V}_h, V_h)$}
    \psfrag{F}{$F(x) = 0$}
    \psfrag{X}{$X \approx x$}
    \includegraphics[width=14cm]{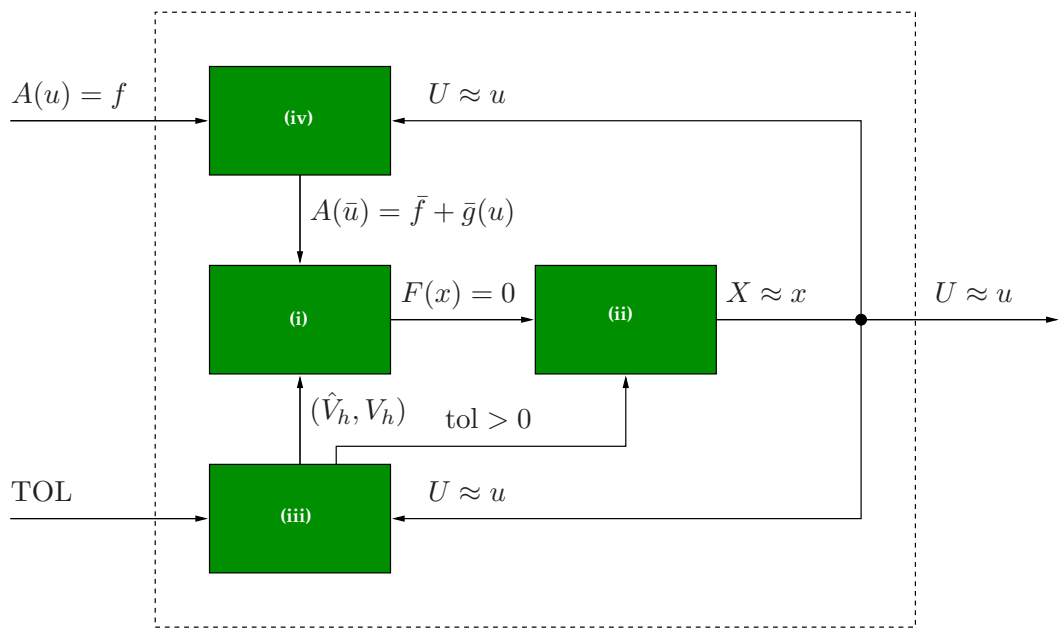}
    \caption{A modularized view of the Automation of Computational
    Mathematical Modeling.}
    \label{fig:automation,1234}
  \end{center}
\end{figure}

\subsection{The Automation of Discretization}

The automation of discretization amounts to automatically generating
a system of discrete equations for the degrees of freedom of a
discrete solution~$U$ approximating the solution~$u$ of the given
model~(\ref{eq:model}), or alternatively, the solution~$u \in V$ of a
corresponding variational problem
\begin{equation} \label{eq:varproblem,again}
  a(u; v) = L(v) \quad \forall v \in \hat{V},
\end{equation}
where as before $a : V \times \hat{V} \rightarrow \R$ is a semilinear
form which is linear in its second argument and $L : \hat{V} \rightarrow \R$
is a linear form. As we saw in Section~\ref{sec:fem}, this process may
be automated by the finite element method, by replacing the function
spaces~$(\hat{V}, V)$ with a suitable pair $(\hat{V}_h, V_h)$ of discrete
function spaces, and an approach to its automation was discussed in
Sections~\ref{sec:automation,tabulation}--\ref{sec:automation,assembly}.
As we shall discuss further below, the pair of discrete function spaces may be
automatically chosen by feed-back control to compute the discrete
solution~$U$ both reliably and efficiently.

\subsection{The Automation of Discrete Solution}

Depending on the model~(\ref{eq:model}) and the method used to
automatically discretize the model, the resulting system of discrete
equations may require more or less work to solve. Typically, the
discrete system is solved by some iterative method such as the
conjugate gradient method (CG) or GMRES, in combination with an
appropriate choice of preconditioner, see for
example~\cite{Saa03,Dem97}.

The resolution of the discretization of~(\ref{eq:model}) may
be chosen automatically by feed-back control from the computed
solution, with the target of minimizing the computational work while
satisfying a given accuracy requirement. As a consequence, see for
example~\cite{logg:article:01}, one obtains an accuracy requirement on
the solution of the system of discrete equations. Thus, the system of
discrete equations does not need to be solved to within machine
precision, but only to within some \emph{discrete
tolerance}~$\mathrm{tol > 0}$ for some error in a functional of the
solution of the discrete system. We shall not pursue this question
further here, but remark that the feed-back control from the computed
solution to the iterative algorithm for the solution of the system of
discrete equations is often weak, and the problem of designing
efficient \emph{adaptive} iterative algorithms for the system of
discrete equations remains open.

\subsection{The Automation of Error Control}

As stated above, the overall task is to produce a solution
of~(\ref{eq:model}) that satisfies a given accuracy requirement with
minimal work. This includes an aspect of \emph{reliability}, that is,
the error in an output quantity of interest depending on the computed
solution should be less than a given tolerance, and an aspect of
\emph{efficiency}, that is, the solution should be computed with
minimal work. Ideally, an algorithm for the solution of
(\ref{eq:model}) should thus have the following properties: Given a
tolerance $\mathrm{TOL} > 0$ and a functional $M$, the algorithm shall
produce a discrete solution $U$ approximating the exact solution $u$
of (\ref{eq:model}), such that
\begin{itemize}
  \item[(A)]
    $\vert M(U) - M(u) \vert \leq \mathrm{TOL}$;
  \item[(B)]
    the computational cost of obtaining the approximation $U$
    is minimal.
\end{itemize}
Conditions (A) and (B) can be satisfied by an \emph{adaptive algorithm},
with the construction of the \emph{discrete
representation}~$(\hat{V}_h,V_h)$ based on feed-back from the computed
solution.

An adaptive algorithm typically involves a stopping criterion,
indicating that the size of the error is less than the given
tolerance, and a modification strategy to be applied if the stopping
criterion is not satisfied. Often, the stopping criterion and the
modification strategy are based on an \emph{a posteriori} error
estimate $E \geq \vert M(U) - M(u) \vert$, estimating the error
in terms of the residual $R(U) = A(U) - f$ and the solution
$\varphi$ of a dual problem connecting to the stability of
(\ref{eq:model}).

\subsubsection{The dual problem}

The dual problem of~(\ref{eq:model}) for the given output
functional~$M$ is given by
\begin{equation} \label{eq:dual}
  \overline{A'}^{*} \varphi = \psi,
\end{equation}
on $\Omega \times [0,T)$, where $\overline{A'}^*$ denotes the
adjoint\footnote{The adjoint is defined by $(Av, w) = (v, A^*w)$ for
  all $v, w \in V$ such that $v = w = 0$ at $t = 0$ and $t = T$.}
of the Fr\'echet derivative~$A'$ of $A$ evaluated at a suitable mean
value of the exact solution~$u$ and the computed solution~$U$,
\begin{equation}
  \overline{A'} = \int_0^1 A' \left( sU + (1-s)u \right) \ds,
\end{equation}
and where $\psi$ is the Riesz representer of a similar mean value of
the Fr\'echet derivative~$M'$ of $M$,
\begin{equation}
  (v, \psi) = \overline{M'} v \quad \forall v \in V.
\end{equation}
By the dual problem~(\ref{eq:dual}), we directly obtain the error
representation
\begin{equation} \label{eq:errorrepresentation}
  \begin{split}
    M(U) - M(u)
    &= \overline{M'} (U - u) = (U - u, \psi) = (U - u, \overline{A'}^* \varphi)
    = (\overline{A'} (U - u), \varphi) \\
    &= (A(U) - A(u), \varphi) = (A(U) - f, \varphi) = (R(U), \varphi).
  \end{split}
\end{equation}
Noting now that if the solution $U$ is computed by a Galerkin method and
thus $(R(U), v) = 0$ for any $v \in \hat{V}_h$, we obtain
\begin{equation} \label{eq:errorrepresentation,interpolant}
  M(U) - M(u) = (R(U), \varphi - \pi_h \varphi),
\end{equation}
where $\pi_h \varphi$ is a suitable approximation of~$\varphi$
in~$\hat{V}_h$. One may now proceed to estimate the error~$M(U) -
M(u)$ in various ways, either by estimating the interpolation error
$\pi_h \varphi - \varphi$ or by directly evaluating the quantity
$(R(U), \varphi - \pi_h \varphi)$. The residual $R(U)$ and the dual
solution~$\varphi$ give precise information about the influence of
the discrete representation~$(\hat{V}_h, V_h)$ on the size of the
error, which can be used in an adaptive feed-back control to choose
a suitable discrete representation for the given output quantity~$M$ of
interest and the given tolerance~$\mathrm{TOL}$ for the error,
see~\cite{EriEst95,BecRan01,EstLar00,logg:article:01}.

\subsubsection{The weak dual problem}

We may estimate the error similarly for the variational
problem~(\ref{eq:varproblem,again}) by considering the
following weak (variational) dual problem: Find $\varphi \in
\hat{V}$ such that
\begin{equation} \label{eq:varproblem,dual}
  \overline{a'}^*(U, u;\, v, \varphi) = \overline{M'}(U, u; v) \quad \forall v \in V,
\end{equation}
where~$\overline{a'}^*$ denotes the adjoint of the bilinear
form~$\overline{a'}$, given as above by an appropriate mean value of
the Fr\'echet derivative of the semilinear form~$a$.
We now obtain the error representation
\begin{equation}
  \begin{split}
    M(U) - M(u)
    &= \overline{M'}(U, u;\, U - u)
     = \overline{a'}^*(U, u;\, U - u, \varphi)
     = \overline{a'}(U, u;\, \varphi, U - u) \\
    &= a(U;\, \varphi) - a(u;\, \varphi)
     = a(U;\, \varphi) - L(\varphi).
  \end{split}
\end{equation}
As before, we use the Galerkin orthogonality to subtract
$a(U; \pi\varphi) - L(\pi\varphi) = 0$ for some $\pi_h \varphi \in \hat{V}_h
\subset \hat{V}$ and obtain
\begin{equation}
  M(U) - M(u) = a(U; \varphi - \pi \varphi) - L(\varphi - \pi \varphi).
\end{equation}

To automate the process of error control, we thus need to
automatically generate and solve the dual problem~(\ref{eq:dual})
or~(\ref{eq:varproblem,dual}) from a given primal
problem~(\ref{eq:model}) or~(\ref{eq:varproblem,again}). We
investigate this question further in~\cite{logg:inprep:errorautomation}.

\subsection{The Automation of Modeling}

The automation of modeling concerns both the problem of finding the
parameters describing the model (\ref{eq:model}) from a given set of
data (inverse modeling), and the automatic construction of a reduced
model for the variation of the solution on resolvable scales (model
reduction). We here discuss briefly the \emph{automation of model
reduction}.

In situations where the solution $u$ of (\ref{eq:model}) varies on
scales of different magnitudes, and these scales are not localized in
space and time, computation of the solution may be very expensive,
even with an adaptive method. To make computation feasible, one may
instead seek to compute an \emph{average} $\bar{u}$ of the solution
$u$ of (\ref{eq:model}) on resolvable scales. Typical examples
include meteorological models for weather prediction, with fast time
scales on the range of seconds and slow time scales on the range of
years, or protein folding represented by a molecular dynamics model,
with fast time scales on the range of femtoseconds and slow time
scales on the range of microseconds.

Model reduction typically involves extrapolation from resolvable
scales, or the construction of a large-scale model from local
resolution of fine scales in time and space. In both cases, a
large-scale model
\begin{equation}
  A(\bar{u}) = \bar{f} + \bar{g}(u),
\end{equation}
for the average $\bar{u}$ is constructed from the
given model (\ref{eq:model}) with a suitable \emph{modeling term}
$\bar{g}(u) \approx A(\bar{u}) - \bar{A}(u)$.

Replacing a given model with a computable reduced model by taking
averages in space and time is sometimes referred to as \emph{subgrid
modeling}. Subgrid modeling has received much attention in recent
years, in particular for the incompressible Navier--Stokes equations,
where the subgrid modeling problem takes the form
of determining the \emph{Reynolds stresses} corresponding to
$\bar{g}$. Many subgrid models have been proposed for the averaged
Navier--Stokes equations, but no clear answer has been given.
Alternatively, the subgrid model may take the form of a least-squares
stabilization, as suggested in \cite{HofJoh04a,HofJoh04b,HofJoh06}.
In either case, the validity of a proposed subgrid model may be
verified computationally by solving an appropriate dual problem and
computing the relevant residuals to obtain an error estimate for the
modeling error, see~\cite{logg:article:06}.

\subsection{The Automation of Optimization}

The \emph{automation of optimization} relies on the automation of
(i)--(iv), with the solution of the primal problem (\ref{eq:model})
and an associated dual problem being the key steps in the minimization
of a given cost functional. In particular, the automation of
optimization relies on the automatic generation of the dual problem.

The optimization of a given cost functional~$\mathcal{J} =
\mathcal{J}(u,p)$, subject to the constraint~(\ref{eq:model}), with
$p$ a function (the control variables) to be determined, can be
formulated as the problem of finding a stationary point of the
associated Lagrangian,
\begin{equation}
  L(u,p,\varphi) = \mathcal{J}(u,p) + (A(u,p) - f(p), \varphi),
\end{equation}
which takes the form of a system of differential equations,
involving the primal and dual problems, as well as an equation
expressing stationarity with respect to the control variables~$p$,
\begin{equation}
  \begin{split}
    A(u,p) &= f(p), \\
    (A')^*(u,p) \varphi &= - \partial \mathcal{J} / \partial u, \\
    \partial \mathcal{J} / \partial p &=
    (\partial f / \partial p)^* \varphi - (\partial A / \partial p)^* \varphi. \\
  \end{split}
\end{equation}
It follows that the optimization problem may be solved by the solution
of a system of differential equations. Note that the first equation is
the given model~(\ref{eq:model}), the second equation is the dual
problem and the third equation gives a direction for the update of the
control variables. The automation of optimization thus relies on the
automated solution of both the primal problem (\ref{eq:model}) and the
dual problem (\ref{eq:dual}), including the automatic generation of
the dual problem.

%------------------------------------------------------------------------------
\section{CONCLUDING REMARKS}
\label{sec:conclusion}

With the FEniCS project~\cite{logg:www:03}, we have the beginnings of
a working system automating (in part) the finite element method, which
is the first step towards the Automation of Computational Mathematical
Modeling, as outlined in~\cite{logg:thesis:03}. As part of this work,
a number of key components, FIAT, FFC and DOLFIN, have been developed.
These components provide reference implementations of the algorithms
discussed in
Sections~\ref{sec:fem}--\ref{sec:automation,assembly}.

As the current toolset, focused mainly on an automation of the finite
element method (the automation of discretization), is becoming more
mature, important new areas of research and development emerge,
including the remaining key steps towards the Automation of CMM. In
particular, we plan to explore the possibility of automatically
generating dual problems and error estimates in an effort to automate
error control.

%------------------------------------------------------------------------------
\section*{ACKNOWLEDGMENT}

This work is in large parts based on the joint efforts of the people
behind the FEniCS project, in particular Todd Dupont, Johan Hoffman,
Johan Jansson, Claes Johnson, Robert C. Kirby, Matthew G. Knepley, Mats G. Larson,
L. Ridgway Scott, Andy Terrel and Garth N. Wells.

%------------------------------------------------------------------------------
\newpage
\bibliographystyle{siam}
\bibliography{bibliography}

%------------------------------------------------------------------------------
\newpage
\section*{NOTATION}

\linespread{1.2}
\begin{table}[H]
  \begin{tabular}{ccl}
    $A$ &--&
    the differential operator of the model $A(u) = f$ \\
    $A$ &--&
    the \emph{global tensor} with entries $\{A_i\}_{i\in\mathcal{I}}$ \\
    $A^0$ &--&
    the \emph{reference tensor} with entries $\{A^0_{i\alpha}\}_{i\in\mathcal{I}_K,\alpha\in\mathcal{A}}$ \\
    $\bar{A}^0$ &--&
    the matrix representation of the (flattened) reference tensor $A^0$ \\
    $A^K$ &--&
    the \emph{element tensor} with entries $\{A^K_i\}_{i\in\mathcal{I}_K}$ \\
    $a$ &--&
    the semilinear, multilinear or bilinear form \\
    $a_K$ &--&
    the local contribution to a multilinear form $a$ from $K$ \\
    $a^K$ &--&
    the vector representation of the (flattened) element tensor $A^K$ \\
    $\mathcal{A}$ &--&
    the set of \emph{secondary indices} \\
    $\mathcal{B}$ &--&
    the set of \emph{auxiliary indices} \\
    $e$ &--&
    the \emph{error}, $e = U - u$ \\
    $F_K$ &--&
    the mapping from $K_0$ to $K$ \\
    $G_K$ &--&
    the \emph{geometry tensor} with entries $\{G_K^{\alpha}\}_{\alpha\in\mathcal{A}}$ \\
    $g_K$ &--&
    the vector representation of the (flattened) geometry tensor $G_K$ \\
    $\mathcal{I}$ &--&
    the set $\prod_{j=1}^r [1,N^j]$ of indices for the global tensor $A$ \\
    $\mathcal{I}_K$ &--&
    the set $\prod_{j=1}^r [1,n_K^j]$ of indices for the element tensor
    $A^K$ (\emph{primary indices}) \\
    $\iota_K$ &--&
    the \emph{local-to-global mapping} from $\mathcal{N}_K$ to $\mathcal{N}$ \\
    $\hat{\iota}_K$ &--&
    the local-to-global mapping from $\hat{\mathcal{N}}_K$ to $\hat{\mathcal{N}}$ \\
    $\iota_K^j$ &--&
    the local-to-global mapping from $\mathcal{N}_K^j$ to $\mathcal{N}^j$ \\
    $K$ &--&
    a \emph{cell} in the mesh~$\mathcal{T}$ \\
    $K_0$ &--&
    the \emph{reference cell} \\
    $L$ &--&
    the linear form (functional) on $\hat{V}$ or $\hat{V}_h$ \\
    $m$ &--&
    the number of discrete function spaces used in the definition of~$a$ \\
    $N$ &--&
    the dimension of $\hat{V}_h$ and $V_h$ \\
    $N^j$ &--&
    the dimension $V_h^j$ \\
    $N_q$ &--&
    the number of quadrature points on a cell \\
    $n_0$ &--&
    the dimension of $\mathcal{P}_0$ \\
    $n_K$ &--&
    the dimension of $\mathcal{P}_K$ \\
    $\hat{n}_K$ &--&
    the dimension of $\hat{\mathcal{P}}_K$ \\
    $n_K^j$ &--&
    the dimension of $\mathcal{P}_K^j$ \\
    $\mathcal{N}$ &--&
    the set of global nodes on $V_h$ \\
    $\hat{\mathcal{N}}$ &--&
    the set of global nodes on $\hat{V}_h$ \\
    $\mathcal{N}^j$ &--&
    the set of global nodes on $V_h^j$ \\
    $\mathcal{N}_0$ &--&
    the set of local nodes on $\mathcal{P}_0$ \\
    $\mathcal{N}_K$ &--&
    the set of local nodes on $\mathcal{P}_K$ \\
    $\hat{\mathcal{N}}_K$ &--&
    the set of local nodes on $\hat{\mathcal{P}}_K$ \\
    $\mathcal{N}_K^j$ &--&
    the set of local nodes on ${\mathcal{P}_K^j}$ \\
  \end{tabular}
\end{table}
\linespread{1.0}

\linespread{1.2}
\begin{table}[H]
  \begin{tabular}{ccl}
    $\nu^0_i$ &--&
    a \emph{node} on $\mathcal{P}_0$ \\
    $\nu^K_i$ &--&
    a node on $\mathcal{P}_K$ \\
    $\hat{\nu}^K_i$ &--&
    a node on $\hat{\mathcal{P}}_K$ \\
    $\nu^{K,j}_i$ &--&
    a node on $\mathcal{P}_K^j$ \\
    $\mathcal{P}_0$ &--&
    the function space on $K_0$ for $V_h$ \\
    $\hat{\mathcal{P}}_0$ &--&
    the function space on $K_0$ for $\hat{V}_h$ \\
    $\mathcal{P}_0^j$ &--&
    the function space on $K_0$ for $V_h^j$ \\
    $\mathcal{P}_K$ &--&
    the local function space on $K$ for $V_h$ \\
    $\hat{\mathcal{P}}_K$ &--&
    the local function space on $K$ for $\hat{V}_h$ \\
    $\mathcal{P}_K^j$ &--&
    the local function space on $K$ for $V_h^j$ \\
    $P_q(K)$ &--&
    the space of polynomials of degree $\leq q$ on $K$ \\
    $\overline{\mathcal{P}}_K$ &--&
    the local function space on $K$ generated by $\{\mathcal{P}_K^j\}_{j=1}^m$ \\
    $R$ &--&
    the \emph{residual}, $R(U) = A(U) - f$ \\
    $r$ &--&
    the arity of the multilinear form $a$ (the rank of $A$ and $A^K$) \\
    $U$ &--&
    the discrete approximate solution, $U \approx u$ \\
    $(U_i)$ &--&
    the vector of expansion coefficients for $U = \sum_{i=1}^N U_i \phi_i$ \\
    $u$ &--&
    the exact solution of the given model $A(u) = f$ \\
    $V$ &--&
    the space of trial functions on $\Omega$ (the trial space) \\
    $\hat{V}$ &--&
    the space of test functions on $\Omega$ (the test space) \\
    $V_h$ &--&
    the space of discrete trial functions on $\Omega$ (the discrete trial space) \\
    $\hat{V}_h$ &--&
    the space of discrete test functions on $\Omega$ (the discrete test space) \\
    $V_h^j$ &--&
    a discrete function space on $\Omega$ \\
    $|V|$ &--&
    the dimension of a vector space $V$ \\
    $\Phi_i$ &--&
    a basis function in $\mathcal{P}_0$ \\
    $\hat{\Phi}_i$ &--&
    a basis function in $\hat{\mathcal{P}}_0$ \\
    $\Phi_i^j$ &--&
    a basis function in $\mathcal{P}_0^j$ \\
    $\phi_i$ &--&
    a basis function in $V_h$ \\
    $\hat{\phi}_i$ &--&
    a basis function in $\hat{V}_h$ \\
    $\phi_i^j$ &--&
    a basis function in $V_h^j$ \\
    $\phi_i^K$ &--&
    a basis function in $\mathcal{P}_K$ \\
    $\hat{\phi}_i^K$ &--&
    a basis function in $\hat{\mathcal{P}}_K$ \\
    $\phi_i^{K,j}$ &--&
    a basis function in $\mathcal{P}_K^j$ \\
    $\varphi$ &--&
    the \emph{dual solution} \\
    $\mathcal{T}$ &--&
    the \emph{mesh} \\
    $\Omega$ &--&
    a bounded domain in~$\R^d$ \\
  \end{tabular}
\end{table}
\linespread{1.0}

\end{document}